\documentclass[12pt,a4paper]{article}
\usepackage{latexsym,epsfig,fullpage,amsfonts,amssymb,amsmath,amscd,epic,graphics,rotating,theorem,
yhmath,pxfonts}
\usepackage[all]{xy}
\usepackage[usenames]{color}
\theoremstyle{plain}
\newtheorem{lemma}{Lemma}[section]
\newtheorem{proposition}[lemma]{Proposition}
\newtheorem{remark}[lemma]{Remark}
\newtheorem{example}[lemma]{Example}
\newtheorem{theorem}[lemma]{Theorem}
\newtheorem{definition}[lemma]{Definition}

{\theorembodyfont{\rmfamily}  \font\ncsc=cmcsc10
 \font\ntt=cmtt12

\begin{document}
\newcommand{\pperp}{\hbox{$\perp\hskip-6pt\perp$}}
\newcommand{\ssim}{\hbox{$\hskip-2pt\sim$}}
\newcommand{\aleq}{{\ \stackrel{3}{\le}\ }}
\newcommand{\ageq}{{\ \stackrel{3}{\ge}\ }}
\newcommand{\aeq}{{\ \stackrel{3}{=}\ }}
\newcommand{\bleq}{{\ \stackrel{n}{\le}\ }}
\newcommand{\bgeq}{{\ \stackrel{n}{\ge}\ }}
\newcommand{\beq}{{\ \stackrel{n}{=}\ }}
\newcommand{\cleq}{{\ \stackrel{2}{\le}\ }}
\newcommand{\cgeq}{{\ \stackrel{2}{\ge}\ }}
\newcommand{\ceq}{{\ \stackrel{2}{=}\ }}
\newcommand{\N}{{\mathbb N}}
\newcommand{\A}{{\mathbb A}}
\newcommand{\K}{{\mathbb K}}
\newcommand{\Z}{{\mathbb Z}}\newcommand{\F}{{\mathbf F}}
\newcommand{\R}{{\mathbb R}}
\newcommand{\C}{{\mathbb C}}
\newcommand{\Q}{{\mathbb Q}}
\newcommand{\PP}{{\mathbb P}}
\newcommand{\mnote}{\marginpar}\newcommand{\red}{{\operatorname{red}}}
\newcommand{\Pic}{{\operatorname{Pic}}}\newcommand{\Sym}{{\operatorname{Sym}}}
\newcommand{\oeps}{{\overline\eps}}\newcommand{\Div}{{\operatorname{Div}}}
\newcommand{\oDel}{{\widetilde\Del}}
\newcommand{\real}{{\operatorname{Re}}}\newcommand{\Aut}{{\operatorname{Aut}}}
\newcommand{\conv}{{\operatorname{conv}}}\newcommand{\BG}{{\operatorname{BG}}}
\newcommand{\Span}{{\operatorname{Span}}}\newcommand{\GS}{{\operatorname{RB}}}
\newcommand{\Ker}{{\operatorname{Ker}}}\newcommand{\GGS}{{\operatorname{GS}}}
\newcommand{\Ann}{{\operatorname{Ann}}}\newcommand{\mt}{{\operatorname{wt}}}
\newcommand{\Fix}{{\operatorname{Fix}}}\newcommand{\Ima}{{\operatorname{Im}}}
\newcommand{\sign}{{\operatorname{sign}}}
\newcommand{\Tors}{{\operatorname{Tors}}}\newcommand{\lab}{{\operatorname{lab}}}
\newcommand{\Card}{{\operatorname{Card}}}\newcommand{\Def}{{\operatorname{Def}}}
\newcommand{\alg}{{\operatorname{alg}}}\newcommand{\ord}{{\operatorname{ord}}}
\newcommand{\oi}{{\overline i}}
\newcommand{\oj}{{\overline j}}
\newcommand{\ob}{{\overline b}}
\newcommand{\os}{{\overline s}}
\newcommand{\oa}{{\overline a}}
\newcommand{\oy}{{\overline y}}
\newcommand{\ow}{{\overline w}}
\newcommand{\ot}{{\overline t}}
\newcommand{\oz}{{\overline z}}
\newcommand{\eps}{{\varepsilon}}
\newcommand{\proofend}{\hfill$\Box$\bigskip}
\newcommand{\Int}{{\operatorname{Int}}}
\newcommand{\pr}{{\operatorname{pr}}}
\newcommand{\Hom}{{\operatorname{Hom}}}
\newcommand{\Ev}{{\operatorname{Ev}}}
\newcommand{\im}{{\operatorname{Im}}}\newcommand{\br}{{\operatorname{br}}}
\newcommand{\sk}{{\operatorname{sk}}}\newcommand{\DP}{{\operatorname{DP}}}
\newcommand{\const}{{\operatorname{const}}}
\newcommand{\Sing}{{\operatorname{Sing}}\hskip0.06cm}
\newcommand{\conj}{{\operatorname{Conj}}}
\newcommand{\Cl}{{\operatorname{Cl}}}
\newcommand{\Crit}{{\operatorname{Crit}}}
\newcommand{\Ch}{{\operatorname{Ch}}}
\newcommand{\discr}{{\operatorname{discr}}}
\newcommand{\Tor}{{\operatorname{Tor}}}
\newcommand{\Conj}{{\operatorname{Conj}}}
\newcommand{\vol}{{\operatorname{vol}}}
\newcommand{\defect}{{\operatorname{def}}}
\newcommand{\codim}{{\operatorname{codim}}}
\newcommand{\tmu}{{\C\mu}}
\newcommand{\ov}{{\overline v}}
\newcommand{\ox}{{\overline{x}}}
\newcommand{\bw}{{\boldsymbol w}}\newcommand{\bn}{{\boldsymbol n}}
\newcommand{\bx}{{\boldsymbol x}}
\newcommand{\bd}{{\boldsymbol d}}
\newcommand{\bz}{{\boldsymbol z}}\newcommand{\bp}{{\boldsymbol p}}
\newcommand{\tet}{{\theta}}
\newcommand{\Del}{{\Delta}}
\newcommand{\bet}{{\beta}}
\newcommand{\kap}{{\kappa}}
\newcommand{\del}{{\delta}}
\newcommand{\sig}{{\sigma}}
\newcommand{\alp}{{\alpha}}
\newcommand{\Sig}{{\Sigma}}
\newcommand{\Gam}{{\Gamma}}
\newcommand{\gam}{{\gamma}}\newcommand{\idim}{{\operatorname{idim}}}
\newcommand{\Lam}{{\Lambda}}
\newcommand{\lam}{{\lambda}}
\newcommand{\SC}{{SC}}
\newcommand{\MC}{{MC}}
\newcommand{\nek}{{,...,}}
\newcommand{\cim}{{c_{\mbox{\rm im}}}}
\newcommand{\clM}{\tilde{M}}
\newcommand{\clV}{\bar{V}}

\title{Refined elliptic tropical enumerative invariants}
\author{Franziska Schroeter
\and Eugenii Shustin}
\date{}
\maketitle
\begin{abstract}
We suggest a new refined (i.e., depending on a parameter)
tropical enumerative invariant of toric surfaces. This is the first
known enumerative invariant that counts tropical curves of positive genus with marked
vertices. Our invariant extends the refined rational broccoli invariant
invented by L. G\"ottsche and the first author, though there is a serious difference between
the invariants: our elliptic invariant counts weights assigned partly to
individual tropical curves and partly to collections of tropical curves, and out invariant is not
always multiplicative over the vertices of the counted tropical curves as was the case for all
tropical enumerative invariants known before. As a consequence
we define elliptic broccoli curves and elliptic broccoli invariants as well as
elliptic tropical descendant invariants for any toric surface.
\end{abstract}

\tableofcontents

\section{Introduction}

Refined (i.e. depending on a formal parameter $y$) tropical enumerative invariants
were introduced by F. Block and L. G\"ottsche \cite{BG}
(see also \cite{IM} for proof of the invariance), and 
they showed that, for $y=1$, its value is the Gromov-Witten tropical invariant, while for
$y=-1$ it coincides with the Welschinger tropical invariant associated with totally real
point constraints.
Under appropriate conditions, enumeration
of plane trivalent tropical curves with Mikhalkin and
Welschinger weights gives Gromov-Witten and Welschinger invariants of toric del Pezzo surfaces, respectively
(see \cite[Theorems 1 and 6]{Mi}). Mikhalkin \cite{Mi1} observed that in specific situations,
the refined invariant itself has an algebro-geometric enumerative meaning.
Note the Block-G\"ottsche invariants count tropical curves with marked points on the edges, and they are
uniformly defined for any genus (cf. \cite[Theorem 4.8]{GM} and \cite[Theorem 1]{IKS} for the complex
and real specifications discovered earlier).


In \cite{GMS}, A. Gathmann,  H. Markwig, and the first author introduced rational tropical broccoli invariants,
tightly related to Welschinger invariants for mixed, real and complex conjugate point constraints.
L. G\"ottsche jointly with the first author \cite{GS}
suggested a refinement of the rational broccoli invariant, which evaluates to the rational
broccoli invariant for $y=-1$ and to a certain rational logarithmic descendant invariant
of toric surfaces as $y=1$.
(see \cite{KM,MR} for detailed treatment of tropical descendant invariants).

The {\bf main goal} of this paper is to define elliptic (tropical) broccoli curves,
elliptic broccoli invariants and their refinement.
An attempt to define elliptic broccoli curves
was undertaken by the first author in \cite[Section 6]{S},
which basically led to indication of difficulties in this task.
These difficulties come from the fact
that the structure of the moduli space of plane tropical curves of a
positive genus with marked vertices is
much more complicated than that for genus zero. This stands in contrast to the
moduli spaces of plane tropical curves without marked vertices, where the difference
between genus zero and positive genera is quite simple (see \cite{GM}). In particular, this
allows one to easily
define invariants counting tropical curves of any genus with markings only on edges
\cite{GM,Mi,IKS,IM}.

The {\bf main outcome} of our work is
\begin{itemize}\item
a refined tropical invariant that counts elliptic plane tropical curves having
markings on edges and at vertices and
passing through appropriately many generic points; the new invariant naturally extends the refined
rational broccoli invariant (see Theorem \ref{t1} in Section \ref{sec1});
\item an elliptic broccoli invariant and a definition of elliptic broccoli curves
(geometric description in Definition \ref{dd1}, Section \ref{broccoli}, and numerical
characterization in
Proposition \ref{cor1new}, Section \ref{char});
\item an elliptic tropical descendant invariant $\big\langle\tau_0(2)^{n_e}
\tau_1(2)^{n_v}\big\rangle^1_\Delta$ (Proposition \ref{cor2} and Remark \ref{rr2} in Section
\ref{char}) that naturally extends the rational
tropical descendant invariant $\big\langle\tau_0(2)^{n_e}
\tau_1(2)^{n_v}\big\rangle^0_\Delta$ studied in \cite{MR}.
\end{itemize}

We show
(Proposition \ref{p5} in Section \ref{sec-com}) that 
our refined invariant is a symmetric
Laurent polynomial and explain how to compute it via a lattice path algorithm
very similar to that in \cite[Section 9]{MR}.

We also point out two main differences
between the Block-G\"ottsche and refined rational broccoli invariants on one side
and the refined elliptic broccoli invariant
on the other side:
\begin{itemize}
\item Refined weights of elliptic tropical curves in count are no longer
products of refined weights of
vertices.
\item Some tropical curves are counted with a joint refined weight, and it
is not clear whether
one can split such a joint weight into individual weights so that the invariant
will become local
(in the sense that, assigning to a top-dimensional face of the moduli space of the counted curves
the weight of a generic element of the
face, one obtains a topological cycle relative to infinity).
\end{itemize}

The 
questions 
on the enumerative meaning of the elliptic invariants and on possible extension of
such invariants to all positive genera
still remain open, and we plan to address it in a forthcoming paper.

\medskip

{\bf Acknowledgements}. A part of this work has been performed during the stay of the
first author
at the Tel Aviv University as a postdoc position supported from the Hermann-Minkowski-Minerva
Center for Geometry at the Tel Aviv University and from the Israel Science Foundation
grant no. 178/13. The second author has been supported by the German-Israeli Foundation
grant no. 1174-197.6/2011 and by the Israel Science Foundation grant no. 176/15.
A substantial part
of this work has been done during the stay of the second author at the
Max-Planck-Institut f\"ur Mathematik (MPIM), Bonn, in September 2015. The second
author is very grateful to the MPIM for hospitality and excellent working conditions.

\section{Plane marked tropical curves}
We shortly recall some basic definitions concerning tropical curves, adapted to our setting
(for details, see \cite{GM,GS,Mi}).

\subsection{Abstract and plane tropical curves}
An {\it
abstract tropical curve} is a finite connected 
compact graph $\overline\Gamma$ without bivalent
vertices such that
the complement $\Gamma=\overline\Gamma\setminus\overline\Gamma^{\;0}_\infty$ to the set $\overline\Gamma^{\;0}_\infty$
of univalent vertices
is a metric graph whose non-closed edges (called {\it ends}) are isometric to $[0,\infty)$
(i.e. univalent vertices are infinitely far from
their neighboring vertices). Denote by $\Gamma^0$, $\Gamma^1$, and
$\Gamma^1_\infty$ the sets of vertices, edges, and ends of
$\Gamma$, respectively. The {\it
genus} of $\Gamma$ is $g(\Gamma)=b_1(\Gamma)$. A tropical curve
$\overline\Gamma$ is called {\it
trivalent} if all the vertices of $\Gamma$ are trivalent. A {\it
marked tropical curve} is a pair
$(\overline\Gamma,\bp)$, where $\bp$ is an ordered sequence of distinct
points of $\Gamma$.
We say that the complement $\overline\Gamma\setminus\bp$ is {\it regular}, if each component
if this set is
simply connected and contains exactly one univalent vertex. The edges of the closure of
a component of a regular set $\overline\Gamma\setminus\bp$ admit a unique orientation
(called {\it canonical}) such that the marked points are the only sources and the univalent
vertex of $\overline\Gamma$ is the only sink. A{\it labeled abstract marked tropical curve}
is an abstract tropical curve with an ordered set
$\overline\Gamma^{\;0}_\infty$ of univalent vertices. We denote a labeled marked
tropical curve by $(\overline
\Gamma^{\;\lab},\bp)$.

An marked plane tropical curve is a tuple $(\overline\Gamma,\bp,h)$, where
$(\overline\Gamma,\bp)$ is a marked abstract tropical curve, and $h:\Gamma\to\R^2$ is
a proper map
such that
\begin{itemize}
\item the restriction $h\big|_E$ to any edge $E\in\Gamma^1$ is a
nonconstant affine-integral map;
\item for any vertex $V\in\Gamma^0$ holds the following balancing condition
\begin{equation}\sum_{V\in E,\ E\in\Gamma^1}D(h\big|_E)(\overline a_V(E))=0\ ,
\label{e1}\end{equation} where
$\overline a_V(E)$ is the unit tangent vector to $E$
at $V$ (oriented along $E$);
\item if, for some vertex $V\in\Gamma^0$, $\dim\Span\big\{D(h\big|_E)(\overline a_V(E)),\ V\in E,\ E\in\Gamma^1\big\}=1$,
then $V\in\bp$ (we call such vertices {\it collinear}).
\end{itemize} The latter condition excludes parasitic parameters, coordinates of images of collinear unmarked vertices.

Notice that each vector $D(h\big|_E)(\overline
a_V(E))$ has integral coordinates, and (if nonzero) it can
be written as $D(h\big|_E)(\overline a_V(E))=m\overline v$, where $\overline v\in\Z^2\setminus\{0\}$ is primitive and $m$ is a positive integer
(called the {\it weight of the edge} $E$ and denoted $\mt(E,h)$).
The {\it degree} of the plane tropical curve
$(\overline\Gamma,h)$ is the (unordered) multiset of vectors
$$\Delta(\overline\Gamma,h)=\big\{D(h\big|_E)(\overline a_V(E)),\ E\in\Gamma^1_\infty\big\}\ .$$
The balancing condition yields that $\Delta(\overline\Gamma,h)$ is a {\it balanced}
multiset, i.e.
$\sum_{\overline a\in\Delta(\overline\Gamma,h)}\overline a=0$.
We call $\Delta(\overline\Gamma,h)$ {\it primitive} if it contains only primitive
vectors, and we call
$\Delta(\overline\Gamma,h)$ {\it nondegenerate}, if $\Span_\R\Delta(\overline\Gamma,h)=\R^2$.

The push-forward $T=h_*(\Gamma)\subset\R^2$ is an embedded tropical curve with the Newton
polygon $P(\Delta)$ built of vectors $\overline a\in\Delta(\overline\Gamma,h)$
rotated by $\pi/2$. By $p_a(\Delta)$ we denote the number of interior integral points
of $P(\Delta)$, the arithmetic genus of a curve in the tautological linear system
on the toric surface associated with the polygon $P(\Delta)$.
There is a natural duality between the edges and vertices of $T$ on one side and the
edges and
polygons of a certain (dual) subdivision of $P(\Delta)$. We denote the dual object by
${\mathcal D}(*)$.
Any edge $E$ of $T$ possesses a natural weight $\mt(E)$, which can be
viewed as the lattice length of the dual edge ${\mathcal D}(E)$.

Observe also that our definition of a plane tropical curves
does not allow loops in $\Gamma$ formed by one edge.
If a cycle contains exactly two edges, then these edges join the same two
vertices of $\Gamma$, and they are mapped by $h$ to the same segment (though may
have different weights). We call such cycles {\it collinear}.

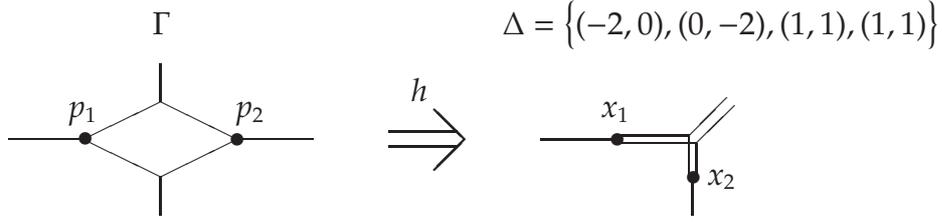
\begin{figure}
\setlength{\unitlength}{1cm}
\begin{picture}(12,4)(-1,1)
\thinlines

\put(0,2){\line(1,0){1}}\put(1,2){\line(2,-1){1}}\put(1,2){\line(2,1){1}}
\put(2,1.5){\line(2,1){1}}\put(2,2.5){\line(2,-1){1}}\put(3,2){\line(1,0){1}}
\put(2,1.5){\line(0,-1){0.5}}\put(2,2.5){\line(0,1){0.5}}
\put(8,2.05){\line(1,0){0.95}}\put(8,1.95){\line(1,0){1.05}}
\put(8.95,1.5){\line(0,1){0.55}}\put(9.05,1.5){\line(0,1){0.45}}
\put(8.95,2.05){\line(1,1){0.5}}\put(9.05,1.95){\line(1,1){0.5}}

\thicklines

\put(5,2.1){\line(1,0){0.9}}\put(5,1.9){\line(1,0){0.9}}\put(5.6,2.4){\line(1,-1){0.4}}
\put(5.6,1.6){\line(1,1){0.4}}\put(7,2){\line(1,0){1}}\put(9,1){\line(0,1){0.5}}

\put(0.9,1.9){$\bullet$}\put(2.9,1.88){$\bullet$}\put(7.9,1.88){$\bullet$}\put(8.9,1.38){$\bullet$}

\put(0.8,2.3){$p_1$}\put(3,2.3){$p_2$}\put(7.8,2.3){$x_1$}\put(9.2,1.4){$x_2$}
\put(1.9,3.4){$\Gamma$}\put(6.5,3.4){$\Delta=\big\{(-2,0),(0,-2),(1,1),(1,1)\big\}$}\put(5.3,2.5){$h$}

\end{picture}
\caption{Fake elliptic conic}\label{fec}
\end{figure}

\begin{remark}
When considering curves of genus $g>0$, we should assume that $g\le p_a(\Delta)$.
Otherwise, one may encounter fake
high genus curves like an ``elliptic conic" shown in Figure \ref{fec}.
\end{remark}

\subsection{Moduli spaces of tropical curves}
Denote by ${\mathcal M}_{g,n,k}$, where $n\ge0$, the moduli space of marked tropical curves
$(\overline\Gamma,\bp)$ of genus $g$
with an ordered configuration $\bp$ of $n$ distinct points of $\Gamma$ and
such that $|\overline\Gamma^{\;0}_\infty|=k$. It is a polyhedral complex, whose cells
parameterize $n$-marked curves of a given
combinatorial type, while the parameters (possibly linearly dependent) are lengths
between neighboring vertices and/or marked points.
The natural closure $\overline{\mathcal M}_{g,n,k}$ includes the curves
obtained by vanishing of some
parameters and may
contain curves of genus $<g$ (if $g$ is positive) as some cycles contract to points.

Furthermore, for any splitting $n=n_v+n_e$ with $n_v,n_e\ge0$, introduce
$${\mathcal M}_{g,(n_v,n_e),k}=\left\{(\overline
\Gamma,\bp)\in{\mathcal M}_{g,n,k}\ :\ 
\text{the first}\ n_v\ \text{points of}\ \bp\
\text{are vertices of}\ \Gamma
\right\}\ .$$ It is known that ${\mathcal M}_{g,(n_v,n_e),k}$ is a
finite polyhedral complex of dimension
$k+g-1+n_e$, and the open top-dimensional cones parameterize trivalent curves
$(\overline\Gamma,\bp)$ such that
$p_i\not\in\Gamma^0$, $i>n_v$, $p_i\in\bp$.
Denote by $\overline{\mathcal M}_{g,(n_v,n_e),k}$ its closure in
$\overline{\mathcal M}_{g,n,k}$.

Given an nondegenerate multiset $\Delta\subset\Z^2\setminus\{0\}$ such that $\sum_{\overline a\in\Delta}\overline a=0$
(further on referred to as {\it balanced multiset}), we have moduli spaces $\overline{\mathcal M}_{g,n}(\R^2,\Delta)$
and $$\overline{\mathcal M}_{g,(n_v,n_e)}(\R^2,\Delta)=
\left\{(\overline\Gamma,\bp,h)\in\overline{\mathcal M}_{g,n}(\R^2,\Delta)\ :
\ \renewcommand{\arraystretch}{1}
\begin{array}{c}(\overline\Gamma,\bp)\in\overline{\mathcal M}_{g,(n_v,n_e),k}\\ k=|\Delta|,\ n_v+n_e=n\end{array}\right\},$$
of marked plane tropical curves
of degree $\Delta$ and genus $g$. We also have a natural
evaluation map
$\Ev_n:\overline{\mathcal M}_{g,n}(\R^2,\Delta)\to\R^{2n}$ that takes
$(\overline\Gamma,\bp,h)$ to the vector $h(\bp)\in\R^{2n}$.
The following fact is well-known (see \cite{GM,Mi,S}):

\begin{lemma}\label{l5}
The space $\overline{\mathcal M}_{g,n}(\R^2,\Delta)$ is a finite polyhedral complex
of (not necessarily pure) dimension $|\Delta|+g-1+n$.
The space $\overline{\mathcal M}_{g,(n_v,n_e)}(\R^2,\Delta)$ is a finite polyhedral complex of
(not necessarily pure) dimension $|\Delta|+g-1+n_e$.
\end{lemma}

Suppose that $n=|\Delta|+g-1$ and denote by $\overline{\mathcal M}^{\:e}_{g,}(\R^2,\Delta)$
the closure in
$\overline{\mathcal M}_{g,n}(\R^2,\Delta)$ of the union of open cells of dimension
$2n=|\Delta|+g-1+n$ whose
$\Ev_{n}$-images have dimension $2n$ ({\it enumeratively essential
cells}). Respectively,
suppose that $2n_v+n_e=|\Delta|+g-1$ and
denote by $\overline{\mathcal M}^{\:e}_{g,(n_v,n_e)}(\R^2,\Delta)$ the closure in
$\overline{\mathcal M}_{g,(n_v,n_e)}(\R^2,\Delta)$ of the union of open cells of dimension
$2n=|\Delta|+g-1+n_e$ whose $\Ev_n$-images have dimension $2n$. By
$\Ev^e_n$,
resp. $\Ev^e_{n_v,n_e}$ we denote the
restriction of $\Ev_n$ to $\overline{\mathcal M}^{\:e}_{g,n}(\R^2,\Delta)$ and
$\overline{\mathcal M}^{\:e}_{g,(n_v,n_e)}(\R^2,\Delta)$, respectively.

\subsection{Labeled tropical curves and their moduli spaces}\label{lab}
A labeling of an abstract tropical curve $\overline\Gamma$ is a linear order on the set of the ends $\Gamma^1_\infty$. The moduli spaces of tropical curves considered in the preceding section, admit
natural labeled counterparts, in particular, fixed an arbitrary order on $\Delta$, we have
$\overline{\mathcal M}^{\;\lab}_{g,(n_v,n_e)}(\R^2,\Delta)$, the space of labeled marked
plane tropical curves of genus $g$, considered up to automorphisms respecting the labeling and the order of $\Delta$. The forgetful map
$$
\pi_{g,(n_v,n_e)}:\overline{\mathcal M}^{\:\lab}_{g,(n_v,n_e)}(\R^2,\Delta)\to\overline{\mathcal M}_{g,(n_v,n_e)}(\R^2,\Delta)$$
is surjective and finite, and, for any element $(\overline\Gamma,\bp,h)\in
\overline{\mathcal M}_{g,(n_v,n_e)}(\R^2,\Delta)$, we have
\begin{equation}\left|(\pi_{g,(n_v,n_e)})^{-1}(\overline\Gamma,\bp,h)\right|=\frac{|\Delta|!}{|\Aut(\overline
\Gamma,\bp,h)|}\ .\label{elab}\end{equation}
If, $2n_v+n_e=|\Delta|+g-1$, we denote by
$\Ev^{\;e,\lab}_{n_v,n_e}:\overline{\mathcal M}^{\:e,\lab}_{g,(n_v,n_e)}(\R^2,\Delta)\to\R^{2n}$ the evaluation map restricted to the (closure of the) enumeratively essential cells.

\subsection{Evaluation map in codimension zero and one} Below we will introduce a refined
enumerative quantity related to the count of elliptic tropical curves
passing through a generic configuration of points in $\R^2$, and we will prove its invariance
via the study of its behavior along generic paths in the space of point constraints. So, we stratify the space of constraints identified with $\R^{2n}$
for an appropriate $n$ into cells of full dimension $2n$, walls of codimension one,
and the complement of codimension two which can be avoided by generic paths.

Throughout this section we
suppose that $\Delta\subset\Z^2\setminus\{0\}$ is a balanced multiset.
Furthermore, to simplify formulations we introduce the following notion:
Given an abstract marked tropical curve $(\overline\Gamma,\bp)$, we call $(\Gamma\setminus\overline\Gamma',\bp
\setminus\overline\Gamma')$, where $\overline\Gamma'\subset\overline\Gamma$ is a closed subgraph
(possibly empty),
{\it simple} if either $\Gamma\setminus\overline\Gamma'$ is trivalent, or all but one vertices of
$\Gamma\setminus\overline\Gamma'$ are trivalent and there is a collinear cycle joining
a trivalent marked vertex with a four-valent unmarked vertex.


\begin{lemma}[\cite{GM,GS,Mi}]\label{l1} Let $n=|\Delta|+g-1$. Then, for each
element $\bx$ belonging to some open dense subset of $\R^{2n}$,
the preimage $\Ev_n^{-1}(\bx)$ is non-empty and finite.
Furthermore, for each curve $(\overline\Gamma,\bp,h)\in\Ev_n^{-1}(\bx)$, the graph $\Gamma$ is
trivalent, $\bp\cap\Gamma^0=\emptyset$, and the set $\overline\Gamma\setminus\bp$
is regular.
\end{lemma}

\begin{lemma}[\cite{GM}, Section 4, and \cite{GS}, Section 2]\label{l2}
Under the hypotheses of Lemma \ref{l1}, the map $\Ev_n^e$ is onto, and its target space
splits into the disjoint union
$$\R^{2n}=X^{2n}\cup X^{2n-1}\cup X^{2n-2}\ ,$$ where
\begin{enumerate}\item[(1)] $X^{2n}$ is the union of open polyhedra of dimension $2n$, and,
for each element $\bx\in X^{2n}$, its preimage $(\Ev_n^e)^{-1}(\bx)$ satisfies the conclusions of
Lemma \ref{l1};
\item[(2)] $X^{2n-1}$ is the union of open polyhedra of dimension $2n
-1$, for each $\bx\in
X^{2n-1}$ the preimage $(\Ev_n^e)^{-1}(\bx)$ is non-empty, finite, and the curves
$(\overline\Gamma,\bp,h)\in(\Ev_n^e)^{-1}(\bx)$ are as follows:
\begin{enumerate}
\item[(2i)] either $\Gamma$ is trivalent, having precisely one marked vertex,
\item[(2ii)] or all but one vertices of $\Gamma$ are trivalent, one vertex is four-valent,
unmarked, $\bp\cap\Gamma^0=\emptyset$,
and the set $\overline\Gamma\setminus\bp$ is regular,
\item[(2iii)] or, in case $g>0$, all but two vertices of $\Gamma$ are trivalent, the remaining
two vertices are four-valent, and they are joined by a couple of edges that are mapped to
the same segment in $\R^2$, and the set $\overline\Gamma\setminus\bp$ is regular.
\end{enumerate}
\item[(3)] $\dim X^{2n-2}\le2n-2$.
\end{enumerate}
\end{lemma}

\begin{lemma}\label{l4}
Let $n=n_v+n_e$, $n_v>0$, $n_e\ge0$, and $2n_v+n_e=|\Delta|+g-1$. Then, for a generic element
$\bx\in\R^{2n}$, the preimage
$(\Ev^e_{n_v,n_e})^{-1}(\bx)$ is finite. Furthermore, one has:

(1) If $g=0$, then each curve $(\overline\Gamma,\bp,h)\in(\Ev^e_{n_v,n_e})^{-1}(\bx)$
satisfies the following conditions
\begin{enumerate}
\item[(1i)] $|\bp\cap\Gamma^0|=n_v$, and the set $\overline\Gamma
\setminus\bp$ is regular;
\item[(1ii)] $\Gamma$ is trivalent.
\end{enumerate}

(2) If $g=1$, then each curve $(\overline\Gamma,\bp,h)\in(\Ev^e_{n_v,n_e})^{-1}(\bx)$
satisfies the following conditions
\begin{enumerate}
\item[(2i)] $|\bp\cap\Gamma^0|=n_v$, and the set $\overline\Gamma\setminus\bp$ is regular;
\item[(2ii)] $(\Gamma,\bp)$ is simple.
\end{enumerate}
\end{lemma}

{\bf Proof.} The case of $g=0$ is settled in \cite{GMS}. The second statement is a part of
Lemma \ref{l6} below.
\proofend

\begin{lemma}\label{l3}
Let $n=n_v+n_e$, $n_v>0$, $n_e\ge0$, and $2n_v+n_e=|\Delta|-1$.
Then the target space of $\Ev^e_{n_v,n_e}:\overline{\mathcal M}^{\;e}_{0,(n_v,n_e)}(\R^2,\Sigma)
\to\R^{2n}$, splits into the disjoint union
$$\R^{2n}=X^{2n}\cup X^{2n-1}\cup X^{2n-2}\ ,$$ where
\begin{enumerate}
\item[(1)] $X^{2n}$ is the union of open polyhedra of dimension $2n$, and,
for each element $\bx\in X^{2n}$, its preimage $(\Ev_{n_v,n_e}^e)^{-1}(\bx)$ is finite and satisfies the
conclusions of Lemma \ref{l4}(1);
\item[(2)] $X^{2n-1}$ is the union of open polyhedra of dimension $2n
-1$, for each $\bx\in
X^{2n-1}$ the preimage $(\Ev_{n_v,n_e}^e)^{-1}(\bx)$ is finite, and the curves
$(\overline\Gamma,\bp,h)\in(\Ev_{n_v,n_e}^e)^{-1}(\bx)$ are as follows:
\begin{enumerate}
\item[(2i)] either $\Gamma$ is trivalent, and $|\bp\cap\Gamma^0|=n_v+1$,
\item[(2ii)] or all but one vertices of $\Gamma$ are trivalent, one vertex is four-valent,
unmarked, $|\bp\cap\Gamma^0|=n_v$, and the set $\overline\Gamma\setminus
\bp$ is regular,
\item[(2iii)] or all but one vertices of $\Gamma$ are trivalent, one vertex is four-valent,
marked, and $|\bp\cap\Gamma^0|=n_v$;
\end{enumerate}
\item[(3)] $\dim X^{2n-2}\le2n-2$.
\end{enumerate}
\end{lemma}

{\bf Proof.} See \cite[Section 3]{GMS} and \cite[Section 4]{GS}.\proofend

\begin{lemma}\label{l6}
Let $n=n_v+n_e$, $n_v>0$, $n_e\ge0$, and $2n_v+n_e=|\Delta|$.
Then the target space of $\Ev^e_{n_v,n_e}:\overline{\mathcal M}^{\;e}_{1,(n_v,n_e)}(\R^2,\Delta)\to\R^{2n}$,
splits into the disjoint union
$$\R^{2n}=X^{2n+n_\infty}\cup X^{2n-1}\cup X^{2n-2}\ ,$$ where
\begin{enumerate}
\item[(1)] $X^{2n}$ is the union of open polyhedra of dimension $2n$, and,
for each element $\bx\in X^{2n}$, its preimage $(\Ev_{n_v,n_e}^e)^{-1}(\bx)$ satisfies the
conclusions of Lemma \ref{l4}(2);
\item[(2)] $X^{2n-1}$ is the union of open polyhedra of dimension $2n-1$, for each $\bx\in
X^{2n-1}$ the preimage $(\Ev_{n_v,n_e}^e)^{-1}(\bx)$ is finite and consists of the
following curves
$(\overline\Gamma,\bp,h)$:
\begin{enumerate}
\item[(2i)] either $\Gamma$ is elliptic, $|\bp\cap\Gamma^0|=n_v+1$,
the set $\overline\Gamma\setminus\bp$ contains exactly one bounded component, and
the marked graph $(\Gamma,\bp)$ is simple;
\item[(2ii)] or $\Gamma$ is elliptic, $|\bp\cap\Gamma^0|=n_v$, the set $\overline
\Gamma\setminus\bp$ is
regular, there is an unmarked four-valent vertex $V\in\Gamma^0$ which
neither belongs to a collinear cycle, nor is neighboring to a collinear marked
trivalent vertex belonging to a
collinear cycle, and the marked graph $(\Gamma\setminus\{V\},\bp)$ is simple;
\item[(2iii)] or $\Gamma$ is elliptic, has no collinear cycle, $|\bp\cap\Gamma^0|=n_v$,
all but one vertices of
$\Gamma$ are trivalent,
one vertex is four-valent, marked, and exactly one edge of the unique bounded component of $\overline
\Gamma\setminus\bp$ is incident to that vertex;
\item[(2iv)] or $\Gamma$ is elliptic, has no collinear cycle, $|\bp\cap\Gamma^0|=n_v$,
all but one vertices of
$\Gamma$ are trivalent,
one vertex is four-valent, marked, and exactly two edges of the unique bounded component of $\overline
\Gamma\setminus\bp$ are incident to that vertex;
\item[(2v)] or $\Gamma$ is elliptic, has a collinear cycle, $|\bp\cap\Gamma^0|=n_v$, and
\begin{enumerate}
\item[(2v-a)] either all but two vertices of $\Gamma$ are trivalent,
two unmarked four-valent vertices are joined by a collinear cycle, and the set
$\overline\Gamma\setminus\bp$ is regular,
\item[(2v-b)] or all but two vertices of $\Gamma$ are trivalent,
two four-valent vertices, one marked, the other unmarked, are joined by a collinear cycle
containing an additional marked point,
\item[(2v-c)] or all but two vertices of $\Gamma$ are trivalent,
two marked four-valent vertices are joined by a collinear cycle;
\end{enumerate}
\item[(2vi)] or $\Gamma$ is elliptic, and all but two vertices of $\Gamma$ are trivalent,
two four-valent vertices, one marked, the other unmarked, are joined by a collinear cycle
containing no extra marked point, and the set
$\overline\Gamma\setminus\bp$ is not regular;
\item[(2vii)] or $\Gamma$ is rational, trivalent,
$|\bp\cap\Gamma^0|=n_v$, and the set $\overline\Gamma\setminus\bp$ contains exactly one bounded component;
\item[(2viii)] or $\Gamma$ is elliptic, has a collinear cycle, $|\bp\cap\Gamma^0|=n_v$, and
\begin{enumerate}
\item[(2viii-a)] either all but one vertices of $\Gamma$ are trivalent, one five-valent
vertex is unmarked, joined with a marked vertex by a collinear cycle,
and the set $\overline\Gamma\setminus\bp$ is regular,
\item[(2viii-b)] or all but two vertices of $\Gamma$ are trivalent, two unmarked four-valent
vertices are neighboring
to a collinear marked trivalent vertex, and the set $\overline\Gamma\setminus\bp$ is regular;
\end{enumerate}
\end{enumerate}
\item[(3)] $\dim X^{2n-2}\le2n-2$.
\end{enumerate}
\end{lemma}

{\bf Proof.} The statement can be extracted from \cite[Chapter 6]{S}.
For the reader's convenience we shall prove it here. Given a curve
$(\overline\Gamma,\bp,h)\in\overline{\mathcal M}^{\;e}_{1,(n_v,n_e)}(\R^2,\Delta)$, by
$\Def(\overline\Gamma,\bp,h)$ we denote the (open) cell of $\overline{\mathcal M}^{\;e}_{1,(n_v,n_e)}(\R^2,\Delta)$
parameterizing curves of the same combinatorial type as $(\overline\Gamma,\bp,h)$.

{\it Step 1.} Suppose that $\bx\in\R^{2n}$ is generic.
Observe that $|\bp\cap\Gamma^0|=n_v$ and that the set $\overline\Gamma\setminus\bp$ must
be regular, since otherwise
it would contain a bounded component, and hence a condition to the images of the endpoints of that component
(cf. \cite[Formula (4)]{Sh}). Next, denote by $s_i$ the valency of the marked vertex $p_i\in\bp$,
$i=1,...,n_v$,
set $u=|\Gamma^0\setminus\bp|$, and denote by $t_j$, $j=1,...,u$, the valencies of the
(somehow ordered) unmarked vertices of $\Gamma$. Counting the marked points $p_i
\in\bp$, $i>n_v$ as
bivalent vertices of $\overline\Gamma$ and respectively defining the set $\Gamma^1$ of
edges of $\Gamma$, we get
$$2|\Gamma^1|=\sum_{i=1}^{n_v}s_i+\sum_{j=1}^ut_j+2n_e+|\Delta|=\sum_{i=1}^{n_v}s_i+
\sum_{j=1}^ut_j+2n_v+3n_e+n_\infty\ ,$$
$$|\overline\Gamma^{\;0}|=u+n_v+n_e+|\Delta|=u+3n_v+2n_e\ ,$$ which together with the genus condition
$|\overline\Gamma^{\;0}|=|\Gamma^1|$ yields
\begin{equation}\sum_{i=1}^{n_v}(4-s_i)+n_e=\sum_{j=1}^u(t_j-2)\ .\label{edop1}\end{equation}
Now, denote by $\overline\Gamma_k$, $k=1,...,|\Delta|$, the closures
of the components of $\Gamma\setminus\bp$ obtained by adding a vertex to each non-closed edge.
Next, denote by $s_{ik}$, $i=1,...,n_v$, resp. $t_{jk}$, $j=1,...,u$, the number of edges of $\overline\Gamma_k$
incident to the respective marked and unmarked vertices, denote by $n_{e,k}$
the number of
bivalent marked vertices of $\overline\Gamma$ belonging to $\overline\Gamma_k$,
and, finally, denote by $u_k$ the number of vertices of $\overline\Gamma_k$
of valency $>2$, $k=1,...,|\Delta|$. Then the regularity condition
$|\overline\Gamma_k^1|+1=|\overline\Gamma_k^0|$, $k=1,...,|\Delta|$, implies that
$$2|\overline\Gamma_k^{\;1}|+2=\sum_{i=1}^{n_v}s_{ik}+\sum_{j=1}^ut_{jk}+
2n_{e,k}+3=2|
\overline\Gamma_k^{\;0}|=
2\sum_{i=1}^{n_v}s_{ik}+2n_{e,k}+2n_{\infty,k}+2u_k+2\ ,$$
$k=1,...,|\Delta|$,
and hence
\begin{equation}\sum_{i=1}^{n_v}(s_i-2)+n_e=\sum_{j=1}^u(t_j-2)\ ,\label{edop2}\end{equation}
which together with (\ref{edop1}) yields that all marked vertices $p_1,...,p_{n_v}$ are trivalent.
Furthermore, in the above notation, $s_{ik}\le2$ for all
$i=1,...,n_v$, $k=1,...,|\Delta|$, and $s_{ik}=2$ for at most one pair $(i,k)$, since each
pair like that yields a cycle of $\overline\Gamma$. If
$s_{ik}=1$ for all $k$, then all unmarked vertices are trivalent (cf. \cite[Proposition 2.23]{Mi}).
If $s_{ik}=2$ for some pair $(i,k)$, and there is an unmarked vertex $V$ of $\overline\Gamma_k$
of valency $\ge4$, then $\overline\Gamma_k\setminus\{V\}$ contains at least three bounded trees, and hence
a restriction to the position of $h(\overline\Gamma_k)\cap\bx$ unless $V$ is
four-valent and
joined by a couple of edges with some $p_i\in\bp$.

{\it Step 2.} It is easy to see that the combinatorial types
listed in item (2) of Lemma are degenerations of the top-dimensional combinatorial types of
item (1) and that
$\dim\Ev^e_{n_v,n_e}(\Def(\overline\Gamma,\bp,h))=2n+n_\infty-1$ for each of them.
We shall show that other degenerations satisfy $\dim\Ev^e_{n_v,n_e}
(\Def(\overline\Gamma,\bp,h))\le2n-2$.
The following two simple observations will be used below.

First, we note that
if $\Gamma\setminus\bp$ contains at least two bounded components then
$\dim\Ev^e_{n_v,n_e}(\Def(\overline\Gamma,\bp,h))\le2n-2$. Indeed,
the closures of two
such bounded components have at most two common points (since the genus is at most $1$), and hence
one gets two linearly independent conditions to the position of $h(\bp)$ (cf.
\cite[Formula (4)]{Sh}).

Second, suppose that $h:\gamma\to\R^2$ is a plane rational tropical curve of degree
$\delta$ having vertices of valencies $s_1,...,s_r$, $r\ge1$, then
$\dim\Def(\gamma,h)=|\delta|-1-\sum_{i=1}^r(s_i-3)$. Indeed, this immediately follows
from the formula $|\delta|-3
-\sum_{i=1}^r(s_i-3)$ for the number of bounded edges.

{\it Step 3.} Let $(\overline\Gamma,\bp,h)\in(\Ev_{n_v,n_e}^e)^{-1}(\bx)$ be rational.
The computation similar to that in Step 1
leads to relations
\begin{equation}\begin{cases}&\sum_{i=1}^{n_v}(4-s_i)+n_e=\sum_{j=1}^u(t_j-2)+2+
\sum_{i=1}^{n_e}(s'_i-2),\\
&\sum_{i=1}^{n_v}(s_i-2)+\sum_{i=1}^{n_e}(s'_i-2)+n_e=\sum_{j=1}^u(t_j-2)+2q,\end{cases}
\label{edop3}\end{equation}
where $s'_i$ is the valency of the point $p_{i+n_v}\in\bp$, $i=1,...,n_e$, and
$q$ is the number of bounded connected components of $\overline\Gamma\setminus\bp$.
If we suppose that $\dim\Ev^e_{n_v,n_e}(\Def(\overline\Gamma,\bp,h))=2n-1$, then
$q\le1$, and from (\ref{edop3}) we immediately get that $q=1$, all marked vertices $p_i\in\bp$,
$i=1,...,n_v$ are trivalent, and all marked points $p_{i+n_v}
\in\bp$, $i=1,...,n_e$, are bivalent.
It can easily be derived that all unmarked vertices are trivalent, since otherwise it would
further reduce dimension of $\Ev^e_{n_v,n_e}(\Def(\overline\Gamma,\bp,h))$ (cf.
Step 2).
This, we fit the conditions of item (2vii).

{\it Step 4.}
Let $(\overline\Gamma,\bp,h)\in(\Ev_{(n_v,n_e)}^e)^{-1}(\bx)$ be elliptic, and some of
the marked points
$p_{i+n_v}$, $1\le i\le n_e$, have valency $\ge 3$. The computation as in Step 1 leads
to relations
\begin{equation}\begin{cases}&\sum_{i=1}^{n_v}(4-s_i)+n_e=\sum_{j=1}^u(t_j-2)+\sum_{i=1}^{n_e}
(s'_i-2),\\
&\sum_{i=1}^{n_v}(s_i-2)+\sum_{i=1}^{n_e}(s'_i-2)+n_e=\sum_{j=1}^u(t_j-2)+2q,\end{cases}
\label{edop4}\end{equation} where $\sum_{i=1}^{n_e}(s'_i-2)\ge1$. Thus, we obtain that $q=1$, $\sum_{i=1}^{n_e}(s'_i-2)=1$
and that all marked vertices $p_i\in\bp$, $i=1,...,n_v$, are trivalent. It also follows that
one point $p_{i+n_v}\in\bp$, $1\le i\le n_e$, is a trivalent vertex, while he others are bivalent.
Finally, as in Step 1, we derive that $(\overline\Gamma,\bp,h)$ meets
the conditions of item (2i).

{\it Step 5.} Suppose that $|\bp\cap\Gamma^0|=n_v$, and $(\overline\Gamma,\bp,h)\in(\Ev_{n_v,n_e}^e)^{-1}(\bx)$ is elliptic with a collinear cycle.
Introduce the rational curve $(\Gamma',\bp,h')$ by replacing the edges of the collinear
cycle with one edge,
whose weight is the sum of weights of the replaced edges. It belongs either to
$\overline{\mathcal M}_{0,(n_v,n_e)}(\R^2,\Delta)$, or to $\overline{\mathcal M}_{0,(n_v-1,n_e+1)}(\R^2,\Delta)$.

In the latter case, $2(n_v-1)+(n_e+1)=|\Delta|-1$, and hence Lemma \ref{l3}(2) applies and
leaves for $(\overline\Gamma,\bp,h)$ the only options
(2ii) and (2viii).

In the former case, the relation $2n_v+n_e=|\Delta|$ means that $\Gamma'\setminus\bp$
contains a bounded component
(which is unique in view of the assumption $\dim\Ev^e_{n_v,n_e}(\Def(\overline\Gamma,
\bp,h))=2n-1$).
In addition, the second observation in Step 2 yields that $\Gamma'$ must be trivalent,
and hence we are left with only options
(2v) and (2vi).

{\it Step 6.} At last suppose that $|\bp\cap\Gamma^{\;0}|=n_v$, and $(\overline\Gamma,\bp,h)$
is elliptic without collinear cycles.
If the set $\overline\Gamma\setminus\bp$
is regular, then (cf. the second observation in Step 2),
all but one vertices of $\Gamma$ are trivalent, one vertex is four-valent.
Furthermore, the four-valent vertex must be unmarked, since otherwise one would encounter
a bounded component of $\Gamma\setminus\bp$.
Thus, we fit conditions of item (2ii).
Suppose now that the set $\Gamma\setminus\bp$ contains a bounded component (unique in view
of the first observation
in Step 2). As in Step 4, we get relations (\ref{edop4}), where $s'_i=2$, $i=1,...,n_e$,
and $q=1$. It follows that
exactly one marked vertex is four-valent, while the others are trivalent. This four-valent
marked vertex must be incident to a bounded component
of $\Gamma\setminus\bp$, since otherwise it would not be a degeneration of a top-dimensional
combinatorial type. On the other hand,
this four-valent vertex cannot be incident to more than two edges of
the bounded component of $\Gamma\setminus\bp$, since otherwise the curve would have genus $>1$.
Thus, we fit either conditions of either item (2iii), or of item (2iv).
\proofend

\section{Refined count of plane marked tropical curves}\label{sec1}
Let $\Delta\subset\Z^2\setminus\{0\}$ be a balanced multiset. Given $\mu\in\Z$, set
$$[\mu]_y^-=\frac{y^{\mu/2}-y^{-\mu/2}}{y^{1/2}-y^{-1/2}},\quad
[\mu]_y^+=\frac{y^{\mu/2}+y^{-\mu/2}}{y^{1/2}+y^{-1/2}},
\quad[\mu]^*_y=\frac{1}{\mu}\cdot\frac{y^{\mu/2}-(-1)^\mu y^{-\mu/2}}{y^{1/2}-
(-1)^\mu y^{-1/2}}\ ,$$
$y$ being a formal parameter.

Let $(\overline\Gamma,h)$ be a plane tropical curve.
For any trivalent vertex
$V\in\Gamma^0$, set its {\it Mikhalkin's weight}:
$$\mu(\overline\Gamma,h,V)=\Big|D(h\big|_{E_1})(\overline
a_V(E_1))\times D(h\big|_{E_2})(\overline
a_V(E_2))\Big|\ ,$$ where $E_1,E_2\in\Gamma^1$ are distinct edges incident to $V$. Due
to the balancing condition
(\ref{e1}), this number does not depend on the choice of a pair of edges incident to $V$
and, in fact,
it equals the lattice area of the triangle ${\mathcal D}(h(V))$, dual to the vertex $h(V)$
of the tropical curve $T=h_*(\Gamma)$. For any univalent vertex $V\in\overline\Gamma^{\;0}_\infty$,
set $\mu(\overline\Gamma,h,V)=\mt(E,h)$, where $E$ is the edge of $\overline\Gamma$ incident to $V$.

\subsection{Trivalent curves without marked vertices}
Let $(\overline\Gamma,\bp,h)\in{\mathcal M}^e_{g,n}(\R^2,\Delta)$ be trivalent, $n=|\Delta|+g-1$,
$\bp\cap\Gamma^0=\emptyset$, and let the set $\overline\Gamma\setminus\bp$ be regular
(equivalently, $\Ev_n^e(\overline\Gamma,\bp,h)\in X^{2n}$ in the notation of
Lemma \ref{l2}). 

Along \cite[Definition 3.5]{BG}
the {\it Block-G\"ottsche weight} of the curve $(\overline\Gamma,\bp,h)$ (with $\bp_\infty=\emptyset$
and $n_\infty=0$)
is
$$\BG_y(\overline\Gamma,\bp,\bp_\infty,h)=
\prod_{V\in\Gamma^0}[\mu(\overline\Gamma,h,V)]_y^-\ .$$
Under some mild condiions on $\Delta$, 
$\BG_y(\Gamma,\bp,h)$ is a symmetric
Laurent polynomial in $y$ (see \cite[Proposition 2.3(4)]{IM}). 


\begin{proposition}[\cite{IM}, Theorem 1]
\label{p1}
For any nondegenerate balanced multiset $\Delta\subset\Z^2\setminus\{0\}$ and integers $g\ge0$,
$n=|\Delta|+g-1$, 
the expression
$$\BG_y(\Delta,g,\bx)=\sum_{(\Gamma,\bp,h)\in
(\Ev_n^e)^{-1}(\bx)}\BG_y(\overline\Gamma,\bp,h)$$
does not depend in the choice of $\bx\in X^{2n}$ (defined in Lemma \ref{l2}).
\end{proposition}


Furthermore,

\begin{proposition}[\cite{IM}, Corollary 2.4, and \cite{Mi}, Theorems 1 and 3]\label{p2}
Under the hypotheses of Proposition \ref{p1}, suppose additionally that
$\Delta$ is primitive. Then
\begin{enumerate}
\item[(1)] $\BG_1(\Delta,g)$ equals the number of irreducible complex curves in the toric
surface $\Tor(\Delta)$, belonging to the tautological linear system, having
genus $g$, and passing through
a generic configuration of $n$ points in $\Tor(\Delta)$; if $\Tor(\Delta)$
is del Pezzo, then
$\BG_1(\Delta,g)=GW_g(\Delta)$, the genus $g$ Gromov-Witten invariant of
$\Tor(\Delta)$ for the tautological
linear system;
\item[(2)] $\BG_{-1}(\Delta,g)$ equals the sum of Welschinger signs of
irreducible real curves in the toric surface
$\Tor(\Delta)$, belonging to the tautological linear system, having genus
$g$, and passing through
a generic configuration of $n$ real points that tropicalize into a configuration
from $X^{2n}$;
if $\Tor(\Delta)$ is del Pezzo and $g=0$, then $\BG_{-1}(\Delta,0)$ equals the
Welschinger invariant $W_0(\Delta)$ of $\Tor(\Delta)$ corresponding to the
tautological
linear system and totally real configurations of points.
\end{enumerate}
\end{proposition}

\subsection{Rational trivalent curves with marked vertices}
Let $(\overline
\Gamma,\bp,h)\in{\mathcal M}^{\;e}_{0,(n_v,n_e)}(\R^2,\Delta)$ be a rational, trivalent
plane tropical curve.
Assume that $n=|\Delta|-1$,
$|\bp\cap\Gamma^0|=n_v$, and the set $\overline
\Gamma\setminus\bp$ is regular. Following \cite[Definition 3.13]{GS}, define the {\it refined broccoli weight}
of $(\overline\Gamma,\bp,h)$ by
\begin{equation}
\GS_y(\overline\Gamma,\bp,h)=\frac{1}{|\Aut(\overline\Gamma,\bp,h)|}
\prod_{V\in\bp\cap\Gamma^0}[\mu(\overline\Gamma,h,V)]_y^+\cdot
\prod_{V\in\Gamma^0\setminus\bp}[\mu(\overline\Gamma,h,V)]_y^-\cdot\prod_{V\in\overline\Gamma^{\;0}_\infty}
[\mu(\overline\Gamma,h,V)]^*_y\ ,
\label{e2}\end{equation}
where $\Aut(\overline\Gamma,\bp,h)$ is the automorphism group of $(\overline\Gamma,\bp,h)$.

\begin{proposition}[\cite{GS}, Theorem 4.1]\label{p3}
For any balanced multiset $\Delta\subset\Z^2\setminus\{0\}$ and nonnegative integers
$n_v,n_e$ such that $2n_v+n_e=|\Delta|-1$, the expression
\begin{equation}\GS_y(\Delta,0,(n_v,n_e),\bx):=\sum_{(\overline\Gamma,\bp,h)\in
(\Ev_{(n_v,n_e)}^e)^{-1}(\bx)}
\GS_y(\overline
\Gamma,\bp,h)\label{ee1}\end{equation}
does not depend on 
the choice of $\bx\in X^{2n}$ (defined in Lemma \ref{l3}).
\end{proposition}

For the proof see \cite[Section 4]{GS}. Thus, $\GS_y(\Delta,0,(n_v,n_e))$ is a genus zero tropical invariant.
Its particular values have an important enumerative meaning:

\begin{proposition}[\cite{GS}, Corollary 3.14 and Lemmas 3.27, 3.29]\label{p4}
If $\Delta$ consists of the vectors $(1,0),(0,1),(-1,-1)$, each one appearing $d$ times,
then $\Tor(\Delta)\simeq\PP^2$, and
$\GS_1(\Delta,0,(n_v,n_e))$ equals the descendant invariant $\big\langle\tau_0(2)^{n_e}
\tau_1(2)^{n_v}\big\rangle^0_d$\ , while $\GS_{-1}(\Delta,0,(n_v,n_e))$ equals the Welschinger invariant
$W_{n_v}(\PP^2,d)$, counting (with signs) real plane rational curves of degree $d$ that pass through
$n_e$ real points and $n_v$ pairs of complex conjugate points in general position.
\end{proposition}

\subsection{Elliptic curves with marked vertices}\label{sec2.3}
Now we define refined weights of generic elliptic curves with marked vertices,
i.e., tropical curves
described in Lemma \ref{l4}(2).

Suppose that $\Delta\subset\Z^2\setminus\{0\}$ is a balanced, nondegenerate multiset, 
$2n_v+n_e=|\Delta|$, where $n_v>0$, $n_e\ge0$.

Denote by ${\mathcal M}'_{1,(n_v,n_e)}(\R^2,\Delta)\subset\overline{\mathcal M}^{\;e}_{1,(n_v,n_e)}(\R^2,\Delta)$
the subset formed by the elements $(\overline\Gamma,\bp,h)$ such that
$\Gamma$ is trivalent, elliptic, $|\bp\cap\Gamma^0|=n_v$, and the
set $\overline\Gamma\setminus\bp$ is regular.
Denote by ${\mathcal M}''_{1,(n_v,n_e)}(\R^2,\Delta)\subset\overline{\mathcal M}^{\;e}_{1,(n_v,n_e)}(\R^2,\Delta)$
the subset formed by the elements $(\overline\Gamma,\bp,h)$ such that
$\Gamma$ is elliptic with all but one vertices trivalent, one unmarked vertex is four-valent
and joined
with a marked vertex by a collinear cycle, $|\bp\cap\Gamma^0|=n_v$, and the set
$\overline\Gamma\setminus\bp$ is regular.

By Lemma \ref{l4}(2), the complement of ${\mathcal M}'_{1,(n_v,n_e)}(\R^2,\Delta)
\cup{\mathcal M}''_{1,(n_v,n_e)}(\R^2,\Delta)$ in
$\overline{\mathcal M}^{\;e}_{1,(n_v,n_e)}(\R^2,\Delta)$ has
positive codimension.

\begin{figure}
\setlength{\unitlength}{1cm}
\begin{picture}(12,7)(-1,0)
\thinlines

\put(5,0){\line(1,1){0.5}}\put(5.5,0.5){\line(1,-1){0.5}}
\put(5.5,0.5){\line(0,1){1}}\put(5.5,1.5){\line(1,1){0.5}}\put(5.5,1.5){\line(-1,1){0.5}}
\put(5.43,0.92){$\bullet$}
\put(6.7,1){$T\subset\R^2$}

\put(1,4){\line(1,1){0.5}}\put(1.5,4.5){\line(1,-1){0.5}}
\put(1.5,5){\line(0,1){0.5}}\put(1,4){\line(1,1){0.5}}
\put(1.5,5.5){\line(1,1){0.5}}\put(1.5,5.5){\line(-1,1){0.5}}
\qbezier(1.5,4.5)(1.3,4.75)(1.5,5)\qbezier(1.5,4.5)(1.7,4.75)(1.5,5)
\put(1.43,4.92){$\bullet$}
\put(1.3,6.3){$\Gamma$}

\put(5,4){\line(1,1){0.5}}\put(5.5,4.5){\line(1,-1){0.5}}
\put(5.5,4.5){\line(0,1){1}}\put(5.5,5.5){\line(1,1){0.5}}\put(5.5,5.5){\line(-1,1){0.5}}
\put(5.43,4.92){$\bullet$}
\put(5.3,6.3){$\Gamma'$}

\put(9,4){\line(1,1){0.5}}\put(9.5,4.5){\line(1,-1){0.5}}
\put(9.5,4.5){\line(0,1){0.5}}\put(9,4){\line(1,1){0.5}}
\put(9.5,5.5){\line(1,1){0.5}}\put(9.5,5.5){\line(-1,1){0.5}}
\qbezier(9.5,5)(9.3,5.25)(9.5,5.5)\qbezier(9.5,5)(9.7,5.25)(9.5,5.5)
\put(9.43,4.92){$\bullet$}
\put(9.3,6.3){$\Gamma$}

\thicklines

\put(3,5){\vector(1,0){1}}\put(8,5){\vector(-1,0){1}}
\put(3,3.5){\vector(1,-1){1}}\put(8,3.5){\vector(-1,-1){1}}
\put(5.5,3.5){\vector(0,-1){1}}
\put(3.4,5.2){$\pi$}\put(7.5,5.2){$\pi$}
\put(3.6,3.2){$h$}\put(5.7,3){$h'$}\put(7.3,3.2){$h$}

\put(1.7,5.3){$V_1$}\put(1.7,4.4){$V_2$}\put(1.1,4.9){$p$}
\put(9.7,5.3){$V_1$}\put(9.7,4.4){$V_2$}\put(9.1,4.9){$p$}
\put(5.7,5.3){$V_1$}\put(5.7,4.4){$V_2$}\put(5.1,4.9){$p$}
\put(5.7,1.3){$v_1$}\put(5.7,0.4){$v_2$}\put(5.1,0.9){$x$}

\end{picture}
\caption{Tropical curves with collinear cycles}\label{fig1}
\end{figure}
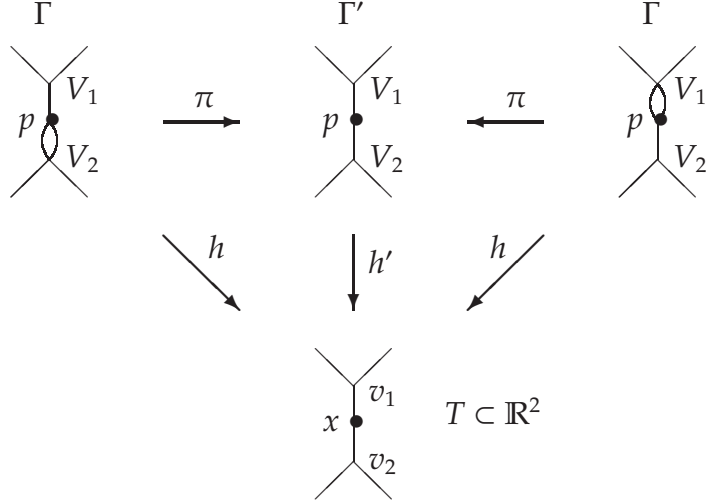

Let $(\overline
\Gamma,\bp,h)\in{\mathcal M}''_{1,(n_v,n_e)}(\R^2,\Delta)$. Then the marked trivalent
vertex belonging to
the collinear cycle is collinear, too. 
The map $\pi_{(\Gamma,h)}:\Gamma\to\Gamma'$, identifying edges of the
collinear cycle according to the $h$-image, defines a curve
$(\overline\Gamma',\bp,h')\in\overline{\mathcal M}^e_{0,(n_v-1,n_e+1)}(\R^2,\Delta)$ such that $h=h'\circ\pi_{(\Gamma,h)}$
(see Figure \ref{fig1}).
Furthermore, the corresponding morphism
$$(\overline
\Gamma,\bp,h)\in{\mathcal M}''_{1,(n_v,n_e)}(\R^2,\Delta)\overset{\pi}{\mapsto}(\overline
\Gamma',\bp,h')\in
\overline{\mathcal M}^{\;e}_{0,(n_v-1,n_e+1)}(\R^2,\Delta)$$ satisfies
$\Ev^e_{(n_v,n_e)}=\Ev^e_{(n_v-1,n_e+1)}\circ\pi$. A curve $(\overline\Gamma',\bp,h')\in
\Ima(\pi)$ has a distinguished marked point $p\in\bp$, the image of the marked vertex in the
collinear cycle. Denote by ${\mathcal M}''_{0,(n_v-1,n_e+1)}
(\R^2,\Delta)$ the moduli space of the tuples $(\overline\Gamma',\bp,p,h')$ such that
$(\overline\Gamma',\bp,h')\in\Ima(\pi)\subset\overline{\mathcal M}^{\;e}_{0,(n_v-1,n_e+1)}
(\R^2,\Delta)$, and $p\in\bp$ is the image
of the marked collinear vertex belonging to the collinear cycle. Note that the point $p$
belongs to some edge $E$ of $\Gamma'$ with
$h'(E)=e$ an edge of the embedded curve $(T,\bx)=h'_*(\Gamma',\bp)$ of
weight
$\mt(e)>1$. Now we introduce the following expressions.

If $e=h'(E)$ is the image of an unbounded edge $E$ of $\Gamma'$, $V_1$ the vertex of $E$,
$v_1=h'(V_1)$, we let
$$
\GS_y(\overline\Gamma',\bp,p,h')=\frac{1}{|\Aut(\overline\Gamma',\bp,p,h')|}\cdot
\Psi^{(1)}_z(m,\nu)\qquad\qquad\qquad\qquad$$
\begin{equation}\qquad\qquad\quad\times\prod_{V\in\bp\cap(\Gamma')^0}[\mu
(\overline
\Gamma',h',V)]_y^+\cdot
\prod_{\renewcommand{\arraystretch}{0.6}
\begin{array}{c}
\scriptstyle{V\in(\Gamma')^0\setminus\bp}\\
\scriptstyle{V\ne V_1}
\end{array}}[\mu(\overline\Gamma',h',V)]_y^-
\cdot\prod_{V\in(\overline\Gamma')^0_\infty}
[\mu(\overline\Gamma',h',V)]^*_y\ ,\label{e8a}\end{equation}
where
$$m=\mt(e),\quad\nu=\frac{\mu(\overline\Gamma',h',V_1)}{m},\quad y=z^2\ ,$$
and
$$\Psi^{(1)}_z(m,\nu)=\frac{2}{(z-z^{-1})^2(z+z^{-1})}$$
\begin{equation}\times\Bigg[\frac{z^{\nu m-1}-
z^{1-\nu m}}{z-z^{-1}}
-m\frac{z^{\nu m-m}-z^{m-\nu m}}{z^m-z^{-m}}
-\frac{z^{\nu m-\nu}-z^{\nu-\nu m}}{z^{\nu}-z^{-\nu}}
\Bigg]\ .
\label{e9a}\end{equation}

If $e=[v_1,v_2]=h'[V_1,V_2]$ (the image of a bounded edge), we let
$$\GS_y(\overline\Gamma',\bp,p,h')=\frac{1}{|\Aut(\overline\Gamma',\bp,p,h')|}\cdot
\Psi^{(2)}_z(m,\nu_1,\nu_2)\qquad\qquad\qquad\qquad$$
\begin{equation}\qquad\qquad\times\prod_{V\in\bp\cap(\Gamma')^0}
[\mu(\overline
\Gamma',h',V)]_y^+\cdot
\prod_{\renewcommand{\arraystretch}{0.6}
\begin{array}{c}
\scriptstyle{V\in(\Gamma')^0\setminus\bp}\\
\scriptstyle{V\not\in\{V_1,V_2\}}
\end{array}}[\mu(\overline\Gamma',h',V)]_y^-\cdot\prod_{V\in(\overline\Gamma')^0_\infty}
[\mu(\overline\Gamma',h',V)]^*_y\ ,\label{e8}\end{equation}
where
$$m=\mt(e),\quad\nu_i=\frac{\mu(\overline\Gamma',h',V_i)}{m},\ i=1,2,\quad y=z^2\ ,$$ and
$$\Psi^{(2)}_z(m,\nu_1,\nu_2)=\frac{1}{(z-z^{-1})^3(z+z^{-1})}$$
$$\times\Bigg[\frac{2(z^{\nu_2m}-z^{-\nu_2m})(z^{\nu_1m-1}-z^{1-\nu_1m})}{z-z^{-1}}
-\frac{2m(z^{\nu_2m}-z^{-\nu_2m})(z^{\nu_1m-m}-z^{m-\nu_1m})}{z^m-z^{-m}}$$
$$+(m-1)(z^{\nu_1m}-z^{-\nu_1m})(z^{\nu_2m}+z^{-\nu_2m})$$
\begin{equation}
-\frac{2(z^{\nu_2m}-z^{-\nu_2m})(z^{\nu_1m-\nu_1}-z^{\nu_1-\nu_1m})}{z^{\nu_1}-z^{-\nu_1}}
-\frac{2(z^{\nu_1m}-z^{-\nu_1m})(z^{\nu_2m-\nu_2}-z^{\nu_2-\nu_2m})}{z^{\nu_2}-z^{-\nu_2}}\Bigg]\ .
\label{e9}\end{equation}

In both cases, $\Aut(\overline\Gamma',\bp,p,h')$ is the automorphism group of $(\overline
\Gamma',\bp,p,h')$.

\begin{remark}\label{rr1} For the elliptic curves $(\overline\Gamma,\bp,h)$ counted in the right-hand side of
(\ref{ee1}), and rational curves $(\overline\Gamma',\bp,p,h')$ counted in the
right-hand side of (\ref{e8a}), the automorphism group is $(\Z/2)^r$, where $r$ is the number of collinear trivalent
vertices incident to a pair of ends of the same weight.
\end{remark}

We postpone study of functions $\Psi^{(1)}_y(m,\nu)$ and $\Psi^{(2)}_y(m,\nu_1,\nu_2)$ till Section
\ref{sec-com}, and here present only the following non-evident property.

\begin{lemma}\label{l7}
The function $\Psi^{(2)}_z(m,\nu_1,\nu_2)$ is symmetric with respect to $\nu_1$ and $\nu_2$, i.e.,
$\Psi^{(2)}_z(m,\nu_1,\nu_2)=\Psi^{(2)}_z(m,\nu_2,\nu_1)$.
\end{lemma}

{\bf Proof.} Since the two last summands in formula for $\Psi^{(2)}_z(m,\nu_1,\nu_2)$ are
symmetric with respect to
$\nu_1,\nu_2$, one only has to show that
$$\frac{2(z^{\nu_2m}-z^{-\nu_2m})(z^{\nu_1m-1}-z^{1-\nu_1m})}{z-z^{-1}}
-\frac{2m(z^{\nu_2m}-z^{-\nu_2m})(z^{\nu_1m-m}-z^{m-\nu_1m})}{z^m-z^{-m}}$$
$$+(m-1)(z^{\nu_1m}-z^{-\nu_1m})(z^{\nu_2m}+z^{-\nu_2m})$$
$$=\frac{2(z^{\nu_1m}-z^{-\nu_1m})(z^{\nu_2m-1}-z^{1-\nu_2m})}{z-z^{-1}}
-\frac{2m(z^{\nu_1m}-z^{-\nu_1m})(z^{\nu_2m-m}-z^{m-\nu_2m})}{z^m-z^{-m}}$$
$$+(m-1)(z^{\nu_2m}-z^{-\nu_2m})(z^{\nu_1m}+z^{-\nu_1m})$$
what can easily be done by a routine direct computation.
\proofend

Our main result is

\begin{theorem}\label{t1} Given a balanced, nondegenerate multiset
$\Delta\subset\Z^2\setminus\{0\}$ and integers $n_v>0$, $n_e\ge0$ such that
$2n_v+n_e=|\Delta|$, $n=n_v+n_e$, the expression
$$\GS_y(\Delta,1,(n_v,n_e),\bx):=\sum_{\renewcommand{\arraystretch}{0.6}
\begin{array}{c}
\scriptstyle{(\overline\Gamma,\bp,h)\in{\mathcal M}'_{1,(n_v,n_e)}
(\R^2,\Delta)}\\
\scriptstyle{h(\bp)=\bx}
\end{array}}\GS_y(\overline\Gamma,\bp,h)$$
\begin{equation}+\sum_{\renewcommand{\arraystretch}{0.6}
\begin{array}{c}
\scriptstyle{(\overline\Gamma',\bp,p,h')\in
{\mathcal M}''_{0,(n_v-1,n_e+1)}(\R^2,\Delta)}\\
\scriptstyle{h'(\bp)=\bx}
\end{array}}\GS_y(\overline\Gamma',\bp,p,h')\label{enew1}\end{equation}
does not depend on the choice of a generic $\bx\in\R^{2n}$.
\end{theorem}

\begin{remark}\label{r2}
In view of relation (\ref{elab}) and the absence of non-trivial automorphisms of
labeled plane tropical curves under consideration, the statement of Theorem \ref{t1} is equivalent to
the invariance of the expression
$$\GS^{\lab}_y(\Delta,1,(n_v,n_e),\bx):=\sum_{\renewcommand{\arraystretch}{0.6}
\begin{array}{c}
\scriptstyle{(\overline\Gamma^{\;\lab},\bp,h)\in({\mathcal M}^{\lab}_{1,(n_v,n_e)}
(\R^2,\Delta))'}\\
\scriptstyle{h(\bp)=\bx}
\end{array}}\GS_y(\overline\Gamma^{\;\lab},\bp,h)$$
$$+\sum_{\renewcommand{\arraystretch}{0.6}
\begin{array}{c}
\scriptstyle{((\overline\Gamma')^{\lab},\bp,p,h')\in
({\mathcal M}^{\lab}_{0,(n_v-1,n_e+1)}(\R^2,\Delta))''}\\
\scriptstyle{h'(\bp)=\bx}
\end{array}}\GS_y((\overline\Gamma')^{\lab},\bp,p,h')$$
with respect to the choice of a generic $\bx\in\R^{2n}$.
\end{remark}

\subsection{Elliptic broccoli curves}\label{broccoli}
Using the refined invariant, we are able to define elliptic broccoli curves.

Fix a balanced, nondegenerate multiset $\Delta\subset\Z^2\setminus\{0\}$
and integers $n_v>0$, $n_e\ge0$ such that
$2n_v+n_e=|\Delta|$, $n=n_v+n_e$. 
Fix also a generic point $\bx\in\R^{2n}$.

\begin{definition}\label{dd1}
(1) Suppose that $(\overline\Gamma,\bp,h)\in{\mathcal M}'_{1,(n_v,n_e)}
(\R^2,\Delta)$ and $h(\bp)=\bx$. Introduce the subgraph $\Gamma^{even}\subset\overline\Gamma$
containing all edges of even weight and their endpoints. We call the curve
$(\overline\Gamma,\bp,h)$ an {\bf elliptic broccoli curve}, if
each marked vertex of $\overline\Gamma$ is adjacent to at most one edge of even weight,
and, in each component of $\Gamma^{even}$, all but one univalent vertices belong to $\overline\Gamma_\infty^0\cup\bp$.

(2) Suppose that $(\overline\Gamma,\bp,h)\in{\mathcal M}''_{1,(n_v,n_e)}
(\R^2,\Delta)$ and $h(\bp)=\bx$, and let $\pi(\overline\Gamma,\bp,h)=(\overline\Gamma',
\bp,p,h')\in{\mathcal M}''_{0,(n_v-1,n_e+1)}
(\R^2,\Delta)$. Denote by $E_p$ the edge of $\overline\Gamma$ that includes the point $p$.
Introduce the subgraph $\overline\Gamma^{\;\prime,even}\subset\overline\Gamma'$ that contains all
edges of even weight except for the edge $E_p$, if
the weight of $E_p$ is even and $E_p$ has finite length, and contains all endpoints of these
edges. We call the curve $(\overline\Gamma,\bp,h)$ an {\bf elliptic broccoli curve} 
if
each marked vertex of $\overline\Gamma'$ is adjacent to at most one edge of even weight,
and, in each component of $\Gamma^{\;\prime,even}$, all but one univalent vertices belong to $\overline\Gamma_\infty^0\cup\bp$.
\end{definition}

Part (1) of Definition \ref{dd1} matches the definition of unoriented rational broccoli curves
\cite[Definition 3.1]{GMS}.

In Proposition \ref{cor1new}, Section \ref{char}, we show that the elliptic broccoli curves are characterized by the property that their refined broccoli invariant, evaluated at
$y=-1$, does not vanish, which, in particular, agrees with \cite[Corollary 3.9]{GS}.

\section{Proof of Theorem \ref{t1}}\label{sec3.7}

For the proof, we choose two generic configurations $\bx^{(0)},\bx^{(1)}\in X^{2n}
\subset\R^{2n}$, join them by a generic
path $\bx^{(t)}\in\R^{2n}$, and check the invariance as the path crosses top-dimensional
cells of $X^{2n-1}$ at their generic points. Clearly, we can
move points of the configuration one by one. Furthermore, by Remark \ref{r2}, we can work with
labeled curves. So, all curves in this section are labeled, but to simplify notations, we skip
the index ``$\lab$" everywhere.

For completeness, we shall consider all possible bifurcations so that as byproduct we prove the
invariance stated in Propositions \ref{p1} and \ref{p3}. In the sequel, we label
bifurcations as in the list presented in Lemma \ref{l6}(2).

Let $\bx^{(t^*)}$ be generic in an $(2n-1)$-dimensional cell of $X^{2n-1}$. Denote by $H_0$
the germ of this cell at $\bx^{(t^*)}$ and by $H_+,H_-\subset\R^{2n}$ the germs of the halfspaces
with common boundary $H_0$. Let $C^*=(\overline\Gamma,\bp,h)\in(\Ev_{(n_v,n_e)}^e)^{-1}
(\bx^{(t^*)})$
be as described in Lemma \ref{l6}(2), and let $F_0\subset\overline{\mathcal M}^{\;e}_{1,(n_v,
n_e)}(\R^2,\Delta)$ be the germ at $C^*$ of the $(2n-1)$-cell projecting by $\Ev^e_{(n_v,n_e}$
onto $H_0$. We shall analyze the $2n$-cells of $\overline{\mathcal M}^{\;e}_{1,(n_v,
n_e)}(\R^2,\Delta)$ attached to $F_0$, their projections onto $H_+,H_-$, and prove the
invariance of the expression $\GS_y(\Delta,1,(n_v,n_e),\bx)$.

\subsection{Degeneration of type (2i)}\label{sec-a} 
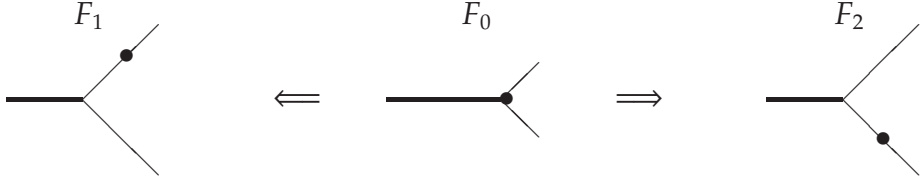
\begin{figure}
\setlength{\unitlength}{1cm}
\begin{picture}(14,3)(0,0)
\thinlines

\put(2,1){\line(1,-1){1}}\put(2,1){\line(1,1){1}}
\put(2.47,1.47){$\bullet$}

\put(12,1){\line(1,-1){1}}\put(12,1){\line(1,1){1}}
\put(12.42,0.38){$\bullet$}

\put(7.5,1){\line(1,-1){0.5}}\put(7.5,1){\line(1,1){0.5}}
\put(7.47,0.90){$\bullet$}

\put(1.9,2){$F_1$}\put(7,2){$F_0$}\put(11.9,2){$F_2$}
\put(4.5,0.9){$\Longleftarrow$}\put(9,0.9){$\Longrightarrow$}

\linethickness{0.5mm}

\put(1,1){\line(1,0){1}}\put(6,1){\line(1,0){1.5}}\put(11,1){\line(1,0){1}}

\end{picture}
\caption{Degeneration of type (2i)}\label{fig2}
\end{figure}

Let $C^*$ be as in
Lemma \ref{l6}(2i).
Notice that the
set $\Gamma^0_{sing}$ of the marked trivalent non-collinear vertices of $\Gamma$ belonging
to the closure
of the bounded component of $\overline\Gamma\setminus\bp$, is nonempty. Each marked point $p\in\Gamma^0_{sing}$
is incident to two unbounded and one bounded component of $\overline\Gamma\setminus\bp$
(see Figure \ref{fig2}, where the bounded component of $\overline\Gamma\setminus\bp$ is
shown fat), and
the cell $F_0$ lies on the boundary of exactly two $2n$-dimensional cells $F_+,F_-$ of
$\overline{\mathcal M}^{\;e}_{1,(n_v,n_e)}(\R^2,\Delta)$ along which the marked point $p$
moves towards one of
the unbounded components of $\overline\Gamma\setminus\bp$. Clearly, $\Ev^e_{(n_v-1,n_e+1)}$
takes $F_+$, $F_-$ onto the
germs $H_+$, $H_-$ of the half spaces in $\R^{2n}$ with $H_0$ as a common boundary
(see Figure \ref{fig2}).
Clearly, for any point $p\in\Gamma^0_{sing}$, such a bifurcation does not affect the value of
$\GS_y(\Delta,1,(n_v,n_e),\bx^{(t)})$, $t\in(\R,t^*)$.

\subsection{Degeneration of type (2ii)}\label{sec-b}
\begin{figure}
\setlength{\unitlength}{1cm}
\begin{picture}(14,7)(0.5,0)
\thinlines

\put(1.5,3.6){\line(1,1){1}}\put(1.5,3.6){\line(1,-1){1}}
\put(1.5,3.6){\line(-1,2){0.5}}\put(1.5,3.6){\line(-1,-2){0.5}}

\put(5.2,5.7){\line(1,-1){0.8}}\put(5.2,5.7){\line(-1,-2){0.4}}
\put(5.2,6.1){\line(1,1){0.8}}\put(5.2,6.1){\line(-1,2){0.4}}
\put(5.2,5.7){\line(0,1){0.4}}

\put(5,3.6){\line(-1,-2){0.4}}\put(5,3.6){\line(-1,2){0.4}}
\put(5.4,3.6){\line(1,1){0.6}}\put(5.4,3.6){\line(1,-1){0.6}}
\put(5,3.6){\line(1,0){0.4}}

\put(5,1.5){\line(-1,-2){0.3}}\put(5,1.9){\line(-1,2){0.3}}
\put(5,1.5){\line(1,1){0.7}}\put(5,1.9){\line(1,-1){0.7}}
\put(5,1.5){\line(0,1){0.4}}

\put(3,3.5){$\Longrightarrow$}
\put(3,4.8){$\Nearrow$}\put(3,2.3){$\Searrow$}

\put(3.3,0){(a)}\put(10.8,0){(b)}

\put(7,5.7){\line(2,1){2}}\put(7,5.7){\line(2,-1){2}}
\put(9,6.7){\line(1,-1){1}}\put(9,4.7){\line(1,1){1}}
\put(7,5.7){\line(1,0){3}}

\put(7,3.6){\line(2,1){2}}\put(7,3.6){\line(2,-1){2}}
\put(9,4.6){\line(1,-1){1}}\put(9,2.6){\line(1,1){1}}
\put(9,2.6){\line(0,1){2}}

\put(7,1.5){\line(2,1){2}}\put(7,1.5){\line(2,-1){2}}
\put(9,2.5){\line(1,-1){1}}\put(9,0.5){\line(1,1){1}}
\put(9,2.5){\line(-1,-1){1}}\put(9,0.5){\line(-1,1){1}}
\put(7,1.5){\line(1,0){1}}

\put(11.7,3.6){\line(2,1){2}}\put(11.7,3.6){\line(2,-1){2}}
\put(13.7,4.6){\line(1,-1){1}}\put(13.7,2.6){\line(1,1){1}}

\put(10.5,3.5){$\Longleftarrow$}
\put(10.7,4.8){$\Nwarrow$}\put(10.7,2.3){$\Swarrow$}

\put(0.5,3.5){$F_0$}\put(4.2,5.6){$F_1$}\put(4.2,3.5){$F_2$}\put(4.2,1.6){$F_3$}

\put(8.6,5.9){$\mu_1$}\put(8.6,5.3){$\mu_2$}
\put(8.1,3.5){$\mu_3$}\put(9.2,3.5){$\mu_4$}
\put(7.5,2.1){$\mu_5$}\put(7.5,0.8){$\mu_6$}

\end{picture}
\caption{Degeneration of type (2ii)}\label{fig3}
\end{figure}
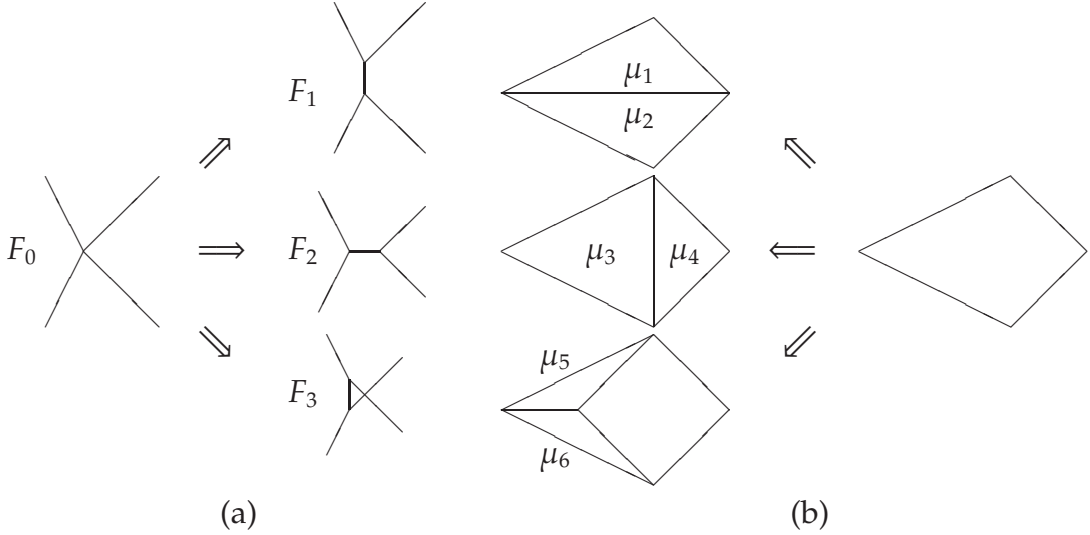

Let $C^*$ be as in Lemma \ref{l6}(2ii). The curve
$C^*$ admits three types of deformation according to three types of
splitting of a four-valent vertex into a pair of trivalent vertices as shown in Figure
\ref{fig3}, where the right part of the
figure exhibits fragments of the dual subdivision of the Newton polygon (cf.
\cite[Page 172]{GM}). That is, $F_0$ lies in the boundary of three $2n$-cells $F_1,F_2,F_3
\subset \overline{\mathcal M}^{\;e}_{1,(n_v,n_e)}(\R^2,\Delta)$. We claim that
$\Ev^e_{(n_v,n_e)}$ takes
$F_1$ onto the germ of the halfspace $H_+$, and takes $F_2,F_3$ onto the germ of the
halfspace $H_-$. Indeed, this holds in the situation considered in
\cite{GM} when $n_v=0$, and it corresponds to the local invariance of the
count of tropical curves with Mikhalkin's weights, which in turn reduces here to the
elementary geometric relation (see Figure \ref{fig3}(b))
$$\mu_1\mu_2=\mu_3\mu_4+\mu_5\mu_6\ ,$$
where $\mu_i$'s, in fact, are equal to the lattice areas of the dual triangles.
The same holds in our case,
if we replace each marked vertex by a couple of close marked points on edges (see, for example
such a replacement in the case of a collinear marked vertex in Figure \ref{fig3}(c)).
Thus, the invariance of $\GS_y(\Delta,1,(n_v,n_e),\bx^{(t)})$,
$t\in(\R,t^*)$, reduces to the relation
$$\varphi(y)(z^{\mu_1}-z^{-\mu_1})(z^{\mu_2}-z^{-\mu_2})=\varphi(y)(z^{\mu_3}-z^{-\mu_3})
(z^{\mu_4}-z^{\mu_4})+\varphi(y)(z^{\mu_5}-z^{-\mu_5})(z^{\mu_6}-
z^{-\mu_6})$$ with some expression $\varphi(y)$ and $y=z^2$. So, we have to show
$$(z^{\mu_1}-z^{-\mu_1})(z^{\mu_2}-z^{-\mu_2})=(z^{\mu_3}-z^{-\mu_3})
(z^{\mu_4}-z^{\mu_4})+(z^{\mu_5}-z^{-\mu_5})(z^{\mu_6}-
z^{-\mu_6}),$$ which immediately follows from the following elementary geometric observations:
\begin{equation}\mu_1+\mu_2=\mu_3+\mu_4,\quad\mu_1-\mu_2=\mu_6-\mu_5,\quad
\mu_3-\mu_4=\mu_5+\mu_6\ .\label{e3}\end{equation}

\subsection{Degeneration of type (2iii)}\label{sec3.3}
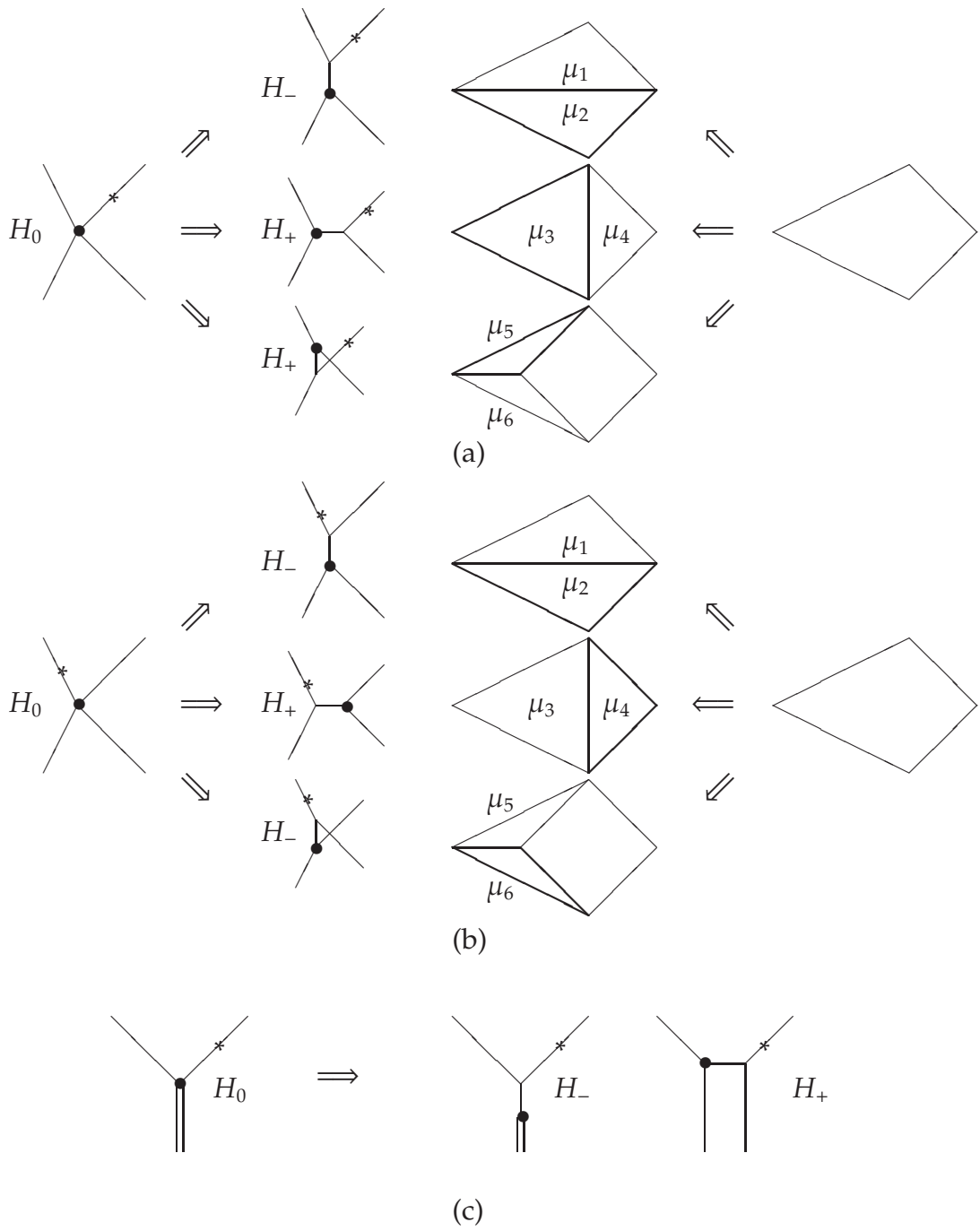
\begin{figure}
\setlength{\unitlength}{1cm}
\begin{picture}(14,18)(0.5,-4)
\thinlines

\put(1.5,3.6){\line(1,1){1}}\put(1.5,3.6){\line(1,-1){1}}
\put(1.5,3.6){\line(-1,2){0.5}}\put(1.5,3.6){\line(-1,-2){0.5}}

\put(5.2,5.7){\line(1,-1){0.8}}\put(5.2,5.7){\line(-1,-2){0.4}}
\put(5.2,6.1){\line(1,1){0.8}}\put(5.2,6.1){\line(-1,2){0.4}}
\put(5.2,5.7){\line(0,1){0.4}}

\put(5,3.6){\line(-1,-2){0.4}}\put(5,3.6){\line(-1,2){0.4}}
\put(5.4,3.6){\line(1,1){0.6}}\put(5.4,3.6){\line(1,-1){0.6}}
\put(5,3.6){\line(1,0){0.4}}

\put(5,1.5){\line(-1,-2){0.3}}\put(5,1.9){\line(-1,2){0.3}}
\put(5,1.5){\line(1,1){0.7}}\put(5,1.9){\line(1,-1){0.7}}
\put(5,1.5){\line(0,1){0.4}}

\put(3,3.5){$\Longrightarrow$}
\put(3,4.8){$\Nearrow$}\put(3,2.3){$\Searrow$}

\put(7,0){(b)}\put(7,7.2){(a)}\put(7,-4){(c)}

\put(7,5.7){\line(2,1){2}}
\put(9,6.7){\line(1,-1){1}}

\put(7,3.6){\line(2,1){2}}\put(7,3.6){\line(2,-1){2}}

\put(7,1.5){\line(2,1){2}}
\put(9,2.5){\line(1,-1){1}}\put(9,0.5){\line(1,1){1}}
\put(9,2.5){\line(-1,-1){1}}

\put(11.7,3.6){\line(2,1){2}}\put(11.7,3.6){\line(2,-1){2}}
\put(13.7,4.6){\line(1,-1){1}}\put(13.7,2.6){\line(1,1){1}}

\put(10.5,3.5){$\Longleftarrow$}
\put(10.7,4.8){$\Nwarrow$}\put(10.7,2.3){$\Swarrow$}

\put(0.5,3.5){$H_0$}\put(4.2,5.6){$H_-$}\put(4.2,3.5){$H_+$}\put(4.2,1.6){$H_-$}

\put(8.6,5.9){$\mu_1$}\put(8.6,5.3){$\mu_2$}
\put(8.1,3.5){$\mu_3$}\put(9.2,3.5){$\mu_4$}
\put(7.5,2.1){$\mu_5$}\put(7.5,0.8){$\mu_6$}

\put(1.5,10.6){\line(1,1){1}}\put(1.5,10.6){\line(1,-1){1}}
\put(1.5,10.6){\line(-1,2){0.5}}\put(1.5,10.6){\line(-1,-2){0.5}}

\put(5.2,12.7){\line(1,-1){0.8}}\put(5.2,12.7){\line(-1,-2){0.4}}
\put(5.2,13.1){\line(1,1){0.8}}\put(5.2,13.1){\line(-1,2){0.4}}
\put(5.2,12.7){\line(0,1){0.4}}

\put(5,10.6){\line(-1,-2){0.4}}\put(5,10.6){\line(-1,2){0.4}}
\put(5.4,10.6){\line(1,1){0.6}}\put(5.4,10.6){\line(1,-1){0.6}}
\put(5,10.6){\line(1,0){0.4}}

\put(5,8.5){\line(-1,-2){0.3}}\put(5,8.9){\line(-1,2){0.3}}
\put(5,8.5){\line(1,1){0.7}}\put(5,8.9){\line(1,-1){0.7}}
\put(5,8.5){\line(0,1){0.4}}

\put(3,10.5){$\Longrightarrow$}
\put(3,11.8){$\Nearrow$}\put(3,9.3){$\Searrow$}

\put(7,12.7){\line(2,1){2}}
\put(9,13.7){\line(1,-1){1}}

\put(9,11.6){\line(1,-1){1}}\put(9,9.6){\line(1,1){1}}

\put(7,8.5){\line(2,-1){2}}
\put(9,9.5){\line(1,-1){1}}\put(9,7.5){\line(1,1){1}}
\put(9,7.5){\line(-1,1){1}}

\put(11.7,10.6){\line(2,1){2}}\put(11.7,10.6){\line(2,-1){2}}
\put(13.7,11.6){\line(1,-1){1}}\put(13.7,9.6){\line(1,1){1}}

\put(10.5,10.5){$\Longleftarrow$}
\put(10.7,11.8){$\Nwarrow$}\put(10.7,9.3){$\Swarrow$}

\put(0.5,10.5){$H_0$}\put(4.2,12.6){$H_-$}\put(4.2,10.5){$H_+$}\put(4.2,8.6){$H_+$}

\put(8.6,12.9){$\mu_1$}\put(8.6,12.3){$\mu_2$}
\put(8.1,10.5){$\mu_3$}\put(9.2,10.5){$\mu_4$}
\put(7.5,9.1){$\mu_5$}\put(7.5,7.8){$\mu_6$}

\put(1.42,3.5){$\bullet$}\put(5.1,5.54){$\bullet$}
\put(5.35,3.47){$\bullet$}\put(4.9,1.37){$\bullet$}

\put(1.42,10.5){$\bullet$}\put(5.1,12.54){$\bullet$}
\put(4.9,10.47){$\bullet$}\put(4.9,8.77){$\bullet$}

\put(1.2,4){$*$}\put(5,6.3){$*$}\put(4.8,3.8){$*$}\put(4.8,2.1){$*$}
\put(1.95,11){$*$}\put(5.5,13.35){$*$}\put(5.7,10.8){$*$}\put(5.4,8.85){$*$}

\thicklines

\put(7,5.7){\line(2,-1){2}}\put(9,4.7){\line(1,1){1}}
\put(7,5.7){\line(1,0){3}}\put(9,4.6){\line(1,-1){1}}\put(9,2.6){\line(1,1){1}}
\put(9,2.6){\line(0,1){2}}\put(7,1.5){\line(2,-1){2}}
\put(9,0.5){\line(-1,1){1}}
\put(7,1.5){\line(1,0){1}}\put(7,12.7){\line(2,-1){2}}\put(9,11.7){\line(1,1){1}}
\put(7,12.7){\line(1,0){3}}\put(7,10.6){\line(2,1){2}}\put(7,10.6){\line(2,-1){2}}
\put(9,9.6){\line(0,1){2}}\put(7,8.5){\line(2,1){2}}\put(9,9.5){\line(-1,-1){1}}
\put(7,8.5){\line(1,0){1}}

\thinlines

\put(3,-2){\line(-1,1){1}}\put(3,-2){\line(1,1){1}}
\put(2.95,-2){\line(0,-1){1}}\put(3.05,-2){\line(0,-1){1}}

\put(5,-2){$\Longrightarrow$}

\put(8,-2){\line(-1,1){1}}\put(8,-2){\line(1,1){1}}\put(8,-2){\line(0,-1){0.5}}
\put(7.95,-2.5){\line(0,-1){0.5}}\put(8.05,-2.5){\line(0,-1){0.5}}

\put(10.7,-1.7){\line(-1,1){0.7}}\put(11.3,-1.7){\line(1,1){0.7}}\put(11.3,-1.7){\line(-1,0){0.6}}
\put(10.7,-1.7){\line(0,-1){1.3}}\put(11.3,-1.7){\line(0,-1){1.3}}

\put(2.9,-2.1){$\bullet$}\put(7.93,-2.6){$\bullet$}\put(10.6,-1.8){$\bullet$}
\put(3.5,-1.55){$*$}\put(8.5,-1.55){$*$}\put(11.5,-1.55){$*$}
\put(3.5,-2.2){$H_0$}\put(8.5,-2.2){$H_-$}\put(12,-2.2){$H_+$}

\end{picture}
\caption{Degeneration of type (2iii)}\label{fig4}
\end{figure}

Let $C^*=(\overline\Gamma,\bp,h)$ be as in Lemma \ref{l6}(2iii).
As above (cf. \cite[Page 24]{GS}) we
derive that $F_0$ lies in the boundary of three $2n$-cells $F_1,F_2,F_3$ of
$\overline{\mathcal M}^{\;e}_{1,(n_v,n_e)}(\R^2,\Delta)$ according to three types of
splitting of the four-valent vertex into a
pair of  trivalent vertices, shown in Figure \ref{fig4}. We have to study two cases
according as the edge of a bounded component of $\overline\Gamma\setminus\bp$
is dual to a side of
the parallelogram inscribed into the quadrangle
or not
(see Figures \ref{fig4}(a,b), where the edge belonging to the bounded component of $\overline
\Gamma\setminus
\bp$ and its dual are labeled by asterisk, and the triangles dual to the marked trivalent
vertices are shown
by fat lines).
Assuming that along the germ of the path $\left(\bx^{(t)}\right)_{t\in[0,1]}$
at $\bx^{(t^*)}$, the image of the marked point at the four-valent vertex of $C^*$ is moving
and the rest of $\bp$ is fixed, we decide which cells $F_i$, $i=1,2,3$, project onto $H_+$
or onto $H_-$ according as the moving marked vertex belongs to the halfplane $\R^2_+$ or
$\R^2_-$. In the notation of Section \ref{sec-b} for the Mikhalkin's weights of the
trivalent vertices
(see Figures \ref{fig4}(a,b)), we have the following additional geometric relations
\begin{equation}\begin{cases}\mu_3=\mu_1+\mu_5,\quad&\text{in Figure \ref{fig4}(a)},\\
\mu_1=\mu_4+\mu_6,\quad&\text{in Figure \ref{fig4}(b)}.\end{cases}\label{e4}\end{equation}
Similarly to Section \ref{sec-b}, the invariance of
$\GS_y(\Delta,1,(n_v,n_e),\bx^{(t)})$, $t\in(\R,t^*)$, reduces to the relation
$$(z^{\mu_3}-z^{-\mu_3})(z^{\mu_4}+z^{-\mu_4})=
(z^{\mu_1}-z^{-\mu_1})(z^{\mu_2}+z^{-\mu_2})+
(z^{\mu_5}-z^{-\mu_5})(z^{\mu_6}+z^{-\mu_6})$$
in case of Figure \ref{fig4}(a), and the the relation
$$(z^{\mu_1}-z^{-\mu_1})(z^{\mu_2}+z^{-\mu_2})=
(z^{\mu_3}+z^{-\mu_3})(z^{\mu_4}-z^{-\mu_4})+
(z^{\mu_5}+z^{-\mu_5})(z^{\mu_6}-z^{-\mu_6})$$
in case of Figure \ref{fig4}(b). Both relations
immediately follow from (\ref{e3}) and (\ref{e4}).


\subsection{Degeneration of type (2iv)}
\begin{figure}
\setlength{\unitlength}{1cm}
\begin{picture}(14,15)(0,1)
\thinlines

\put(1,15){\line(1,1){1}}\put(1,15){\line(1,-1){1}}
\put(1,15){\line(-1,2){0.5}}\put(1,15){\line(-1,-2){0.5}}

\put(5,15){\line(1,1){1}}\put(5,14.6){\line(1,-1){1}}
\put(5,15){\line(-1,2){0.5}}\put(5,14.6){\line(-1,-2){0.5}}
\put(5,14.6){\line(0,1){0.4}}

\put(7.5,15){\line(1,1){1}}\put(7.5,14.6){\line(1,-1){1}}
\put(7.5,15){\line(-1,2){0.5}}\put(7.5,14.6){\line(-1,-2){0.5}}
\put(7.5,14.6){\line(0,1){0.4}}

\put(10,14.6){\line(1,1){1}}\put(10,15){\line(1,-1){1}}
\put(10,15){\line(-1,2){0.5}}\put(10,14.6){\line(-1,-2){0.5}}
\put(10,14.6){\line(0,1){0.4}}

\put(12.5,14.6){\line(1,1){1}}\put(12.5,15){\line(1,-1){1}}
\put(12.5,15){\line(-1,2){0.5}}\put(12.5,14.6){\line(-1,-2){0.5}}
\put(12.5,14.6){\line(0,1){0.4}}

\put(0.7,15.4){$*$}\put(0.7,14.45){$*$}
\put(4.7,15.4){$*$}\put(4.7,14.05){$*$}
\put(7.2,15.4){$*$}\put(7.2,14.05){$*$}
\put(9.7,15.4){$*$}\put(9.7,14.05){$*$}
\put(12.2,15.4){$*$}\put(12.2,14.05){$*$}

\put(0.9,14.9){$\bullet$}\put(4.9,14.9){$\bullet$}
\put(7.4,14.5){$\bullet$}\put(9.9,14.9){$\bullet$}
\put(12.4,14.5){$\bullet$}

\put(1.3,14.9){$H_0$}\put(5.3,14.7){$H_+$}\put(7.8,14.7){$H_-$}
\put(10.6,14.7){$H_+$}\put(13.1,14.7){$H_-$}

\put(0.9,13){(a)}\put(4.8,13){(b$_1$)}\put(7.3,13){(b$_2$)}
\put(9.8,13){(b$_3$)}\put(12.3,13){(b$_4$)}

\put(0,11){\line(2,1){2}}\put(0,11){\line(2,-1){2}}
\put(2,12){\line(1,-1){1}}\put(2,10){\line(1,1){1}}

\dashline{0.2}(1,11)(2,12)\dashline{0.2}(1,11)(2,10)

\put(0.9,11.35){$*$}\put(0.9,10.4){$*$}

\put(5,11){\line(2,1){2}}\put(5,11){\line(2,-1){2}}
\put(7,12){\line(1,-1){1}}\put(7,10){\line(1,1){1}}
\put(5,11){\line(1,0){3}}

\put(9,11){\line(2,1){2}}\put(9,11){\line(2,-1){2}}
\put(11,12){\line(1,-1){1}}\put(11,10){\line(1,1){1}}
\put(9,11){\line(1,0){1}}\put(10,11){\line(1,1){1}}\put(10,11){\line(1,-1){1}}

\put(6.7,11.2){$\mu_1$}\put(6.7,10.6){$\mu_2$}
\put(9.5,11.6){$\mu_5$}\put(9.5,10.3){$\mu_6$}

\put(1.4,9.3){(a')}\put(8.3,9.3){(b')}

\put(1,7.2){\line(1,1){1}}\put(1,7.2){\line(1,-1){1}}
\put(1,7.2){\line(-1,2){0.5}}\put(1,7.2){\line(-1,-2){0.5}}

\put(5.3,7.2){\line(1,1){1}}\put(5.3,7.2){\line(1,-1){1}}
\put(5,7.2){\line(-1,2){0.5}}\put(5,7.2){\line(-1,-2){0.5}}
\put(5,7.2){\line(1,0){0.3}}

\put(7.8,7.2){\line(1,1){1}}\put(7.8,7.2){\line(1,-1){1}}
\put(7.5,7.2){\line(-1,2){0.5}}\put(7.5,7.2){\line(-1,-2){0.5}}
\put(7.5,7.2){\line(1,0){0.3}}

\put(10,6.8){\line(1,1){1}}\put(10,7.2){\line(1,-1){1}}
\put(10,7.2){\line(-1,2){0.5}}\put(10,6.8){\line(-1,-2){0.5}}
\put(10,6.8){\line(0,1){0.4}}

\put(12.5,6.8){\line(1,1){1}}\put(12.5,7.2){\line(1,-1){1}}
\put(12.5,7.2){\line(-1,2){0.5}}\put(12.5,6.8){\line(-1,-2){0.5}}
\put(12.5,6.8){\line(0,1){0.4}}

\put(1.4,7.6){$*$}\put(0.7,7.6){$*$}
\put(5.7,7.6){$*$}\put(4.7,7.6){$*$}
\put(8.2,7.6){$*$}\put(7.2,7.6){$*$}
\put(10.6,7.3){$*$}\put(9.7,7.6){$*$}
\put(13.1,7.3){$*$}\put(12.2,7.6){$*$}

\put(0.9,7.1){$\bullet$}\put(4.9,7.1){$\bullet$}
\put(7.7,7.1){$\bullet$}\put(9.9,7.1){$\bullet$}
\put(12.4,6.7){$\bullet$}

\put(1.3,7.1){$H_0$}\put(5.7,7.1){$H_+$}\put(8.2,7.1){$H_-$}
\put(10.6,6.9){$H_+$}\put(13.1,6.9){$H_+$}

\put(0.9,5.2){(c)}\put(4.8,5.2){(d$_1$)}\put(7.3,5.2){(d$_2$)}
\put(9.8,5.2){(d$_3$)}\put(12.3,5.2){(d$_4$)}

\put(0,3.2){\line(2,1){2}}\put(0,3.2){\line(2,-1){2}}
\put(2,4.2){\line(1,-1){1}}\put(2,2.2){\line(1,1){1}}

\dashline{0.2}(1,3.2)(2,4.2)\dashline{0.2}(1,3.2)(2,2.2)

\put(0.9,3.55){$*$}\put(2.45,3.55){$*$}

\put(5,3.2){\line(2,1){2}}\put(5,3.2){\line(2,-1){2}}
\put(7,4.2){\line(1,-1){1}}\put(7,2.2){\line(1,1){1}}
\put(7,2.2){\line(0,1){2}}

\put(9,3.2){\line(2,1){2}}\put(9,3.2){\line(2,-1){2}}
\put(11,4.2){\line(1,-1){1}}\put(11,2.2){\line(1,1){1}}
\put(9,3.2){\line(1,0){1}}\put(10,3.2){\line(1,1){1}}\put(10,3.2){\line(1,-1){1}}

\put(6.3,3.1){$\mu_3$}\put(7.2,3.1){$\mu_4$}
\put(9.5,3.8){$\mu_5$}\put(9.5,2.5){$\mu_6$}

\put(1.4,1.5){(c')}\put(8.3,1.5){(d')}

\end{picture}
\caption{Degeneration of type (2iv)}\label{fig5}
\end{figure}
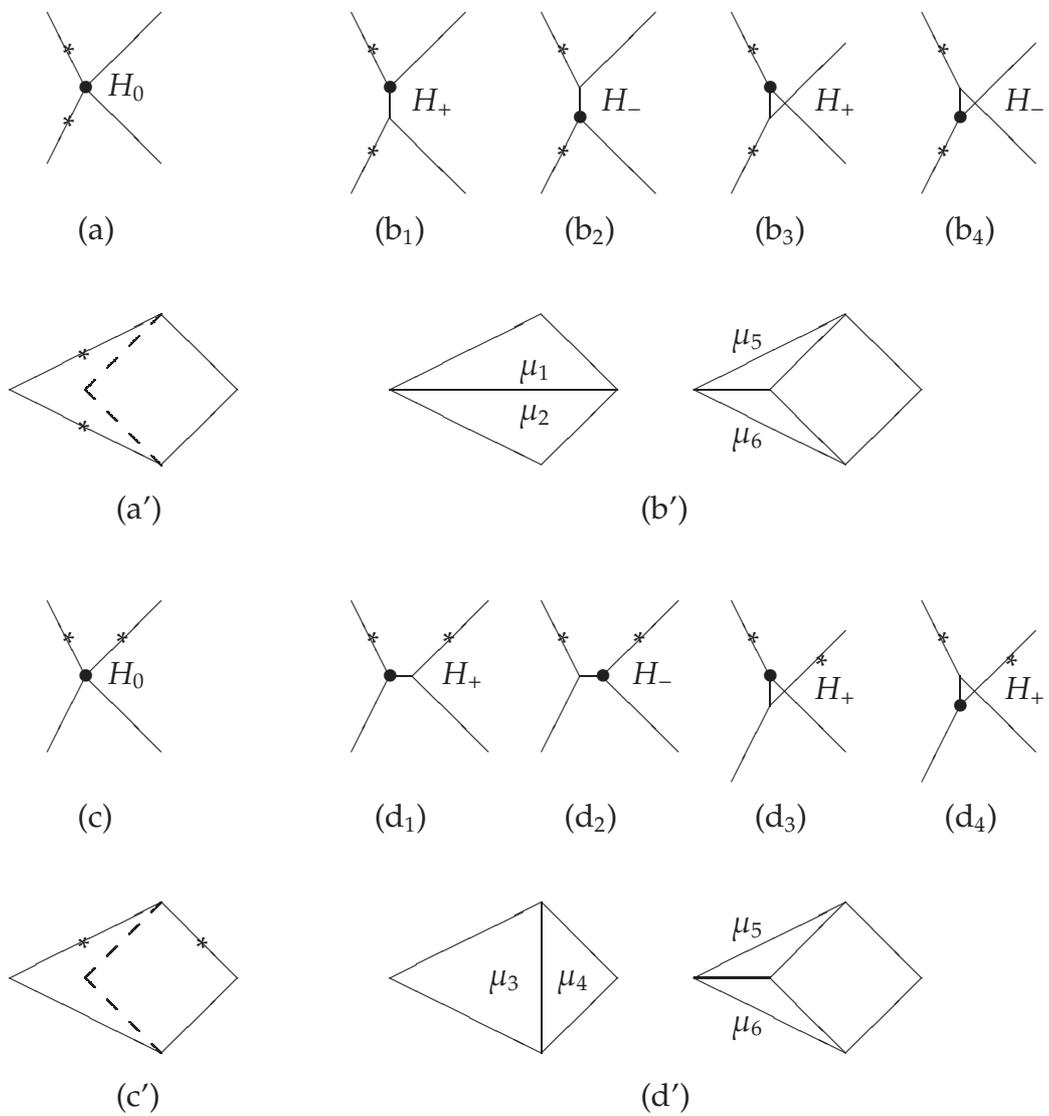

Let $C^*=(\overline\Gamma,\bp,h)$ be as in Lemma \ref{l6}(2iv)
(see Figure \ref{fig5}(a,c)). It admits four different deformations shown in Figures \ref{fig5} (b$_1$-b$_4$) and
\ref{fig5}(d$_1$-d$_4$). It is easy to see that, if both the edges of the bounded component of
$\overline
\Gamma\setminus\bp$ are dual to the sides of the inscribed parallelogram or both are
not (cf. Figures
\ref{fig5}(a,b$_1$-b$_4$) and \ref{fig5}(a',b')), two $2n$-cells of
$\overline{\mathcal M}^{\;e}_{1,
(n_v,n_e)}(\R^2,\Delta)$ attached to $F_0$ are projected onto $H^+$ and the other two onto $H_-$.
In turn, if one edge of the bounded component of $\overline
\Gamma\setminus\bp$ is dual to a side of the inscribed
parallelogram and the other is not (cf. Figures \ref{fig5}(c,d$_1$-d$_4$) and \ref{fig5}(c',d')),
then
three $2n$-cells of $\overline{\mathcal M}^{\;e}_{1,
(n_v,n_e)}(\R^2,\Delta)$ attached to $F_0$ are projected onto $H^+$ and the remaining
one onto $H_-$.
Thus, the invariance of $\GS_y(\Delta,1,(n_v,n_e),\bx^{(t)})$, $t\in(\R,t^*)$, reduces
\begin{itemize}\item in the former case to the equality (cf. Figure \ref{fig5}(b))
\begin{eqnarray}&(z^{\mu_1}+z^{-\mu_1})(z^{\mu_2}-z^{-\mu_2})+(z^{\mu_5}+
z^{-\mu_5})(z^{\mu_6}-z^{-\mu_6})\nonumber\\
&\qquad\qquad=(z^{\mu_1}-z^{-\mu_1})(z^{\mu_2}+z^{-\mu_2})+(z^{\mu_5}-z^{-\mu_5})
(z^{\mu_6}+z^{-\mu_6})\ ,\nonumber\end{eqnarray}
or, equivalently,
$$z^{\mu_2-\mu_1}-z^{\mu_1-\mu_2}=z^{\mu_5-\mu_6}-z^{\mu_6-\mu_5}\ ,$$ which
immediately follows from the second
relation in (\ref{e3});
\item in the latter case to the equality (cf. Figure \ref{fig5}(d))
\begin{eqnarray}&(z^{\mu_3}+z^{-\mu_3})(z^{\mu_4}-z^{-\mu_4})+
(z^{\mu_5}+z^{-\mu_5})(z^{\mu_6}-z^{-\mu_6})\nonumber\\
&\qquad\qquad+
(z^{\mu_5}-z^{-\mu_5})(z^{\mu_6}+z^{-\mu_6})=(z^{\mu_3}-z^{-\mu_3})(z^{\mu_4}+z^{-\mu_4})\ ,\nonumber\end{eqnarray}
or, equivalently,
$$z^{\mu_3-\mu_4}-z^{\mu_4-\mu_3}=z^{\mu_5+\mu_6}-z^{-\mu_5-\mu_6}\ ,$$
which immediately follows from the first relation in (\ref{e3}).
\end{itemize}

\subsection{Degeneration of type (2v)}
Let $C^*$ be as in Lemma \ref{l6}(2v) (see Figure \ref{fig6}(a-c)). The collinear cycle
contains either one or two marked points,
and we can assume that the image of one of them
in $\bx^{(t)}$ moves transversally to the image of the collinear cycle,
while the rest of marked points have fixed images along the path $\bx^{(t)}$,
$t\in(\R,t^*)$. The curve $C^*$ admits exactly two
deformations into generic elements of $\overline{\mathcal M}^{\;e}_{1,(n_v,n_e)}
(\R^2,\Delta)$, and
the corresponding $2n$-cells $F_1,F_2$, attached to $F_0$, project to $H_+,H_-$,
respectively
(see Figure \ref{fig6}(a,b,c), cf. \cite[Page 175, case (d)]{GM}). The invariance of
$\GS_y(\Delta,1,(n_v,n_e),\bx^{(t)})$, $t\in(\R,t^*)$, is evident, since the weights
of the edges
converging to the edges of the collinear cycle and
the Mikhalkin's weights of the trivalent vertices converging to the four-valent vertices
are respectively conserved (see Figure \ref{fig6}(d)).

\begin{figure}
\setlength{\unitlength}{1cm}
\begin{picture}(14,14.5)(0,0)
\thinlines

\put(0,13){\line(1,1){0.5}}\put(0.5,13.5){\line(1,0){2}}\put(0.5,13.5){\line(0,1){0.5}}
\put(0.5,14){\line(-1,1){0.5}}\put(0.5,14){\line(1,0){2}}\put(2.5,13.5){\line(1,-1){0.5}}
\put(2.5,13.5){\line(0,1){0.5}}\put(2.5,14){\line(1,1){0.5}}

\put(12,13){\line(1,1){0.5}}\put(12.5,13.5){\line(1,0){2}}\put(12.5,13.5){\line(0,1){0.5}}
\put(12.5,14){\line(-1,1){0.5}}\put(12.5,14){\line(1,0){2}}\put(14.5,13.5){\line(1,-1){0.5}}
\put(14.5,13.5){\line(0,1){0.5}}\put(14.5,14){\line(1,1){0.5}}

\put(0,9.5){\line(1,1){0.5}}\put(0.5,10){\line(1,0){2}}\put(0.5,10){\line(0,1){0.5}}
\put(0.5,10.5){\line(-1,1){0.5}}\put(0.5,10.5){\line(1,0){2}}\put(2.5,10){\line(1,-1){0.5}}
\put(2.5,10){\line(0,1){0.5}}\put(2.5,10.5){\line(1,1){0.5}}

\put(12,9.5){\line(1,1){0.5}}\put(12.5,10){\line(1,0){2}}\put(12.5,10){\line(0,1){0.5}}
\put(12.5,10.5){\line(-1,1){0.5}}\put(12.5,10.5){\line(1,0){2}}\put(14.5,10){\line(1,-1){0.5}}
\put(14.5,10){\line(0,1){0.5}}\put(14.5,10.5){\line(1,1){0.5}}

\put(0,6){\line(1,1){0.5}}\put(0.5,6.5){\line(1,0){2}}\put(0.5,6.5){\line(0,1){0.5}}
\put(0.5,7){\line(-1,1){0.5}}\put(0.5,7){\line(1,0){2}}\put(2.5,6.5){\line(1,-1){0.5}}
\put(2.5,6.5){\line(0,1){0.5}}\put(2.5,7){\line(1,1){0.5}}

\put(12,6){\line(1,1){0.5}}\put(12.5,6.5){\line(1,0){2}}\put(12.5,6.5){\line(0,1){0.5}}
\put(12.5,7){\line(-1,1){0.5}}\put(12.5,7){\line(1,0){2}}\put(14.5,6.5){\line(1,-1){0.5}}
\put(14.5,6.5){\line(0,1){0.5}}\put(14.5,7){\line(1,1){0.5}}

\put(6,13.3){\line(1,1){0.5}}\put(6.5,13.8){\line(-1,1){0.5}}
\put(6.5,13.75){\line(1,0){2}}\put(6.5,13.9){\line(1,0){2}}
\put(8.5,13.8){\line(1,1){0.5}}\put(8.5,13.8){\line(1,-1){0.5}}

\put(6,9.8){\line(1,1){0.5}}\put(6.5,10.3){\line(-1,1){0.5}}
\put(6.5,10.25){\line(1,0){2}}\put(6.5,10.4){\line(1,0){2}}
\put(8.5,10.3){\line(1,1){0.5}}\put(8.5,10.3){\line(1,-1){0.5}}

\put(6,6.3){\line(1,1){0.5}}\put(6.5,6.8){\line(-1,1){0.5}}
\put(6.5,6.75){\line(1,0){2}}\put(6.5,6.87){\line(1,0){2}}
\put(8.5,6.8){\line(1,1){0.5}}\put(8.5,6.8){\line(1,-1){0.5}}

\put(4,13.6){$\Longleftarrow$}\put(10,13.6){$\Longrightarrow$}
\put(4,10.1){$\Longleftarrow$}\put(10,10.1){$\Longrightarrow$}
\put(4,6.6){$\Longleftarrow$}\put(10,6.6){$\Longrightarrow$}

\put(0.4,13.9){$\bullet$}\put(2.4,13.4){$\bullet$}
\put(6.4,13.7){$\bullet$}\put(8.4,13.7){$\bullet$}
\put(12.4,13.4){$\bullet$}\put(14.4,13.9){$\bullet$}
\put(0.4,9.9){$\bullet$}\put(1.4,10.35){$\bullet$}
\put(6.4,10.2){$\bullet$}\put(7.4,10.25){$\bullet$}
\put(2.7,10.7){$\bullet$}\put(8.7,10.5){$\bullet$}
\put(14.7,10.7){$\bullet$}\put(12.4,10.35){$\bullet$}
\put(13.4,9.9){$\bullet$}\put(0.2,6.15){$\bullet$}
\put(1.4,6.9){$\bullet$}\put(6.2,6.5){$\bullet$}
\put(2.7,6.1){$\bullet$}\put(2.7,7.2){$\bullet$}
\put(7.4,6.75){$\bullet$}\put(8.7,6.4){$\bullet$}
\put(8.7,7){$\bullet$}\put(12.2,6.15){$\bullet$}
\put(14.6,6.2){$\bullet$}\put(14.6,7.1){$\bullet$}
\put(13.4,6.4){$\bullet$}

\put(1.4,14.2){$H_+$}\put(7.4,14.2){$H_0$}
\put(13.4,14.2){$H_-$}
\put(1.4,10.7){$H_+$}\put(7.4,10.7){$H_0$}
\put(13.4,10.7){$H_-$}
\put(1.4,7.2){$H_+$}\put(7.4,7.2){$H_0$}
\put(13.4,7.2){$H_-$}

\put(7.3,12){(a)}\put(7.3,8.5){(b)}
\put(7.3,5){(c)}\put(7.3,0){(d)}

\put(1.5,1){\line(-1,1){1.5}}\put(0,2.5){\line(1,1){1.5}}
\put(1.5,1){\line(1,1){1.5}}\put(1.5,1){\line(0,1){3}}
\put(1.5,4){\line(1,-1){1.5}}\put(0,2.5){\line(3,1){1.5}}
\put(3,2.5){\line(-3,1){1.5}}

\put(7.5,1){\line(-1,1){1.5}}\put(6,2.5){\line(1,1){1.5}}
\put(7.5,1){\line(1,1){1.5}}\put(7.5,1){\line(0,1){3}}
\put(7.5,4){\line(1,-1){1.5}}

\put(13.5,1){\line(-1,1){1.5}}\put(12,2.5){\line(1,1){1.5}}
\put(13.5,1){\line(1,1){1.5}}\put(13.5,1){\line(0,1){3}}
\put(13.5,4){\line(1,-1){1.5}}\put(12,2.5){\line(3,-1){1.5}}
\put(15,2.5){\line(-3,-1){1.5}}

\put(4,2.4){$\Longleftarrow$}\put(10,2.4){$\Longrightarrow$}

\put(0.9,3.1){$\mu_1$}\put(0.9,2){$\mu_2$}
\put(1.7,3.1){$\mu_3$}\put(1.7,2){$\mu_4$}
\put(12.9,2.8){$\mu_2$}\put(12.9,1.8){$\mu_1$}
\put(13.7,2.8){$\mu_4$}\put(13.7,1.8){$\mu_3$}

\end{picture}
\caption{Degeneration of type (2v)}\label{fig6}
\end{figure}

\subsection{Degenerations of types (2vi) and (2vii), I}\label{sec3.6}
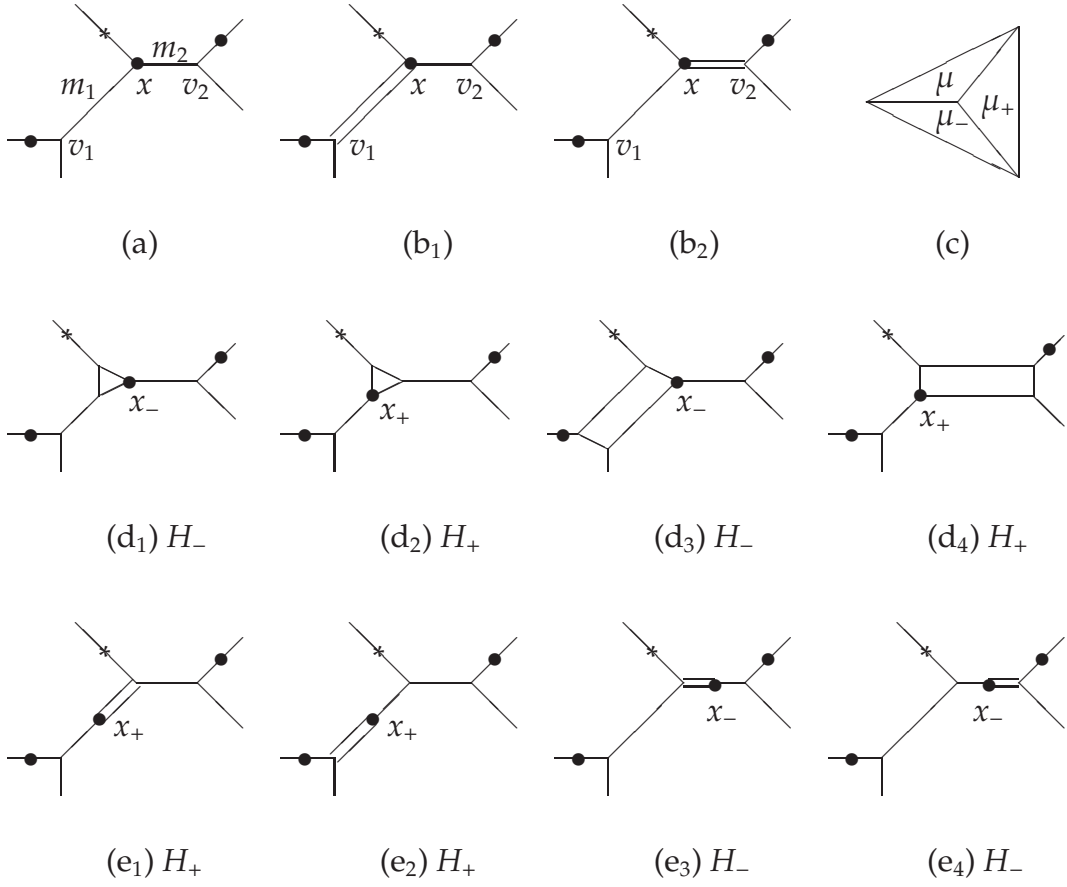
\begin{figure}
\setlength{\unitlength}{1cm}
\begin{picture}(14,11.5)(0,0)
\thinlines

\put(0.7,1){\line(0,1){0.5}}\put(0.7,1.5){\line(-1,0){0.7}}
\put(0.7,1.5){\line(1,1){0.5}}\put(1.25,1.95){\line(1,1){0.5}}
\put(1.15,2.05){\line(1,1){0.5}}\put(1.7,2.5){\line(-1,1){0.8}}
\put(1.7,2.5){\line(1,0){0.8}}\put(2.5,2.5){\line(1,1){0.6}}
\put(2.5,2.5){\line(1,-1){0.6}}

\put(4.3,1){\line(0,1){0.5}}\put(4.3,1.5){\line(-1,0){0.7}}
\put(4.8,2){\line(1,1){0.5}}\put(4.35,1.45){\line(1,1){0.5}}
\put(4.25,1.55){\line(1,1){0.5}}\put(5.3,2.5){\line(-1,1){0.8}}
\put(5.3,2.5){\line(1,0){0.8}}\put(6.1,2.5){\line(1,1){0.6}}
\put(6.1,2.5){\line(1,-1){0.6}}

\put(7.9,1){\line(0,1){0.5}}\put(7.9,1.5){\line(-1,0){0.7}}
\put(7.9,1.5){\line(1,1){1}}\put(8.9,2.5){\line(-1,1){0.8}}
\put(9.3,2.5){\line(1,0){0.4}}\put(8.9,2.55){\line(1,0){0.4}}
\put(8.9,2.45){\line(1,0){0.4}}\put(9.7,2.5){\line(1,1){0.6}}
\put(9.7,2.5){\line(1,-1){0.6}}

\put(11.5,1){\line(0,1){0.5}}\put(11.5,1.5){\line(-1,0){0.7}}
\put(11.5,1.5){\line(1,1){1}}\put(12.5,2.5){\line(-1,1){0.8}}
\put(12.5,2.5){\line(1,0){0.4}}\put(12.9,2.55){\line(1,0){0.4}}
\put(12.9,2.45){\line(1,0){0.4}}\put(13.3,2.5){\line(1,1){0.6}}
\put(13.3,2.5){\line(1,-1){0.6}}

\put(0.2,1.37){$\bullet$}\put(3.8,1.37){$\bullet$}
\put(7.4,1.37){$\bullet$}\put(11,1.37){$\bullet$}
\put(1.2,2.8){$*$}\put(4.8,2.8){$*$}
\put(8.4,2.8){$*$}\put(12,2.8){$*$}
\put(2.7,2.7){$\bullet$}\put(6.3,2.7){$\bullet$}
\put(9.9,2.7){$\bullet$}\put(13.5,2.7){$\bullet$}
\put(1.1,1.9){$\bullet$}\put(4.7,1.9){$\bullet$}
\put(9.2,2.35){$\bullet$}\put(12.8,2.35){$\bullet$}

\put(1.3,0){(e$_1$)\ $H_+$}\put(4.9,0){(e$_2$)\ $H_+$}
\put(8.5,0){(e$_3$)\ $H_-$}\put(12.1,0){(e$_4$)\ $H_-$}

\put(0.7,5.3){\line(0,1){0.5}}\put(0.7,5.8){\line(-1,0){0.7}}
\put(0.7,5.8){\line(1,1){0.5}}\put(1.2,6.3){\line(0,1){0.4}}
\put(1.2,6.3){\line(2,1){0.4}}\put(1.6,6.5){\line(1,0){0.9}}
\put(1.2,6.7){\line(2,-1){0.4}}\put(1.2,6.7){\line(-1,1){0.6}}
\put(2.5,6.5){\line(1,-1){0.5}}\put(2.5,6.5){\line(1,1){0.5}}

\put(4.3,5.3){\line(0,1){0.5}}\put(4.3,5.8){\line(-1,0){0.7}}
\put(4.3,5.8){\line(1,1){0.5}}\put(4.8,6.3){\line(0,1){0.4}}
\put(4.8,6.3){\line(2,1){0.4}}\put(5.2,6.5){\line(1,0){0.9}}
\put(4.8,6.7){\line(2,-1){0.4}}\put(4.8,6.7){\line(-1,1){0.6}}
\put(6.1,6.5){\line(1,-1){0.5}}\put(6.1,6.5){\line(1,1){0.5}}

\put(7.9,5.3){\line(0,1){0.3}}\put(7.5,5.8){\line(-1,0){0.4}}
\put(7.9,5.6){\line(-2,1){0.4}}\put(8.4,6.7){\line(-1,-1){0.9}}
\put(8.8,6.5){\line(-1,-1){0.9}}\put(8.8,6.5){\line(1,0){0.9}}
\put(8.4,6.7){\line(2,-1){0.4}}\put(8.4,6.7){\line(-1,1){0.6}}
\put(9.7,6.5){\line(1,-1){0.5}}\put(9.7,6.5){\line(1,1){0.5}}

\put(11.5,5.3){\line(0,1){0.5}}\put(11.5,5.8){\line(-1,0){0.7}}
\put(11.5,5.8){\line(1,1){0.5}}\put(12,6.3){\line(0,1){0.4}}
\put(12,6.3){\line(1,0){1.5}}\put(13.5,6.3){\line(0,1){0.4}}
\put(12,6.7){\line(1,0){1.5}}\put(12,6.7){\line(-1,1){0.6}}
\put(13.5,6.3){\line(1,-1){0.4}}\put(13.5,6.7){\line(1,1){0.4}}

\put(0.2,5.67){$\bullet$}\put(3.8,5.67){$\bullet$}
\put(7.2,5.67){$\bullet$}\put(11,5.67){$\bullet$}
\put(0.7,7){$*$}\put(4.3,7){$*$}
\put(7.9,7){$*$}\put(11.5,7){$*$}
\put(2.7,6.7){$\bullet$}\put(6.3,6.7){$\bullet$}
\put(9.9,6.7){$\bullet$}\put(13.6,6.8){$\bullet$}
\put(1.5,6.37){$\bullet$}\put(4.7,6.2){$\bullet$}
\put(8.7,6.37){$\bullet$}\put(11.9,6.2){$\bullet$}

\put(1.3,4.3){(d$_1$)\ $H_-$}\put(4.9,4.3){(d$_2$)\ $H_+$}
\put(8.5,4.3){(d$_3$)\ $H_-$}\put(12.1,4.3){(d$_4$)\ $H_+$}

\put(0.7,9.2){\line(0,1){0.5}}\put(0.7,9.7){\line(-1,0){0.7}}
\put(0.7,9.7){\line(1,1){1}}
\put(1.7,10.7){\line(-1,1){0.8}}
\put(1.7,10.7){\line(1,0){0.8}}\put(2.5,10.7){\line(1,1){0.6}}
\put(2.5,10.7){\line(1,-1){0.6}}

\put(4.3,9.2){\line(0,1){0.5}}\put(4.3,9.7){\line(-1,0){0.7}}
\put(4.35,9.65){\line(1,1){1}}
\put(4.25,9.75){\line(1,1){1}}\put(5.3,10.7){\line(-1,1){0.8}}
\put(5.3,10.7){\line(1,0){0.8}}\put(6.1,10.7){\line(1,1){0.6}}
\put(6.1,10.7){\line(1,-1){0.6}}

\put(7.9,9.2){\line(0,1){0.5}}\put(7.9,9.7){\line(-1,0){0.7}}
\put(7.9,9.7){\line(1,1){1}}\put(8.9,10.7){\line(-1,1){0.8}}
\put(8.9,10.75){\line(1,0){0.8}}
\put(8.9,10.65){\line(1,0){0.8}}\put(9.7,10.7){\line(1,1){0.6}}
\put(9.7,10.7){\line(1,-1){0.6}}

\put(0.2,9.57){$\bullet$}\put(3.8,9.57){$\bullet$}
\put(7.4,9.57){$\bullet$}
\put(1.2,11){$*$}\put(4.8,11){$*$}
\put(8.4,11){$*$}
\put(2.7,10.9){$\bullet$}\put(6.3,10.9){$\bullet$}
\put(9.9,10.9){$\bullet$}
\put(1.6,10.6){$\bullet$}\put(5.2,10.6){$\bullet$}
\put(8.8,10.6){$\bullet$}

\put(11.3,10.2){\line(2,1){2}}\put(11.3,10.2){\line(2,-1){2}}
\put(13.3,9.2){\line(0,1){2}}\put(12.5,10.2){\line(4,5){0.8}}
\put(12.5,10.2){\line(4,-5){0.8}}\put(12.5,10.2){\line(-1,0){1.2}}

\put(1.5,8.2){(a)}\put(5.2,8.2){(b$_1$)}\put(8.7,8.2){(b$_2$)}\put(12.2,8.2){(c)}

\put(1.7,10.3){$x$}\put(5.3,10.3){$x$}\put(8.9,10.3){$x$}
\put(0.8,9.5){$v_1$}\put(4.5,9.5){$v_1$}\put(8,9.5){$v_1$}
\put(2.3,10.3){$v_2$}\put(5.9,10.3){$v_2$}\put(9.5,10.3){$v_2$}

\put(1.6,6.1){$x_-$}\put(4.9,6){$x_+$}\put(8.8,6.1){$x_-$}\put(12,5.9){$x_+$}
\put(1.4,1.8){$x_+$}\put(5,1.8){$x_+$}\put(9.2,2){$x_-$}\put(12.7,2){$x_-$}

\put(12.2,10.4){$\mu$}\put(12.2,9.9){$\mu_-$}\put(12.8,10.1){$\mu_+$}
\put(0.7,10.3){$m_1$}\put(1.9,10.8){$m_2$}

\end{picture}
\caption{Degenerations of types (2vi) and (2vii), I}\label{fig7}
\end{figure}

We consider degenerations of types (2vi) and (2vii) together. Let $C^*$ be either as in
as in Lemma \ref{l6}(2vi) or as in Lemma \ref{l6}(2vii) (with the same embedded plane
tropical curve
$T=h_*(C^*)$ in both the situations, see Figure \ref{fig7}(a,b)).
The closure of the bounded component of $\overline\Gamma\setminus\bp$ contains at least
one marked vertex.
Let $p\in\bp_v$ be one of them (see Figures \ref{fig7}(a,b), where $x=h(p)$ and the bounded
component of
$\overline\Gamma\setminus\bp$ is marked by asterisk).
In this section, we consider the case when the edges of the unbounded components of
$\overline\Gamma\setminus\bp$ incident to $p$ have finite length.
Assume also that the weights of the corresponding edges of $T$ satisfy $m_1,m_2>1$.
The case of $m_1=1$ or $m_2=1$ can
be treated in the same (in fact, rather simpler) way.

Assuming that the germ of the path $(\bx(t)$, $t\in(\R,t^*)$, is such that the point
$x$ moves along a generic line, while the rest of $\bx=h(\bp)$ is fixed, we encounter the
deformations of
the above curves
shown in Figures \ref{fig7}(d,e), where we also label $2n$-cells of
$\overline{\mathcal M}^{\;e}_{1,(n_v,n_e)}(\R^2,\Delta)$ projecting onto
$H_+$ or $H_-$.
Notice that the curves shown in Figures \ref{fig7}(d$_1$-d$_4$) belong to
${\mathcal M}'_{1,(n_v,n_e)}(\R^2,\Delta)$ and the curves shown in Figures
\ref{fig7}(e$_1$-e$_4$) belong to ${\mathcal M}''_{1,(n_v,n_e)}(\R^2,\Delta)$.
Introduce the following parameters of the curve $C^*$ shown in Figure \ref{fig7}(a):
$$\nu_i=\frac{\mu(\Gamma,h,V_i)}{m_i},\quad \nu'_i=
\frac{\mu(\Gamma,h,p)}{m_i},\quad i=1,2\ .$$
By definition, we can express the contributions of the curves shown in Figures
\ref{fig7}(d$_1$-d$_4$), respectively, to $\GS_y(\Delta,1,(n_v,n_e),\bx^{(t)})$, $t\in(\R,t^*)$, as
\begin{equation}c_1=\varphi(y)(z^{\nu_1m_1}-z^{-\nu_1m_1})(z^{\nu_2m_2}-z^{-\nu_2m_2})
\sum_{v\in\Int{\mathcal D}(p)\cap\Z^2}\sigma_+(v)\ ,\label{e10}\end{equation}
\begin{equation}c_2=\varphi(y)(z^{\nu_1m_1}-z^{-\nu_1m_1})(z^{\nu_2m_2}-z^{-\nu_2m_2})
\sum_{v\in\Int{\mathcal D}(p)\cap\Z^2}\sigma_-(v)\ ,\label{e11}\end{equation}
$$c_3=\varphi(y)(z^{\nu_2m_2}-z^{-\nu_2m_2})\sum_{i=1}^{m_1-1}\big((z^{\nu_1i}-z^{-\nu_1i})
(z^{\nu_1(m_1-i)}-z^{-\nu_1(m_1-i)})$$
$$\qquad\qquad\times(z^{\nu'_1i}+z^{-\nu'_1i})(z^{\nu'_1(m_1-i)}-z^{-\nu'_1(m_1-i)})\big)\ ,$$
$$c_4=\varphi(y)(z^{\nu_1m_1}-z^{-\nu_1m_1})\sum_{i=1}^{m_2-1}\big((z^{\nu_2i}-z^{-\nu_2i})
(z^{\nu_2(m_2-i)}-z^{-\nu_2(m_2-i)})$$
$$\qquad\qquad\times(z^{\nu'_2i}+z^{-\nu'_2i})(z^{\nu'_2(m_2-i)}-z^{-\nu'_2(m_2-i)})\big)\ ,$$
where $y=z^2$, $\varphi(y)$ is the contribution of the vertices outside the
fragments shown in Figure \ref{fig7}, and
\begin{equation}
\begin{cases}&\sigma_+=(z^{\mu(v)}-z^{-\mu(v)})(z^{\mu_+(v)}+
z^{-\mu_+(v)})(z^{\mu_-(v)}-z^{-\mu_-(v)}),\\
&\sigma_-=(z^{\mu(v)}-z^{-\mu(v)})(z^{\mu_+(v)}-z^{-\mu_+(v)})
(z^{\mu_-(v)}+z^{-\mu_-(v)}),\end{cases}
\label{e12}\end{equation}
$\mu(v),\mu_+(v),\mu_-(v)$ being the lattice areas of the triangles in the
subdivision of ${\mathcal D}(p)$
with vertex $v$ (see Figure \ref{fig7}(c)). Correspondingly, the
contributions of the curves
shown in Figures \ref{fig7}(e$_1$-e$_2$) are
$$d_1=\varphi(y)(z^{\nu_1m_1}-z^{-\mu_1m_1})\Psi_z^{(2)}(\nu_2,\nu'_2,m_2),
\quad d_2=\varphi(y)(z^{\nu_2m_2}-z^{-\mu_2m_2})\Psi_z^{(2)}(\nu_1,\nu'_1,m_1)\ .$$
The desired invariance of $\GS_y(\Delta,1,(n_v,n_e),\bx^{(t)})$, $t\in(\R,t^*)$, reads
\begin{equation}c_1+c_3+d_1=c_2+c_4+d_2\ .\label{e5}\end{equation}

\begin{lemma}\label{l9}
\begin{eqnarray}c_2-c_1&=&2\varphi(y)(z^{\nu_1m_1}-z^{-\mu_1m_1})(z^{\nu_2m_2}-
z^{-\nu_2m_2})\nonumber\\
& &\times\Bigg[m_2\frac{z^{\nu'_2m_2-m_2}-z^{m_2-\nu'_2m_2}}{z^{m_2}-z^{-m_2}}-
m_1\frac{z^{\nu'_1m_1-m_1}-z^{m_1-\nu'_1m_1}}{z^{m_1}-z^{-m_1}}\nonumber\\
& &\qquad+\frac{z^{\nu'_2m_2-\nu'_2}-z^{\nu'_2-\nu'_2m_2}}{z^{\nu'_2}-z^{-\nu'_2}}
-\frac{z^{\nu'_1m_1-\nu'_1}-z^{\nu'_1-\nu'_1m_1}}{z^{\nu'_1}-z^{-\nu'_1}}\Bigg]\ ,
\label{e6}\end{eqnarray}
\begin{eqnarray}c_4-c_3&=&\varphi(y)(z^{\nu'_1m_1}-z^{-\nu'_1m_1})\Bigg[(m_2-1)
(z^{\nu_1m_1}-z^{-\nu_1m_1})
(z^{\nu_2m_2}+z^{-\nu_2m_2})\nonumber\\
& &\qquad
-(m_1-1)(z^{\nu_1m_1}+z^{-\nu_1m_1})(z^{\nu_2m_2}-z^{-\nu_2m_2})\nonumber\\
& &\qquad+2(z^{\nu_2m_2}-z^{-\nu_2m_2})\frac{z^{\nu_1m_1-\nu_1}-
z^{\nu_1-\nu_1m_1}}{z^{\nu_1}-z^{-\nu_1}}
\nonumber\\
& &\qquad-2(z^{\nu_1m_1}-z^{-\nu_1m_1})\frac{z^{\nu_2m_2-\nu_2}
-z^{\nu_2-\nu_2m_2}}{z^{\nu_2}-z^{-\nu_2}}\Bigg]\ .
\label{e7}\end{eqnarray}
\end{lemma}

Before proving Lemma \ref{l9}, we note that (\ref{e5}) follows from (\ref{e6})
and (\ref{e7}), since along (\ref{e8}) and (\ref{e9}) one gets
$$c_3=\varphi(y)(z^{\nu_2m_2}-z^{-\nu_2m_2})$$
$$\times\Bigg[2(z^{\nu_1m_1}-z^{-\nu_1m_1})\frac{z^{\nu'_1m_1-1}-z^{1-\nu'_1m_1}}{z-z^{-1}}
-2m_1(z^{\nu_1m_1}-z^{-\nu_1m_1})\frac{z^{\nu'_1m_1-m_1}-
z^{m_1-\nu'_1m_1}}{z^{m_1}-z^{-m_1}}$$
$$+(m_1-1)(z^{\nu'_1m_1}-z^{\nu'_1m_1})(z^{\nu_1m_1}+z^{-\nu_1m_1})$$
$$-2(z^{\nu_1m_1}-z^{-\nu_1m_1})\frac{z^{\nu'_1m_1-\nu'_1}
-z^{\nu'_1-\nu'_1m_1}}{z^{\nu'_1}-z^{-\nu'_1}}
-2(z^{\nu'_1m_1}-z^{\nu'_1m_1})\frac{z^{\nu_1m_1-\nu_1}-z^{\nu_1-
\nu_1m_1}}{z^{\nu_1}-z^{-\nu_1}}\Bigg]$$
and
$$c_4=\varphi(y)(z^{\nu_1m_1}-z^{-\nu_1m_1})$$
$$\times\Bigg[2(z^{\nu_2m_2}-z^{-\nu_2m_2})\frac{z^{\nu'_2m_2-1}-z^{1-\nu'_2m_2}}{z-z^{-1}}
-2m_2(z^{\nu_2m_2}-z^{-\nu_2m_2})\frac{z^{\nu'_2m_2-m_2}-
z^{m_2-\nu'_2m_2}}{z^{m_2}-z^{-m_2}}$$
$$+(m_2-1)(z^{\nu'_2m_2}-z^{\nu'_2m_2})(z^{\nu_2m_2}+z^{-\nu_2m_2})$$
$$-2(z^{\nu_2m_2}-z^{-\nu_2m_2})\frac{z^{\nu'_2m_2-\nu'_2}
-z^{\nu'_2-\nu'_2m_2}}{z^{\nu'_2}-z^{-\nu'_2}}
-2(z^{\nu'_2m_2}-z^{\nu'_2m_2})\frac{z^{\nu_2m_2-\nu_2}-z^{\nu_2-
\nu_2m_2}}{z^{\nu_2}-z^{-\nu_2}}\Bigg]\ ,$$
and the first summand in the brackets cancels out in $c_4-c_3$ in view of $\nu'_1m_1=
\nu'_2m_2$.

{\bf Proof of Lemma \ref{l9}.} We start with formula (\ref{e7}) which is simpler. We have
$$c_3=\varphi(y)(z^{\nu_2m_2}-z^{-\nu_2m_2})$$ $$\times\sum_{i=1}^{m_1-1}
\left[(z^{\nu'_1i}+z^{-\nu'_1i})
(z^{\nu'_1(m_1-i)}-z^{-\nu'_1(m_1-i)})(z^{\nu_1i}-z^{-\nu_1i})(z^{\nu_1(m_1-i)}
-z^{-\nu_1(m_1-i)})\right]$$
$$=\varphi(y)(z^{\nu_2m_2}-z^{-\nu_2m_2})$$
$$\times\sum_{i=1}^{m_1-1}\left[z^{(\nu_1+\nu'_1)m_1}+z^{(\nu'_1-\nu_1)m_1}
-z^{(\nu_1-\nu'_1)m_1}-z^{-(\nu'_1+\nu_1)m_1}
-2(z^{\nu'_1m_1}-z^{-\nu'_1m_1})z^{(2i-m_1)\nu_1}\right]$$
$$=\varphi(y)(z^{\nu_2m_2}-z^{-\nu_2m_2})(z^{\nu'_1m_1}-z^{-\nu'_1m_1})$$
\begin{equation}\times\left[(m_1-1)(z^{\nu_1m_1}+z^{-\nu_1m_1})-2\frac{z^{\nu_1m_1-m_1}-
z^{m_1-\nu_1m_1}}{z^{m_1}-z^{-m_1}}\right]\ .\label{e14a}\end{equation}
Similarly (notice that $\nu'_1m_1=\nu'_2m_2$),
$$c_4=\varphi(y)(z^{\nu_1m_1}-z^{-\nu_1m_1})(z^{\nu'_1m_1}-z^{-\nu'_1m_1})$$
$$\times\left[(m_2-1)(z^{\nu_2m_2}+z^{-\nu_2m_2})-2\frac{z^{\nu_2m_2-m_2}-
z^{m_2-\nu_2m_2}}{z^{m_2}-z^{-m_2}}\right]\ ,$$
and formula (\ref{e7}) follows.

In view of formulas (\ref{e10}), (\ref{e11}), and (\ref{e12}), and using the equality
$$\mu(v)+\mu_+(v)+\mu_-(v)=\mu:=\mu(\Gamma,h,p),\quad v\in\Int\Delta(p)\cap\Z^2\ ,$$
we obtain
$$c_2-c_1=\varphi(y)(z^{\nu_1m_1}-z^{-\nu_1m_1})(z^{\nu_2m_2}-
z^{-\nu_2m_2})$$
$$\times\sum_{v\in\Int\Delta(p)\cap\Z^2}\Big[(z^{\mu(v)}-z^{-\mu(v)})(z^{\mu_+(v)}-
z^{-\mu_+(v)})(z^{\mu_-(v)}+z^{-\mu_-(v)})$$
$$-(z^{\mu(v)}-z^{-\mu(v)})(z^{\mu_+(v)}+
z^{-\mu_+(v)})(z^{\mu_-(v)}-z^{-\mu_-(v)})\Big]$$
$$=2\varphi(y)(z^{\nu_1m_1}-z^{-\mu_1m_1})(z^{\nu_2m_2}-
z^{-\nu_2m_2})$$
$$\Bigg[\sum_{v\in\Int\Delta(p)\cap\Z^2}(z^{\mu-2\mu_-(v)}+z^{2\mu_-(v)-\mu})
-\sum_{v\in\Int\Delta(p)\cap\Z^2}(z^{\mu-2\mu_+(v)}+z^{2\mu_+(v)-\mu})\Bigg]\ .$$
To complete the proof of (\ref{e6}), and thereby of Lemma \ref{l9}, we use the relations
$$2\sum_{v\in\Int\Delta(p)\cap\Z^2}(z^{\mu-2\mu_-(v)}+z^{2\mu_-(v)-\mu})=
2m_2\frac{z^{\nu'_2m_2-m_2}-z^{m_2-\nu'_2m_2}}{z^{m_2}-z^{-m_2}}$$
\begin{equation}-
\frac{z^{\nu'_1m_1-\nu'_1}-z^{\nu'_1-\nu'_1m_1}}{z^{\nu'_1}-z^{-\nu'_1}}
-\frac{z^{\nu'_3m_3-\nu'_3}-z^{\nu'_3-\nu'_3m_3}}{z^{\nu'_3}-z^{-\nu'_3}}
\label{ea3}\end{equation}
and
$$2\sum_{v\in\Int\Delta(p)\cap\Z^2}(z^{\mu-2\mu_+(v)}+z^{2\mu_+(v)-\mu})=
2m_1\frac{z^{\nu'_1m_1-m_1}-z^{m_1-\nu'_1m_1}}{z^{m_1}-z^{-m_1}}$$
\begin{equation}-
\frac{z^{\nu'_2m_2-\nu'_2}-z^{\nu'_2-\nu'_2m_2}}{z^{\nu'_2}-z^{-\nu'_2}}
-\frac{z^{\nu'_3m_3-\nu'_3}-z^{\nu'_3-\nu'_3m_3}}{z^{\nu'_3}-z^{-\nu'_3}}\ ,
\label{ea4}\end{equation}
where $m_3$ is the weight of the third edge incident to the vertex $p$ of $\Gamma$,
and $\nu'_3=\mu(\Gamma,h,p)/m_3$.
Finally, both relations, (\ref{ea3}) and (\ref{ea4}), follow from Lemma \ref{la1} below.
\proofend

\begin{lemma}\label{la1}
Let ${\mathcal T}\subset\R^2$ be the lattice triangle with vertices
$(0,0)$, $(m,0)$, $(k,l)$, where $l,m\ge1$. Denote by $(x_v,y_v)$ the
coordinates of a point $v\in\R^2$. Then
$$2\sum_{v\in\Int{\mathcal T}\cap\Z^2}(z^{m(l-2y_v)}+z^{m(2y_v-l)})
+\sum_{\renewcommand{\arraystretch}{0.6}
\begin{array}{c}
\scriptstyle{v\in\partial{\mathcal T}\cap\Z^2}\\
\scriptstyle{0<y_v<l}
\end{array}}(z^{m(l-2y_v)}+z^{m(2y_v-l)})$$
\begin{equation}=m\frac{z^{ml-m}-z^{m-ml}}{z^m-z^{-m}}\ .\label{ea1}\end{equation}
\end{lemma}

{\bf Proof.}
(1) Suppose that $m=1$.

First, we claim that the sequences $(l-2y_v)_{v\in{\mathcal T}\cap\Z^2,0<y_v<l}$ and
$(2y_y-l)_{v\in\Int{\mathcal T}\cap\Z^2}$ are disjoint, and each number
$s\equiv l\mod2$, $-l<s<l$, appears once in one of the sequences.
Since in each of the sequences, the elements are distinct, and the total number of
elements in both sequences equals $2l-2$ (Pick's formula), for the above
claim it is enough to verify that the
sequences are disjoint. Indeed, let $v_1=(x_1,y_1)\in{\mathcal T}\cap Z^2$, $0<y_1<l$, $v_2=(x_2,y_2)\in\Int{\mathcal T}\cap\Z^2$,
and $l-2y_1=2y_2-l$, or, equivalently, $y_2=l-y_1$. We have the following relations for slopes:
$$\frac{x_1}{y_1}\ge\frac{k}{l},\quad \frac{x_1-1}{y_1}\le\frac{k-1}{l},\quad\frac{x_2}{y_2}>\frac{k}{l},\quad
\frac{x_2-1}{y_2}<\frac{k-1}{l}\ ,$$
which yield
$$\frac{k}{l}y_1\le x_1\le\frac{k}{l}y_1+\frac{l-y_1}{l},\quad
\frac{k}{l}y_2<x_1<\frac{k}{l}y_2+\frac{l-y_2}{l}\ .$$ Plugging $y_2=l-y_1$ to
the second relation, we obtain
$$\frac{k}{l}y_1-\frac{y_1}{l}<k-x_2<\frac{k}{l}y_1\ ,$$ that is, two
integers $k-x_2<x_1$ in the unit
interval $\big(\frac{k}{l}y_1-\frac{y_1}{l},\frac{k}{l}y_1+\frac{l-y_1}{l}\big]$,
which is a contradiction.

Second, we notice that the sequences $(l-2y_v)_{v\in\partial{\mathcal T}\cap\Z^2}$ and $(2y_v-l)_{v\in\partial{\mathcal T}\cap\Z^2}$ coincide.

Both claims together yield that
$$2\sum_{v\in\Int{\mathcal T}\cap\Z^2}(z^{l-2y_v}+z^{2y_v-l})
+\sum_{\renewcommand{\arraystretch}{0.6}
\begin{array}{c}
\scriptstyle{v\in\partial{\mathcal T}\cap\Z^2}\\
\scriptstyle{0<y_v<l}
\end{array}}(z^{l-2y_v}+z^{2y_v-l})$$
\begin{equation}=2\sum_{i=1}^{l-1}z^{l-2i}=2\frac{z^{l-1}-z^{1-l}}{z-z^{-1}}\ .
\label{ea2}\end{equation}

(2) For an arbitrary $m\ge1$, we divide the triangle ${\mathcal T}$ into the triangles
\mbox{${\mathcal T}_s=\conv\big\{(s-1,0),(s,0),(k,l)\big\}$}, $s=1,...,m$,
sum up the formulas (\ref{ea2})
for ${\mathcal T}_s$, $s=1,...,m$, and substitute $z^m$ for $z$, finally obtaining (\ref{ea1}).
\proofend

\subsection{Degenerations of types (2vi) and (2vii), II}\label{sec3.6a}
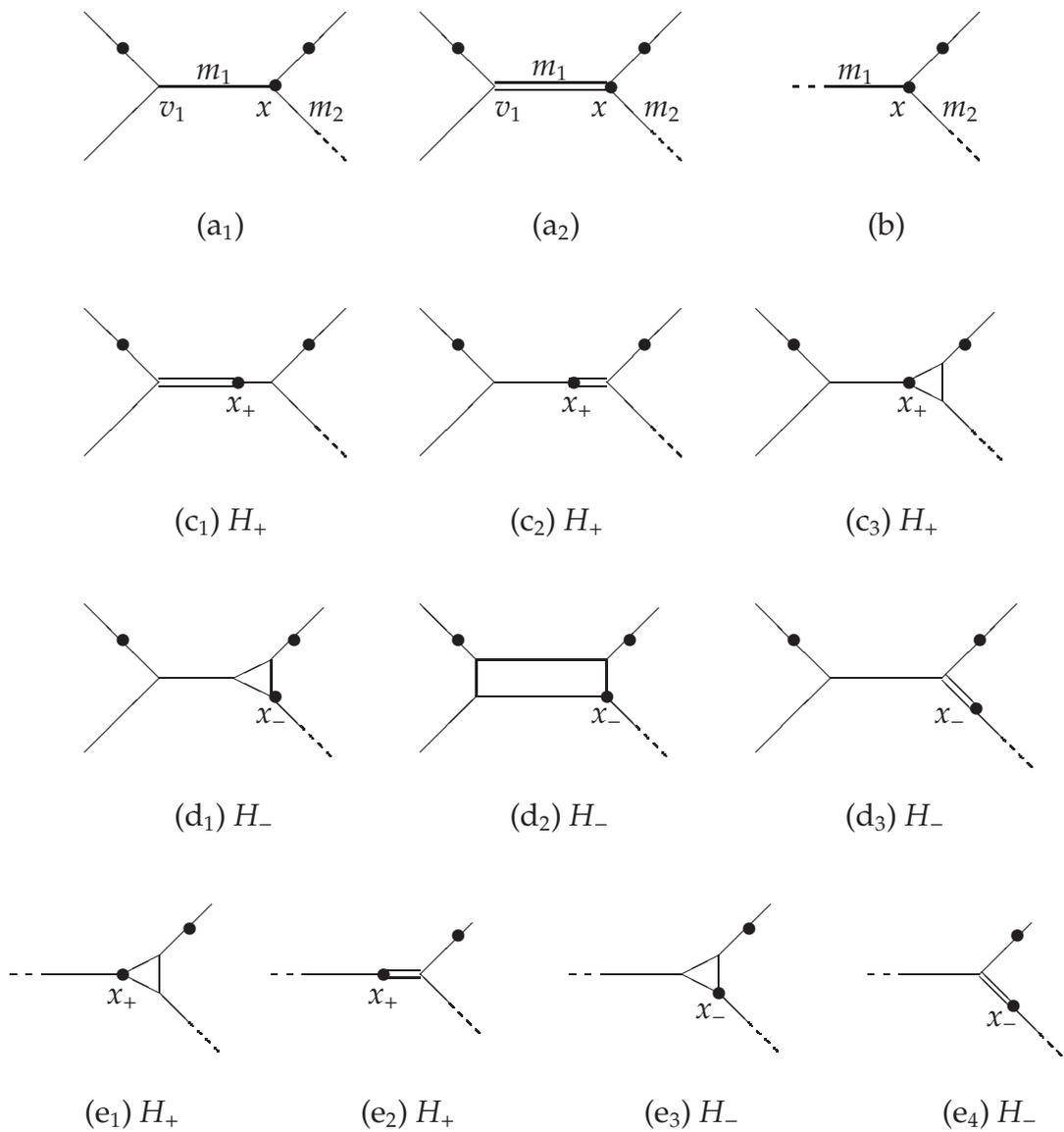
\begin{figure}
\setlength{\unitlength}{1cm}
\begin{picture}(14,15)(0,0)
\thinlines

\put(1,13){\line(1,1){1}}\put(2,14){\line(-1,1){1}}
\put(2,14){\line(1,0){1.5}}\put(3.5,14){\line(1,1){1}}
\put(3.5,14){\line(1,-1){0.6}}
\dashline{0.1}(4.1,13.4)(4.5,13)

\put(5.5,13){\line(1,1){1}}\put(6.5,14){\line(-1,1){1}}
\put(6.5,14.05){\line(1,0){1.5}}\put(6.5,13.95){\line(1,0){1.5}}\put(8,14){\line(1,1){1}}
\put(8,14){\line(1,-1){0.6}}
\dashline{0.1}(8.6,13.4)(9,13)

\put(11,14){\line(1,0){1}}\put(12,14){\line(1,1){1}}\put(12,14){\line(1,-1){0.6}}
\dashline{0.1}(12.6,13.4)(13,13)
\dashline{0.1}(10.5,14)(11,14)

\put(1,9){\line(1,1){1}}\put(2,10){\line(-1,1){1}}
\put(2,10.05){\line(1,0){1}}\put(2,9.95){\line(1,0){1}}
\put(3,10){\line(1,0){0.5}}\put(3.5,10){\line(1,1){1}}
\put(3.5,10){\line(1,-1){0.6}}
\dashline{0.1}(4.1,9.4)(4.5,9)

\put(5.5,9){\line(1,1){1}}\put(6.5,10){\line(-1,1){1}}
\put(7.5,10.05){\line(1,0){0.5}}\put(7.5,9.95){\line(1,0){0.5}}
\put(6.5,10){\line(1,0){1}}\put(8,10){\line(1,1){1}}
\put(8,10){\line(1,-1){0.6}}
\dashline{0.1}(8.6,9.4)(9,9)

\put(10,9){\line(1,1){1}}\put(11,10){\line(-1,1){1}}
\put(12,10){\line(2,1){0.5}}\put(12,10){\line(2,-1){0.5}}
\put(12.5,9.75){\line(0,1){0.5}}\put(12.5,10.25){\line(1,1){0.7}}
\put(12.5,9.75){\line(1,-1){0.4}}\put(11,10){\line(1,0){1}}
\dashline{0.1}(12.9,9.35)(13.3,8.95)

\put(1,5){\line(1,1){1}}\put(2,6){\line(-1,1){1}}
\put(3,6){\line(2,1){0.5}}\put(3,6){\line(2,-1){0.5}}
\put(3.5,5.75){\line(0,1){0.5}}\put(3.5,6.25){\line(1,1){0.7}}
\put(3.5,5.75){\line(1,-1){0.4}}\put(2,6){\line(1,0){1}}
\dashline{0.1}(3.9,5.35)(4.3,4.95)

\put(5.5,5){\line(1,1){0.75}}\put(6.25,6.25){\line(-1,1){0.75}}
\put(8,5.75){\line(0,1){0.5}}\put(8,6.25){\line(1,1){0.7}}
\put(8,5.75){\line(1,-1){0.4}}\put(6.25,6.25){\line(1,0){1.75}}
\put(6.25,5.75){\line(0,1){0.5}}\put(6.25,5.75){\line(1,0){1.75}}
\dashline{0.1}(8.4,5.35)(8.8,4.95)

\put(10,5){\line(1,1){1}}\put(11,6){\line(-1,1){1}}
\put(11,6){\line(1,0){1.5}}\put(12.5,6){\line(1,1){1}}
\put(12.55,6.05){\line(1,-1){0.4}}\put(12.5,5.97){\line(1,-1){0.4}}
\put(12.9,5.6){\line(1,-1){0.4}}
\dashline{0.1}(13.3,5.2)(13.7,4.8)

\dashline{0.1}(0,2)(0.5,2)\dashline{0.1}(2.4,1.35)(2.8,0.95)
\put(0.5,2){\line(1,0){1}}\put(1.5,2){\line(2,1){0.5}}
\put(1.5,2){\line(2,-1){0.5}}\put(2,1.75){\line(0,1){0.5}}
\put(2,1.75){\line(1,-1){0.4}}\put(2,2.25){\line(1,1){0.7}}

\dashline{0.1}(7.5,2)(8,2)\dashline{0.1}(9.9,1.35)(10.3,0.95)
\put(8,2){\line(1,0){1}}\put(9,2){\line(2,1){0.5}}
\put(9,2){\line(2,-1){0.5}}\put(9.5,1.75){\line(0,1){0.5}}
\put(9.5,1.75){\line(1,-1){0.4}}\put(9.5,2.25){\line(1,1){0.7}}

\dashline{0.1}(3.5,2)(4,2)\dashline{0.1}(5.9,1.6)(6.3,1.2)
\put(4,2){\line(1,0){1}}\put(5,2.05){\line(1,0){0.5}}
\put(5,1.95){\line(1,0){0.5}}
\put(5.5,2){\line(1,-1){0.4}}\put(5.5,2){\line(1,1){0.7}}

\dashline{0.1}(11.5,2)(12,2)\dashline{0.1}(13.8,1.2)(14.1,0.9)
\put(12,2){\line(1,0){1}}\put(13,1.97){\line(1,-1){0.4}}\put(13.4,1.6){\line(1,-1){0.4}}
\put(13.03,2.03){\line(1,-1){0.4}}\put(13,2){\line(1,1){0.7}}

\put(1.4,14.4){$\bullet$}\put(5.9,14.4){$\bullet$}\put(12.4,14.4){$\bullet$}
\put(3.9,14.4){$\bullet$}\put(8.4,14.4){$\bullet$}\put(3.45,13.9){$\bullet$}
\put(7.95,13.88){$\bullet$}\put(11.95,13.88){$\bullet$}

\put(1.4,10.4){$\bullet$}\put(5.9,10.4){$\bullet$}\put(10.4,10.4){$\bullet$}
\put(12.7,10.4){$\bullet$}
\put(3.9,10.4){$\bullet$}\put(8.4,10.4){$\bullet$}\put(2.95,9.88){$\bullet$}
\put(7.45,9.88){$\bullet$}\put(11.95,9.88){$\bullet$}

\put(1.4,6.4){$\bullet$}\put(5.9,6.4){$\bullet$}\put(10.4,6.4){$\bullet$}\put(12.9,6.4){$\bullet$}
\put(3.7,6.4){$\bullet$}\put(8.2,6.4){$\bullet$}\put(3.45,5.63){$\bullet$}
\put(7.9,5.63){$\bullet$}\put(12.85,5.48){$\bullet$}

\put(1.4,1.88){$\bullet$}\put(4.9,1.88){$\bullet$}\put(9.4,1.63){$\bullet$}\put(13.35,1.45){$\bullet$}
\put(2.3,2.5){$\bullet$}\put(5.9,2.4){$\bullet$}\put(9.8,2.5){$\bullet$}\put(13.4,2.4){$\bullet$}

\put(2,13.6){$v_1$}\put(6.5,13.6){$v_1$}\put(3.3,13.6){$x$}\put(7.8,13.6){$x$}\put(11.8,13.6){$x$}
\put(2.9,9.6){$x_+$}\put(7.4,9.6){$x_+$}\put(11.9,9.6){$x_+$}
\put(3.3,5.4){$x_-$}\put(7.8,5.4){$x_-$}\put(12.4,5.4){$x_-$}
\put(1.3,1.6){$x_+$}\put(4.8,1.6){$x_+$}\put(9.2,1.4){$x_-$}\put(13.1,1.3){$x_-$}

\put(2.5,14.1){$m_1$}\put(7,14.15){$m_1$}\put(11.1,14.1){$m_1$}
\put(4,13.6){$m_2$}\put(8.5,13.6){$m_2$}\put(12.5,13.6){$m_2$}

\put(2.5,12){(a$_1$)}\put(7,12){(a$_2$)}\put(11.5,12){(b)}
\put(2.2,8){(c$_1$)\ $H_+$}\put(6.7,8){(c$_2$)\ $H_+$}\put(11.2,8){(c$_3$)\ $H_+$}
\put(2.2,4){(d$_1$)\ $H_-$}\put(6.7,4){(d$_2$)\ $H_-$}\put(11.2,4){(d$_3$)\ $H_-$}
\put(1,0){(e$_1$)\ $H_+$}\put(4.7,0){(e$_2$)\ $H_+$}\put(8.5,0){(e$_3$)\ $H_-$}
\put(12.5,0){(e$_4$)\ $H_-$}

\end{picture}
\caption{Degeneration of type (2vi) and (2vii), II}\label{fig11}
\end{figure}

Now, we go back to the hypotheses of Section \ref{sec3.6} with the following modification:
we suppose that one or two edges of the unbounded components of
$\overline\Gamma\setminus\bp$ incident to $p$ are ends of $\Gamma$ (see Figure
\ref{fig11}(a,b), where the ends are shown by solid and dashed lines).
Furthermore, we assume that the weights of the edges of $T=h_*(C^*)$ adjacent to $x$
and covered by the edges of the unbounded components of $\Gamma\setminus\bp$, satisfy $m_1,m_2>1$.
The case of $m_1=1$ or $m_2=1$ can be treated in the same way.

Suppose, first, that exactly one edge of $\overline\Gamma\setminus\bp$ incident to $p$
is an end of $\Gamma$ (see Figure \ref{fig11}(a)).
Without loss of generality, we can
also suppose that the germ of the path $\bx^{(t)}\in\R^{2n}$, $t\in(\R,t^*)$, is such that
the point $x=h(p)$ moves along a line transversal to the $h$-image of the bounded component
of $\Gamma\setminus\bp$ for the original curve $C^*$, while the rest of
$\bx$ stay fixed. Then we observe the deformations of $C^*$ depicted in Figures
\ref{fig11}(c$_1$-c$_3$,d$_1$-d$_3$), respectively corresponding to $\Ev_n$-projections to
the halfspaces $H_+$ and $H_-$.

Introduce the following parameters of the curve $C^*$ shown in Figure \ref{fig11}(a$_1$):
$$m_1=\mt([v_1,v_2]),\quad\nu_1=\frac{\mu(\Gamma,h,V_1)}{m_1},\quad\nu'_i=\frac{
\mu(\Gamma,h,p)}{m_i},\ i=1,2\ ,$$ where
$v_1=h(V_1)$, $V_1\in\Gamma^0$. The constancy of $\GS_y(\Delta,1,(n_v,n_e),\bx(t))$, $t\in(\R,t^*)$, amounts
to the verification of the relation
\begin{equation}c_1+c_2+c_3=d_1+d_2+d_3\ ,\label{e14b}\end{equation} where
$c_1\varphi(y),c_2\varphi(y),c_3\varphi(y)$ and $d_1
\varphi(y),d_2\varphi(y),d_3
\varphi(y)$ are the contributions of the curves shown in Figures \ref{fig11}(c$_1$-c$_3$) and \ref{fig11}(d$_1$-d$_3$), respectively,
with $\varphi(y)$ some expression and (setting $y=z^2$)
$$c_1+c_2=\Psi^{(2)}_z(m_1,\nu_1,\nu'_1)(z-z^{-1})^3(z+z^{-1})\ ,$$
$$c_3-d_1\overset{\text{(\ref{ea3})}}{=}(z^{\nu_1m_1}-z^{-\nu_1m_1})\Bigg[m_1\frac{z^{\nu'_1m_1-m_1}-
z^{m_1-\nu'_1m_1}}{z^{m_1}-z^{-m_1}}-m_2\frac{z^{\nu'_2m_2-m_2}-
z^{m_2-\nu'_2m_2}}{z^{m_2}-z^{-m_2}}$$
$$+\frac{z^{\nu'_1m_1-\nu'_1}-z^{\nu'_1-\nu'_1m_1}}{z^{\nu'_1}-z^{-\nu'_1}}
-\frac{z^{\nu'_2m_2-\nu'_2}-z^{\nu'_2-\nu'_2m_2}}{z^{\nu'_2}-z^{-\nu'_2}}\Bigg]\ ,$$
$$d_2\overset{\text{(\ref{e14a})}}{=}(z^{\nu'_1m_1}-z^{-\nu'_1m_1})\Bigg[(m_1-1)(z^{\nu_1m_1}
+z^{-\nu_1m_1})-
2\frac{z^{\nu_1m_1-\nu_1}-z^{\nu_1-\nu_1m_1}}{z^{\nu_1}-z^{-\nu_1}}\Bigg]\ ,$$
$$d_3=(z^{\nu_1m_1}-z^{\nu_1m_1})\Psi^{(1)}_z(m_2,\nu'_2)(z-z^{-1})^2(z+z^{-1})\ .$$
Substituting these formulas to (\ref{e14b}) and using expressions (\ref{e9}), (\ref{e9a})
for $\Psi^{(2)}$ and $\Psi^{(1)}$, we immediately establish the validity of (\ref{e14b}).

Suppose now that two edges of $\overline\Gamma\setminus\bp$ incident to $p$
are ends of $\Gamma$ (see Figure \ref{fig11}(b)).
We have to consider deformations of $C^*$ shown in Figures \ref{fig11}(e$_1$-e$_4$)
labeled according as the
$\Ev_n$-image of the deformed curves belongs to the halfspace $H_+$ or $H_-$. The constancy of
$\GS_y(\Delta,1,(n_v,n_e),\bx(t))$, $t\in(\R,t^*)$, reduces to the relation
\begin{equation}e_1+e_2=e_3+e_4\ ,\label{e14c}\end{equation}
where $e_i\varphi(y)$ is the contribution of the curves shown in Figure \ref{fig11}(e$_i$), $i=1,2,3,4$,
with some expression $\varphi(y)$ and
$$e_1-e_3\overset{\text{(\ref{ea3})}}{=}m_1\frac{z^{\nu'_1m_1-m_1}-
z^{m_1-\nu'_1m_1}}{z^{m_1}-z^{-m_1}}-m_2\frac{z^{\nu'_2m_2-m_2}-
z^{m_2-\nu'_2m_2}}{z^{m_2}-z^{-m_2}}$$
$$+\frac{z^{\nu'_1m_1-\nu'_1}-z^{\nu'_1-\nu'_1m_1}}{z^{\nu'_1}-z^{-\nu'_1}}
-\frac{z^{\nu'_2m_2-\nu'_2}-z^{\nu'_2-\nu'_2m_2}}{z^{\nu'_2}-z^{-\nu'_2}}\ ,$$
$$e_2=\Psi^{(1)}_z(m_1,\nu'_1)(z-z^{-1})^2(z+z^{-1}),\quad e_4=\Psi^{(1)}_z(m_2,\nu'_2)(z-z^{-1})^2(z+z^{-1})\ .$$
Substituting these formulas to (\ref{e14c}) and using expression (\ref{e9a}) for
$\Psi^{(1)}$, we immediately establish the validity of (\ref{e14c}).

\subsection{Degeneration of type (2viii)}
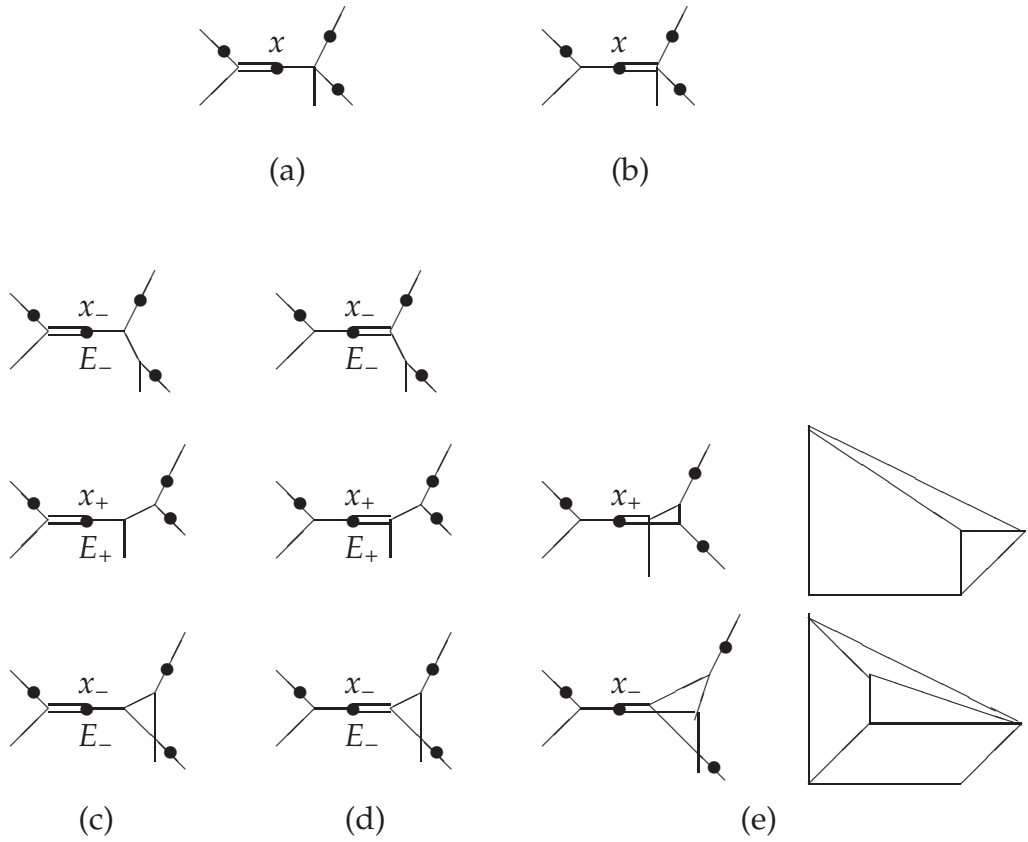
\begin{figure}
\setlength{\unitlength}{1cm}
\begin{picture}(13,11)(0,0)
\thinlines

\put(2.5,9.6){\line(1,1){0.5}}\put(3,10.1){\line(-1,1){0.5}}
\put(3,10.15){\line(1,0){0.5}}\put(3,10.05){\line(1,0){0.5}}
\put(3.5,10.1){\line(1,0){0.5}}\put(4,10.1){\line(1,-1){0.5}}
\put(4,10.1){\line(0,-1){0.5}}\put(4,10.1){\line(1,2){0.4}}

\put(7,9.6){\line(1,1){0.5}}\put(7.5,10.1){\line(-1,1){0.5}}
\put(8,10.15){\line(1,0){0.5}}\put(8,10.05){\line(1,0){0.5}}
\put(7.5,10.1){\line(1,0){0.5}}\put(8.5,10.1){\line(1,-1){0.5}}
\put(8.5,10.1){\line(0,-1){0.5}}\put(8.5,10.1){\line(1,2){0.4}}

\put(2.7,10.2){$\bullet$}\put(7.2,10.2){$\bullet$}
\put(3.4,9.97){$\bullet$}\put(7.9,9.97){$\bullet$}
\put(4.1,10.4){$\bullet$}\put(8.6,10.4){$\bullet$}
\put(4.2,9.7){$\bullet$}\put(8.7,9.7){$\bullet$}

\put(3.4,8.6){(a)}\put(7.9,8.6){(b)}

\put(3.4,10.3){$x$}\put(7.9,10.3){$x$}

\put(0,6.1){\line(1,1){0.5}}\put(0.5,6.6){\line(-1,1){0.5}}
\put(0.5,6.65){\line(1,0){0.5}}\put(0.5,6.55){\line(1,0){0.5}}
\put(1,6.6){\line(1,0){0.5}}\put(1.7,6.2){\line(1,-1){0.4}}
\put(1.7,6.2){\line(0,-1){0.4}}\put(1.5,6.6){\line(1,2){0.4}}
\put(1.5,6.6){\line(1,-2){0.2}}

\put(0,3.6){\line(1,1){0.5}}\put(0.5,4.1){\line(-1,1){0.5}}
\put(0.5,4.15){\line(1,0){0.5}}\put(0.5,4.05){\line(1,0){0.5}}
\put(1,4.1){\line(1,0){0.5}}\put(1.9,4.3){\line(1,-1){0.4}}
\put(1.5,4.1){\line(0,-1){0.5}}\put(1.9,4.3){\line(1,2){0.4}}
\put(1.5,4.1){\line(2,1){0.4}}

\put(0,1.1){\line(1,1){0.5}}\put(0.5,1.6){\line(-1,1){0.5}}
\put(0.5,1.65){\line(1,0){0.5}}\put(0.5,1.55){\line(1,0){0.5}}
\put(1,1.6){\line(1,0){0.5}}\put(1.5,1.6){\line(1,-1){0.8}}
\put(1.9,1.8){\line(0,-1){0.9}}\put(1.9,1.8){\line(1,2){0.4}}
\put(1.5,1.6){\line(2,1){0.4}}

\put(3.5,6.1){\line(1,1){0.5}}\put(4,6.6){\line(-1,1){0.5}}
\put(4.5,6.65){\line(1,0){0.5}}\put(4.5,6.55){\line(1,0){0.5}}
\put(4,6.6){\line(1,0){0.5}}\put(5.2,6.2){\line(1,-1){0.4}}
\put(5.2,6.2){\line(0,-1){0.4}}\put(5,6.6){\line(1,2){0.4}}
\put(5,6.6){\line(1,-2){0.2}}

\put(3.5,3.6){\line(1,1){0.5}}\put(4,4.1){\line(-1,1){0.5}}
\put(4.5,4.15){\line(1,0){0.5}}\put(4.5,4.05){\line(1,0){0.5}}
\put(4,4.1){\line(1,0){0.5}}\put(5.4,4.3){\line(1,-1){0.4}}
\put(5,4.1){\line(0,-1){0.5}}\put(5.4,4.3){\line(1,2){0.4}}
\put(5,4.1){\line(2,1){0.4}}

\put(3.5,1.1){\line(1,1){0.5}}\put(4,1.6){\line(-1,1){0.5}}
\put(4.5,1.65){\line(1,0){0.5}}\put(4.5,1.55){\line(1,0){0.5}}
\put(4,1.6){\line(1,0){0.5}}\put(5,1.6){\line(1,-1){0.8}}
\put(5.4,1.8){\line(0,-1){0.9}}\put(5.4,1.8){\line(1,2){0.4}}
\put(5,1.6){\line(2,1){0.4}}

\put(0.2,6.7){$\bullet$}\put(0.9,6.47){$\bullet$}
\put(0.2,4.2){$\bullet$}\put(0.9,3.97){$\bullet$}
\put(0.2,1.7){$\bullet$}\put(0.9,1.47){$\bullet$}
\put(3.7,6.7){$\bullet$}\put(4.4,6.47){$\bullet$}
\put(3.7,4.2){$\bullet$}\put(4.4,3.97){$\bullet$}
\put(3.7,1.7){$\bullet$}\put(4.4,1.47){$\bullet$}

\put(1.6,6.9){$\bullet$}\put(1.8,5.9){$\bullet$}
\put(1.95,4.5){$\bullet$}\put(2,4){$\bullet$}
\put(1.95,2){$\bullet$}\put(2,0.9){$\bullet$}
\put(5.1,6.9){$\bullet$}\put(5.3,5.9){$\bullet$}
\put(5.45,4.5){$\bullet$}\put(5.5,4){$\bullet$}
\put(5.45,2){$\bullet$}\put(5.5,0.9){$\bullet$}

\put(0.9,6.8){$x_-$}\put(0.9,1.8){$x_-$}
\put(4.4,6.8){$x_-$}\put(4.4,1.8){$x_-$}
\put(0.9,4.3){$x_+$}\put(4.4,4.3){$x_+$}
\put(0.9,6.1){$E_-$}\put(0.9,1.1){$E_-$}
\put(4.4,6.1){$E_-$}\put(4.4,1.1){$E_-$}
\put(0.9,3.6){$E_+$}\put(4.4,3.6){$E_+$}

\put(0.9,0){(c)}\put(4.4,0){(d)}

\put(7,1.1){\line(1,1){0.5}}\put(7.5,1.6){\line(-1,1){0.5}}
\put(8,1.65){\line(1,0){0.4}}\put(8,1.55){\line(1,0){1}}
\put(7.5,1.6){\line(1,0){0.5}}\put(8.4,1.65){\line(1,-1){1}}
\put(9.05,1.55){\line(0,-1){0.8}}\put(9.2,2.05){\line(1,2){0.4}}
\put(8.4,1.65){\line(2,1){0.8}}\put(9.2,2.05){\line(-1,-3){0.2}}

\put(7,3.6){\line(1,1){0.5}}\put(7.5,4.1){\line(-1,1){0.5}}
\put(8,4.15){\line(1,0){0.4}}\put(8,4.05){\line(1,0){0.8}}
\put(7.5,4.1){\line(1,0){0.5}}\put(8.8,4.05){\line(1,-1){0.6}}
\put(8.4,4.15){\line(0,-1){0.8}}\put(8.8,4.3){\line(1,2){0.4}}
\put(8.4,4.1){\line(2,1){0.4}}\put(8.8,4.05){\line(0,1){0.25}}

\put(7.2,4.2){$\bullet$}\put(7.9,3.97){$\bullet$}
\put(7.2,1.7){$\bullet$}\put(7.9,1.47){$\bullet$}

\put(8.9,4.6){$\bullet$}\put(9,3.63){$\bullet$}
\put(9.3,2.3){$\bullet$}\put(9.15,0.7){$\bullet$}

\put(7.9,4.3){$x_+$}\put(7.9,1.8){$x_-$}

\put(10.5,0.6){\line(0,1){2.25}}\put(10.5,0.6){\line(1,0){2}}
\put(12.5,0.6){\line(1,1){0.8}}\put(10.5,2.8){\line(2,-1){2.75}}
\put(10.5,0.6){\line(1,1){0.8}}\put(11.3,1.4){\line(1,0){2}}
\put(10.5,2.8){\line(1,-1){0.8}}\put(11.3,1.4){\line(0,1){0.65}}
\put(11.3,2.05){\line(3,-1){2}}

\put(10.5,3.1){\line(0,1){2.25}}\put(10.5,3.1){\line(1,0){2}}
\put(12.5,3.1){\line(1,1){0.85}}\put(10.5,5.35){\line(2,-1){2.8}}
\put(12.5,3.1){\line(0,1){0.85}}\put(12.5,3.95){\line(1,0){0.85}}
\put(10.5,5.3){\line(3,-2){2}}

\put(9.6,0){(e)}

\end{picture}
\caption{Degeneration of type (2viii)}\label{fig8}
\end{figure}

Let $C^*$ be as in Lemma \ref{l6}(2viii) (see Figure \ref{fig8}(a,b)). We simultaneously consider
the degenerations of types (2viii-a) and (2viii-b) having the same image
$(\overline\Gamma',\bp,h')$ in $\overline{\mathcal M}^{\;e}_{0,(n_v-1,n_e+1)}(\R^2,\Delta)$
(see Section \ref{sec2.3}).
Let $p\in\bp$ be the collinear trivalent vertex, and $E\in(\Gamma')^1$ the edge of
$\Gamma'$ containing $p$.
Note that the only one edge of $\Gamma$ incident to the endpoints of $E$ has
outward canonical orientation;
furthermore, without loss of generality we can suppose that the other edges of $\Gamma$
incident to the endpoints of $E$ contain marked points (see Figure \ref{fig8}(a,b)).
We will verify the
invariance of $\GS_y(\Delta,1,(n_v,n_e),\bx^{(t)})$,
$t\in(\R,t^*)$, simultaneously considering deformations of the type (2viii-a) and (2viii-b).
We can assume that, along the path $\bx^{(t)}$, $t\in(\R,t^*)$,
the image $x=h(p)$ moves transversally to the edge $h'(E)$, whereas
$h(\bp\setminus\{p\})=h'(\bp\setminus\{p\})$ stays fixed.
The three possible deformations of $C^*$ of type (2viii-a) correspond to splittings of the
four-valent vertex
outside the collinear cycle into a pair of trivalent vertices (see Figure \ref{fig8}(c)).
A curve $C^*$ of type (2viii-b) similarly admits three deformations shown in Figure
\ref{fig8}(d) as well as
two more types shown in Figure \ref{fig8}(e). In Figures \ref{fig8}(c,d,e) we present
the subdivisions of the quadrangle $Q={\mathcal D}(V_4)$, $V_4$ being the four-valent
endpoint of $E$, and labeled
with $H_+$, $H_-$ according to the move of $h(p)$ into one or the other halfplane
bounded by the line passing through $h'(E)$. We observe that the $H_+,H_-$-labels of the deformations
shown in Figures \ref{fig8}(c,d) meet the rule described in Section \ref{sec-b}. We also denote by
$\bx_+$, $\bx_-$ the images of $\bp$ in the plane corresponding to deformations labeled by
$H_+$, $H_-$, respectively.

We need to show that the total refined weight of the curves labeled by $H_+$
equals that of the curves labeled by $H_-$. Note that we can restrict our attention to only fragments presented in
Figures \ref{fig8}(c,d,e), since the remaining part of these curves provides the same
multiplicative contribution to each term. Furthermore, we can make two additional assumptions
\begin{itemize}\item the outgoing edge
of the trivalent endpoint $V_3$ of $E$ has weight one; indeed,
the refined weights of all curves depend only
on the Mikhalkin weight of this vertex of $\Gamma'$ and on the weight of $E$,
and we can vary the weights and slopes of
the images of the edges incident to $V$ and different from $E^*$ while keeping the
aforementioned parameters;
\item the marked point $p$ can be chosen arbitrarily close to $V_3$.
\end{itemize}

To prove the required equality,
we do not perform direct tedious computations
but will exhibit a deformation of the configuration $\bx$ along which the $H_+$-types turn into the $H_-$-types so that,
in this way, the considered curves undergo degenerations of types (2i)-(2vii) for which we
have already established the invariance of
$\GS_y(\Delta,1,(n_v,n_e),\bx)$.

For each of the curves presented in figures \ref{fig8}(c,d,e), we denote by $E_+$ or
$E_-$ the edge of the
embedded plane tropical curve $T=h_*(\overline\Gamma)$, containing $h(p)$, so that the
subindex $\pm$
matches the $H_+,H_-$-labels. Next, we deform the configuration $\bx_\eps$, $\eps=\pm$,
so that the point
$x=h(p)$ moves along the line, containing the edge $E_\eps$, to the position $x'$, also in a small neighborhood
of $V$ (see Figure \ref{fig9}(a)), while the rest of $\bx$ stays fixed. Denote by
$\bx'_\eps$ the deformed configuration. Now, for the initial configurations $\bx_\pm$ we introduce
the curves shown in
Figure \ref{fig9}(c), where the triangular fragments are dual to the subdivisions of the ${\mathcal D}(V_3)$, induced by the choice of each of the integral points
$v\in\Int{\mathcal D}(V_3)$
(see Figure \ref{fig9}(b)), and for the configurations $\bx'_\pm$, we introduce the curves shown in
Figures \ref{fig9}(d,e,f).

\begin{figure}
\setlength{\unitlength}{1cm}
\begin{picture}(14,19)(0,0)
\thinlines

\put(0,1.2){\line(1,1){0.5}}\put(0.5,2.1){\line(-1,1){0.5}}
\put(0.5,2.1){\line(1,0){0.9}}\put(0.5,1.7){\line(1,0){1.4}}
\put(1.4,2.1){\line(1,-1){1}}\put(0.5,1.7){\line(0,1){0.4}}
\put(1.95,1.7){\line(0,-1){0.7}}\put(2.2,2.5){\line(1,2){0.4}}
\put(1.4,2.1){\line(2,1){0.8}}\put(2.2,2.5){\line(-1,-3){0.27}}

\put(3.5,1.2){\line(1,1){0.5}}\put(4,2.1){\line(-1,1){0.5}}
\put(4,2.1){\line(1,0){1.55}}\put(4,1.7){\line(1,0){0.85}}
\put(5.55,2.1){\line(0,-1){1.4}}\put(4,1.7){\line(0,1){0.4}}
\put(4.9,1.7){\line(1,-1){1.1}}\put(5.7,2.5){\line(1,2){0.4}}
\put(5.7,2.5){\line(-1,-1){0.8}}\put(5.7,2.5){\line(-1,-3){0.15}}

\put(7,1.2){\line(1,1){0.5}}\put(7.5,2.1){\line(-1,1){0.5}}
\put(7.5,1.7){\line(0,1){0.4}}\put(7.5,2.1){\line(1,0){0.9}}
\put(7.5,1.7){\line(1,0){1.3}}\put(8.8,1.7){\line(1,-1){0.6}}
\put(8.4,2.1){\line(0,-1){1}}\put(8.8,2.3){\line(1,2){0.4}}
\put(8.4,2.1){\line(2,1){0.4}}\put(8.8,1.7){\line(0,1){0.6}}

\put(10.5,1.2){\line(1,1){0.5}}\put(11,2.1){\line(-1,1){0.5}}
\put(11,1.7){\line(0,1){0.4}}\put(11,2.1){\line(1,0){1.5}}
\put(11,1.7){\line(1,0){0.7}}\put(12.5,2.1){\line(1,-1){1}}
\put(11.7,1.7){\line(0,-1){0.5}}\put(12.5,2.5){\line(1,2){0.4}}
\put(12.5,2.5){\line(-1,-1){0.8}}\put(12.5,2.1){\line(0,1){0.4}}

\put(7.2,2.2){$\bullet$}\put(7.4,1.57){$\bullet$}
\put(10.7,2.2){$\bullet$}\put(10.9,1.57){$\bullet$}
\put(0.2,2.2){$\bullet$}\put(0.4,1.6){$\bullet$}
\put(3.7,2.2){$\bullet$}\put(3.9,1.6){$\bullet$}

\put(8.9,2.6){$\bullet$}\put(9,1.28){$\bullet$}
\put(12.6,2.8){$\bullet$}\put(13.1,1.3){$\bullet$}
\put(2.3,2.75){$\bullet$}\put(2.15,1.15){$\bullet$}
\put(5.8,2.75){$\bullet$}\put(5.6,0.8){$\bullet$}

\put(1,0){(f$_1$)\ $H_-$}\put(4.5,0){(f$_2$)\ $H_-$}
\put(8,0){(f$_3$)\ $H_+$}\put(11.5,0){(f$_4$)\ $H_+$}

\put(1,5.3){\line(1,1){0.5}}\put(1.5,6.2){\line(-1,1){0.5}}
\put(1.5,6.2){\line(1,0){1}}\put(1.5,5.8){\line(1,0){1}}
\put(1.5,5.8){\line(0,1){0.4}}
\put(2.7,5.4){\line(1,-1){0.4}}
\put(2.7,5.4){\line(0,-1){0.4}}\put(2.5,6.2){\line(1,2){0.4}}
\put(2.5,5.8){\line(1,-2){0.2}}\put(2.5,5.8){\line(0,1){0.4}}

\put(5,5.3){\line(1,1){0.5}}\put(5.5,6.2){\line(-1,1){0.5}}
\put(5.5,5.8){\line(0,1){0.4}}
\put(5.5,6.2){\line(1,0){0.9}}\put(5.5,5.8){\line(1,0){0.5}}
\put(6.8,6.4){\line(1,-1){1.2}}
\put(6,5.8){\line(0,-1){0.5}}\put(6.8,6.4){\line(1,2){0.4}}
\put(6.4,6.2){\line(2,1){0.4}}\put(6,5.8){\line(1,1){0.4}}

\put(9,5.3){\line(1,1){0.5}}\put(9.5,6.2){\line(-1,1){0.5}}
\put(9.5,6.2){\line(1,0){1}}\put(9.5,5.8){\line(1,0){1}}
\put(9.5,5.8){\line(0,1){0.4}}\put(10.5,5.8){\line(0,1){0.4}}
\put(10,5.8){\line(1,0){0.5}}\put(10.5,5.8){\line(1,-1){1}}
\put(10.9,6.4){\line(0,-1){1.5}}\put(10.9,6.4){\line(1,2){0.4}}
\put(10.5,6.2){\line(2,1){0.4}}

\put(1.4,5.65){$\bullet$}\put(1.2,6.3){$\bullet$}
\put(2.6,6.5){$\bullet$}\put(2.8,5.1){$\bullet$}
\put(5.4,5.65){$\bullet$}\put(5.2,6.3){$\bullet$}
\put(6.9,6.7){$\bullet$}\put(7.6,5.4){$\bullet$}
\put(9.4,5.65){$\bullet$}\put(9.2,6.3){$\bullet$}
\put(10.9,6.5){$\bullet$}\put(11,5.1){$\bullet$}

\put(1.6,4.3){(e$_1$)\ $H_-$}\put(5.8,4.3){(e$_2$)\ $H_+$}\put(9.6,4.3){(e$_3$)\ $H_-$}

\put(1,9.3){\line(1,1){0.5}}\put(1.5,10.2){\line(-1,1){0.5}}
\put(1.5,10.2){\line(2,-1){0.4}}\put(1.5,9.8){\line(2,1){0.4}}
\put(1.5,9.8){\line(0,1){0.4}}\put(1.9,10){\line(1,0){0.6}}
\put(2.7,9.6){\line(1,-1){0.4}}
\put(2.5,10){\line(1,2){0.4}}
\put(2.5,10){\line(1,-2){0.2}}\put(2.7,9.2){\line(0,1){0.4}}

\put(5,9.3){\line(1,1){0.5}}\put(5.5,10.2){\line(-1,1){0.5}}
\put(5.5,10.2){\line(2,-1){0.4}}\put(5.5,9.8){\line(2,1){0.4}}
\put(5.5,9.8){\line(0,1){0.4}}\put(5.9,10){\line(1,0){0.5}}
\put(6.8,10.2){\line(1,-1){1.2}}
\put(6.4,10){\line(0,-1){0.8}}\put(6.8,10.2){\line(1,2){0.4}}
\put(6.4,10){\line(2,1){0.4}}

\put(9,9.3){\line(1,1){0.5}}\put(9.5,10.2){\line(-1,1){0.5}}
\put(9.5,10.2){\line(2,-1){0.4}}\put(9.5,9.8){\line(2,1){0.4}}
\put(9.5,9.8){\line(0,1){0.4}}\put(9.9,10){\line(1,0){0.6}}
\put(10.5,10){\line(1,-1){1}}
\put(10.9,10.2){\line(0,-1){1}}\put(10.9,10.2){\line(1,2){0.4}}
\put(10.5,10){\line(2,1){0.4}}

\put(1.4,9.65){$\bullet$}\put(1.2,10.3){$\bullet$}
\put(2.6,10.3){$\bullet$}\put(2.8,9.3){$\bullet$}
\put(5.4,9.65){$\bullet$}\put(5.2,10.3){$\bullet$}
\put(6.9,10.5){$\bullet$}\put(7.6,9.2){$\bullet$}
\put(9.4,9.65){$\bullet$}\put(9.2,10.3){$\bullet$}
\put(10.9,10.3){$\bullet$}\put(11,9.3){$\bullet$}

\put(1.6,8.3){(d$_1$)\ $H_-$}\put(5.8,8.3){(d$_2$)\ $H_+$}\put(9.6,8.3){(d$_3$)\ $H_-$}

\put(1,13.3){\line(1,1){0.5}}\put(1.5,14.2){\line(-1,1){0.5}}
\put(1.5,14.2){\line(2,-1){0.4}}\put(1.5,13.8){\line(2,1){0.4}}
\put(1.5,13.8){\line(0,1){0.4}}\put(1.9,14){\line(1,0){0.6}}
\put(2.7,13.6){\line(1,-1){0.4}}
\put(2.5,14){\line(1,2){0.4}}
\put(2.5,14){\line(1,-2){0.2}}\put(2.7,13.2){\line(0,1){0.4}}

\put(5,13.3){\line(1,1){0.5}}\put(5.5,14.2){\line(-1,1){0.5}}
\put(5.5,14.2){\line(2,-1){0.4}}\put(5.5,13.8){\line(2,1){0.4}}
\put(5.5,13.8){\line(0,1){0.4}}\put(5.9,14){\line(1,0){0.5}}
\put(6.8,14.2){\line(1,-1){1.2}}
\put(6.4,14){\line(0,-1){0.8}}\put(6.8,14.2){\line(1,2){0.4}}
\put(6.4,14){\line(2,1){0.4}}

\put(9,13.3){\line(1,1){0.5}}\put(9.5,14.2){\line(-1,1){0.5}}
\put(9.5,14.2){\line(2,-1){0.4}}\put(9.5,13.8){\line(2,1){0.4}}
\put(9.5,13.8){\line(0,1){0.4}}\put(9.9,14){\line(1,0){0.6}}
\put(10.5,14){\line(1,-1){1}}
\put(10.9,14.2){\line(0,-1){1}}\put(10.9,14.2){\line(1,2){0.4}}
\put(10.5,14){\line(2,1){0.4}}

\put(1.8,13.85){$\bullet$}\put(1.2,14.3){$\bullet$}
\put(2.6,14.3){$\bullet$}\put(2.8,13.3){$\bullet$}
\put(5.8,13.85){$\bullet$}\put(5.2,14.3){$\bullet$}
\put(6.9,14.5){$\bullet$}\put(7.6,13.2){$\bullet$}
\put(9.8,13.85){$\bullet$}\put(9.2,14.3){$\bullet$}
\put(10.9,14.3){$\bullet$}\put(11,13.3){$\bullet$}

\put(1.6,12.3){(c$_1$)\ $H_-$}\put(5.8,12.3){(c$_2$)\ $H_+$}\put(9.6,12.3){(c$_3$)\ $H_-$}

\put(3,17){\line(1,1){1}}\put(4,18){\line(1,0){2}}\put(3,19){\line(1,-1){1}}
\dashline{0.2}(3,18)(4,18)
\put(3.2,17.85){$\bullet$}\put(4.6,17.85){$\bullet$}
\put(3.1,18.2){$x'_{\pm}$}\put(4.5,18.2){$x_{\pm}$}
\put(5.5,17.5){$E_{\pm}$}
\put(4.4,16){(a)}

\put(7,18){\line(1,1){1}}\put(7,18){\line(2,-1){2}}
\put(8,19){\line(1,-2){1}}\put(7,18){\line(1,0){1}}
\put(8,18){\line(1,-1){1}}\put(8,18){\line(0,1){1}}
\put(7.7,18.1){$v$}
\put(8,16){(b)}

\put(1.8,13.5){$x_-$}\put(5.8,13.5){$x_+$}\put(9.8,13.5){$x_-$}
\put(1.7,9.5){$x'_-$}\put(5.7,9.5){$x'_+$}\put(9.7,9.5){$x'_-$}
\put(1.5,5.4){$x'_-$}\put(5.5,5.4){$x'_+$}\put(9.5,5.4){$x'_-$}
\put(0.5,1.3){$x'_-$}\put(4,1.3){$x'_-$}\put(7.5,1.3){$x'_+$}\put(11,1.3){$x'_+$}

\end{picture}
\caption{Deformations in case (2viii), I}\label{fig9}
\end{figure}

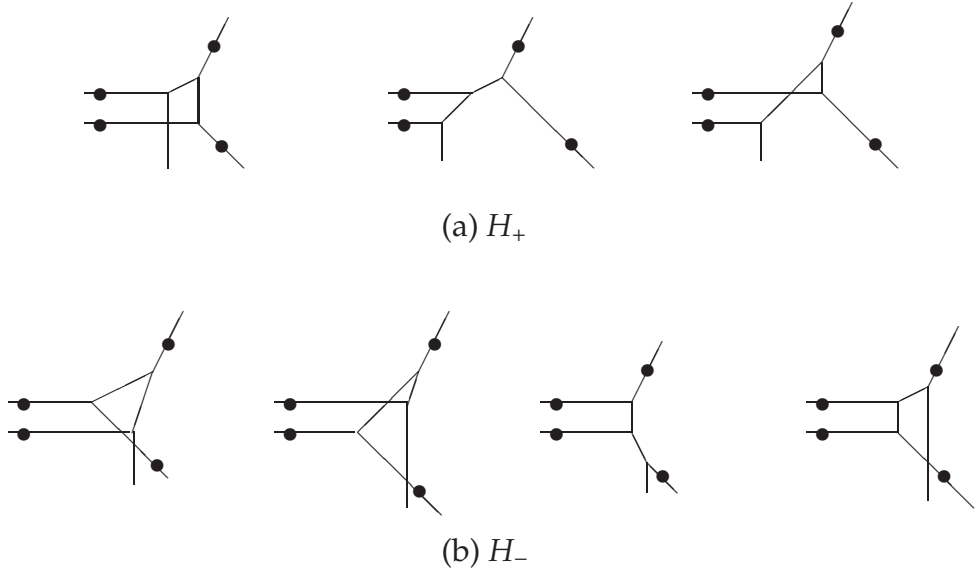
\begin{figure}
\setlength{\unitlength}{1cm}
\begin{picture}(14,8)(0,0)
\thinlines

\put(0.3,2.1){\line(1,0){1.1}}\put(0.3,1.7){\line(1,0){1.6}}
\put(1.4,2.1){\line(1,-1){1}}
\put(1.95,1.7){\line(0,-1){0.7}}\put(2.2,2.5){\line(1,2){0.4}}
\put(1.4,2.1){\line(2,1){0.8}}\put(2.2,2.5){\line(-1,-3){0.27}}

\put(3.8,2.1){\line(1,0){1.75}}\put(3.8,1.7){\line(1,0){1.05}}
\put(5.55,2.1){\line(0,-1){1.4}}
\put(4.9,1.7){\line(1,-1){1.1}}\put(5.7,2.5){\line(1,2){0.4}}
\put(5.7,2.5){\line(-1,-1){0.8}}\put(5.7,2.5){\line(-1,-3){0.15}}

\put(7.3,2.1){\line(1,0){1.2}}\put(7.3,1.7){\line(1,0){1.2}}
\put(8.7,1.3){\line(1,-1){0.4}}
\put(8.7,1.3){\line(0,-1){0.4}}\put(8.5,2.1){\line(1,2){0.4}}
\put(8.5,1.7){\line(1,-2){0.2}}\put(8.5,1.7){\line(0,1){0.4}}

\put(10.8,2.1){\line(1,0){1.2}}\put(10.8,1.7){\line(1,0){1.2}}
\put(12,1.7){\line(0,1){0.4}}
\put(11.5,1.7){\line(1,0){0.5}}\put(12,1.7){\line(1,-1){1}}
\put(12.4,2.3){\line(0,-1){1.5}}\put(12.4,2.3){\line(1,2){0.4}}
\put(12,2.1){\line(2,1){0.4}}

\put(7.4,1.96){$\bullet$}\put(7.4,1.56){$\bullet$}
\put(10.9,1.96){$\bullet$}\put(10.9,1.56){$\bullet$}
\put(0.4,1.96){$\bullet$}\put(0.4,1.56){$\bullet$}
\put(3.9,1.96){$\bullet$}\put(3.9,1.56){$\bullet$}

\put(2.9,6.7){$\bullet$}\put(3,5.38){$\bullet$}
\put(11.1,6.9){$\bullet$}\put(11.6,5.4){$\bullet$}
\put(2.3,2.75){$\bullet$}\put(2.15,1.15){$\bullet$}
\put(5.8,2.75){$\bullet$}\put(5.6,0.8){$\bullet$}

\put(6,0){(b)\ $H_-$}

\put(1.3,6.2){\line(1,0){1.1}}
\put(1.3,5.8){\line(1,0){1.5}}\put(2.8,5.8){\line(1,-1){0.6}}
\put(2.4,6.2){\line(0,-1){1}}\put(2.8,6.4){\line(1,2){0.4}}
\put(2.4,6.2){\line(2,1){0.4}}\put(2.8,5.8){\line(0,1){0.6}}

\put(5.3,6.2){\line(1,0){1.1}}\put(5.3,5.8){\line(1,0){0.7}}
\put(6.8,6.4){\line(1,-1){1.2}}
\put(6,5.8){\line(0,-1){0.5}}\put(6.8,6.4){\line(1,2){0.4}}
\put(6.4,6.2){\line(2,1){0.4}}\put(6,5.8){\line(1,1){0.4}}

\put(9.3,6.2){\line(1,0){1.7}}
\put(9.3,5.8){\line(1,0){0.9}}\put(11,6.2){\line(1,-1){1}}
\put(10.2,5.8){\line(0,-1){0.5}}\put(11,6.6){\line(1,2){0.4}}
\put(11,6.6){\line(-1,-1){0.8}}\put(11,6.2){\line(0,1){0.4}}

\put(1.4,5.65){$\bullet$}\put(1.4,6.06){$\bullet$}
\put(8.6,2.4){$\bullet$}\put(8.8,1){$\bullet$}
\put(5.4,5.65){$\bullet$}\put(5.4,6.06){$\bullet$}
\put(6.9,6.7){$\bullet$}\put(7.6,5.4){$\bullet$}
\put(9.4,5.65){$\bullet$}\put(9.4,6.06){$\bullet$}
\put(12.4,2.4){$\bullet$}\put(12.5,1){$\bullet$}

\put(6,4.3){(a)\ $H_+$}

\end{picture}
\caption{Deformations in case (2viii), II}\label{fig10}
\end{figure}

Then we introduce the following quantities:
\begin{itemize}\item $\Sig^{(1)}_\eps$, the total refined weight of the curves of the shape shown
in Figure \ref{fig8}(c,d), labeled by $H_\eps$, and projected to the same
rational curve $(\overline\Gamma',\bp,h')\in{\mathcal M}^e_{0,4}(\R^2,\Delta)$ with
$h'(\bp)=\bx_\eps$, $\eps=\pm$,
\item $\Sig^{(2)}_\eps$, the total refined weight of the
elliptic curves of shapes shown
in Figure \ref{fig8}(e), labeled by $H_\eps$, having the given degree, and passing through the
configuration $\bx_\eps$, $\eps=\pm$,
\item $\Sig^{(3)}_\eps$, the total refined weight of the elliptic curves of the shape shown
in Figure \ref{fig9}(c), labeled by $H_\eps$, having the given degree, and passing
through $\bx_\eps$, $\eps=\pm$,
\item $\Sig^{(4)}_\eps$, the total refined weight of the elliptic curves of the shape shown
in Figure \ref{fig9}(d), labeled by $H_\eps$, having the given degree, and passing
through $\bx'_\eps$, $\eps=\pm$,
\item $\Sig^{(5)}_\eps$, the total refined weight of the elliptic curves of the shape shown
in Figure \ref{fig9}(e), labeled by $H_\eps$, having the given degree, and passing
through $\bx'_\eps$, $\eps=\pm$,
\item $\Sig^{(6)}_\eps$, the total refined weight of the elliptic curves of the shape shown
in Figure \ref{fig9}(e), labeled by $H_\eps$, having the given degree, and passing
through $\bx'_\eps$, $\eps=\pm$.
\end{itemize}
The required invariance reads
\begin{equation}\Sig^{(1)}_++\Sig^{(2)}_+=\Sig^{(1)}_-+\Sig^{(2)}_-\ .
\label{e-inv}\end{equation}
On the other hand, according to the result of Section \ref{sec3.6},
$$\Sig^{(1)}_\eps=\Sig^{(4)}_\eps+\Sig^{(5)}_\eps-\Sig^{(3)}_\eps,\quad\eps=\pm$$
(recall that, by our assumption, the weight of the left lower unbounded edge equals one),
according to the result of Section \ref{sec3.3},
$$\Sig^{(2)}_\eps=\Sig^{(6)}_\eps,\quad\eps=\pm\ ,$$
according to the result of Section \ref{sec-b},
$$\Sig^{(3)}_+=\Sig^{(3)}_-\quad\text{and}\quad\Sig^{(4)}_+=\Sig^{(4)}_-\ ,$$
and hence (\ref{e-inv}) converts to the following relation:
\begin{equation}\Sig^{(5)}_++\Sig^{(6)}_+=\Sig^{(5)}_-+\Sig^{(6)}_-\ .
\label{e-inv1}\end{equation}
If we fix weights of the horizontal edges for curves shown in Figure \ref{fig9}(e,f),
equality (\ref{e-inv1}) reduces to the comparison of the refined
weight of the rational curves shown in Figure
\ref{fig10}(a) on one side, and of the refined weight of the rational curves shown in
Figure \ref{fig10}(b) (in both figures we assume that the marked points are far away from the vertices).
However, it amounts to 
the comparison of the two values of $\BG_y$ for rational curves of the same degree
hitting two generic configurations of points in $\R^2$, and hence the equality by Proposition
\ref{p1} (alternatively, it follows from the consideration of degenerations in Sections
\ref{sec-a} and \ref{sec-b}).

\section{Computations} 
\label{sec-com}
Here we explicitly answer the question on what are elliptic broccoli curves. Namely, we show that
the refined elliptic broccoli invariant $\GS_y(\Delta,1,(n_v,n_e))$ is a Laurent polynomial
in $y$ (Proposition \ref{p5}), evaluate it at $y=1$ and $y=-1$, and using the latter evaluation

\subsection{Evaluation at $y=1$ and $y=-1$, and the Laurent property}\label{sec-local}

\begin{lemma}\label{l8} (1) The functions $[m]_{z^2}^{\pm}$, $m^*_{z^2}$ and $\Psi^{(1)}_z(m,\nu)$,
$\Psi^{(2)}_z(m,\nu_1,\nu_2)$ behave regularly at $z=1$ and attain the following values
$$[m]_1^+=1,\quad\lim_{z\to1}[m]_{z^2}^-=m\quad\lim_{z\to1}[m]^*_{z^2}=\begin{cases}1,\quad&\ m\equiv0\mod2\\
\frac{1}{m},\quad&\ m\equiv1\mod2\end{cases}\ ,$$
$$\lim_{z\to1}\Psi^{(1)}_z(m,\nu)=\frac{\nu(\nu-1)m(m^2-1)}{12}\ ,$$
$$\lim_{z\to1}\Psi^{(2)}_z(m,\nu_1,\nu_2)=\frac{\nu_1\nu_2(\nu_1+\nu_2-1)m^2(m^2-1)}{12}\ .$$

(2) The functions $[m]_{z^2}^{\pm}$, $[m]^*_{z^2}$, and $\Psi^{(1)}_z(m,\nu)$,
$\Psi^{(2)}_z(m,\nu_1,\nu_2)$ reveal the following behavior in a neighborhood of $z=i$:
\begin{enumerate}\item[(2i)] for $m\in\Z$, we have $$\lim_{z\to i}[m]_{z^2}^+
=m\cdot\lim_{z\to i}[m]^*_{z^2}=m(-1)^{(m-1)/2},\quad m\equiv1\mod2\ ,$$
$$\lim_{z\to i}(z+z^{-1})[m]_{z^2}^+=2(-1)^{m/2},\quad m\equiv0\mod2\ ,$$ $$\lim_{z\to i}[m]_{z^2}^-=(-1)^{(m-1)/2},\quad m\equiv1\mod2\ ,$$
$$\lim_{z\to i}\frac{[m]_{z^2}^-}{z+z^{-1}}=m\cdot
\lim_{z\to i}\frac{[m]^*_{z^2}}{z+z^{-1}}=(-1)^{m/2-1}m,\quad m\equiv0\mod2\ ,$$
\item[(2ii)]
$\Psi^{(1)}_z(m,\nu)$ is regular at $z=i$ if $m$ is odd; furthermore, if $\nu$ is odd one has
\begin{equation}\lim_{z\to i}\Psi^{(1)}_z(m,\nu)
=\frac{1}{4}\cdot
\begin{cases}(1-\nu)(m^2-1),\quad & \nu\equiv m\equiv 1\mod4,\\
(1-\nu)(m^2+1),\quad & \nu\equiv-m\equiv1\mod4,\\
(m-1)(\nu(m-1)-m-1),\quad & \nu\equiv -m\equiv-1\mod4,\\
(\nu-1)(m^2+1)-2\nu,\quad & \nu\equiv m\equiv-1\mod4,\end{cases}\label{ef1}\end{equation}
and if $\nu$ is even, $\Psi^{(1)}_z(m,\nu)$  vanishes at $z=i$ and
$$\lim_{z\to i}\frac{\Psi^{(1)}_z(m,\nu)}{z+z^{-1}}=\frac{(-1)^{\nu/2}}{24}
m(m^2-1)\nu(\nu-3)\ ;$$
if $m$ is even, then $\Psi^{(1)}_z$ has a pole at $z=i$, and one has
$$\lim_{z\to i}(z+z^{-1})\Psi^{(1)}_z(m,\nu)=\begin{cases}\frac{m\nu}{2},\quad&\nu\equiv0\mod2\\
(-1)^{m/2}\cdot\frac{m(\nu-1)}{2},\quad&\nu\equiv1\mod2\end{cases}$$
\item[(2iii)] the function $\Psi^{(2)}_z(m,\nu_1,\nu_2)$ is regular at $z=i$; furthermore, if $m=2m'$ is even then
$$\lim_{z\to i}\Psi^{(2)}_z(m,\nu_1,\nu_2)=\frac{(-1)^{m'(\nu_1+\nu_2)}}{2}\cdot
\begin{cases}m'(\nu_1+\nu_2),\quad & \nu_1\equiv\nu_2\equiv1\mod2,\\
m'(\nu_1+\nu_2-\nu_1\nu_2),\quad &\nu_1+\nu_2\equiv1\mod2,\\
m'(\nu_1+\nu_2-2\nu_1\nu_2),\quad & \nu_1\equiv\nu_2\equiv0\mod2,\end{cases}$$ if
$m$ is odd and $\nu_1=2\nu'_1+1$, $\nu_2=2\nu'_2+1$ are odd, then
$$\lim_{z\to i}\Psi^{(2)}_z(m,\nu_1,\nu_2)=\frac{(-1)^{m(\nu'_1+\nu'_2+1)}}{4}(m^2-1)(2\nu'_1+2\nu'_2+1)\ ,$$
if $m$ is odd, $\nu_1$ is odd, and $\nu_2$ is even, then
$\Psi^{(2)}_z(m,\nu_1,\nu_2)$ vanishes at $z=i$, and
$$\lim_{z\to i}\frac{\Psi^{(2)}_z(m,\nu_1,\nu_2)}{z+z^{-1}}=
\frac{(-1)^{(1+m(\nu_1+\nu_2))/2}}{16}\big(2m^2\nu_1\nu_2-2m\nu_2-2m^3\nu_1\nu_2+2m^2\nu_2$$
$$-2m^2\nu_1\nu_2+2m\nu_1\nu_2+m(m^2-1)\nu_2^2\big)\ ,$$
if $m=2m'+1$ is odd and $\nu_1=2\nu'_1$, $\nu_2=2\nu'_2$ are even, then $\Psi^{(2)}_z(m,\nu_1,\nu_2)$ has a double zero at
$z=i$, and
$$\lim_{z\to i}\frac{\Psi^{(2)}_z(m,\nu_1,\nu_2)}{(z+z^{-1})^2}=\frac{(-1)^{\nu'_1+\nu'_2+1}}{3}
\nu'_1\nu'_2(2\nu'_1+2\nu'_2-3)
m'(m'+1)(2m'+1)^2\ .$$
\end{enumerate}
\end{lemma}

{\bf Proof.} All formulas for $[m]_{z^2}^{\pm}$ and $[m]^*_{z^2}$
are elementary and known \cite{BG,GS}. All the statements on
$\Psi^{(1)}_z(m,\nu)$ and
$\Psi^{(2)}_z(m,\nu_1,\nu_2)$ 
can be obtained by a routine direct computation. 
\proofend

\begin{proposition}\label{p5}
Given a balanced, nondegenerate multiset
$\Delta\subset\Z^2\setminus\{0\}$
and integers $n_v>0$, $n_e\ge0$ such that
$2n_v+n_e=|\Delta|$, $n=n_v+n_e$, the refined invariant
$\GS_y(\Delta,1,(n_v,n_e))$ is a symmetric Laurent polynomial in $y$, i.e., a polynomial
of degree $|P(\Delta)\cap\Z^2|-|\Delta|-1$
in $y+y^{-1}$
with a positive leading coefficient.
\end{proposition}

{\bf Proof.} We shall show that each weight
$\GS_y(\overline\Gamma,\bp,h)$ and $\GS_y(\overline\Gamma',\bp,p,h')$ in the right-hand
side of
formula (\ref{enew1}) is a symmetric Laurent polynomial in $y$ with a positive leading
coefficient.
Substitute $y=z^2$ in the formulas for $\GS_y(\overline\Gamma,\bp,h)$ and $\GS_y(\overline\Gamma',\bp,p,h')$.

Show, first, that $\GS_{z^2}(\overline\Gamma,\bp,h)$ and $\GS_{z^2}(\overline\Gamma',
\bp,p,h')$ are even functions of $z$
(equivalent to the fact that the weights of tropical curves are functions of $y$).
Indeed, the expressions
$[\mu]_{z^2}^{\pm}$ $[m]^*_{z^2}$ are even, resp. odd functions of $z$ according as $\mu$ is odd,
resp. even. By \cite[Proposition 2.3(4)]{IM}, the total number of trivalent and univalent
vertices $V$ of $\overline\Gamma$
with even $\mu(\overline\Gamma,h,V)$ is even.
Hence, $\GS_{z^2}(\overline\Gamma,\bp,h)$ is an even function of $z$ for $(\overline
\Gamma,\bp,h)
\in{\mathcal M}'_{1,(n_v,n_e)}(\R^2,\Delta)$.
The same argument
works in the case
of $(\overline\Gamma',\bp,p,h')\in{\mathcal M}''_{0,(n_v-1,n_e+1)}(\R^2,\Delta)$
if we notice that
\begin{itemize}\item the expression $\Psi^{(2)}_z(m,\nu_1,\nu_2)$ is an even,
resp. odd function of $z$
according as the pair
$(\nu_1m,\nu_2m)$ contains an even, resp. odd number of even values,
\item if $m$ is odd, then the expression $\Psi^{(1)}_z(m,\nu)[\mu(\overline\Gamma',h',V)]_{z^2}^*$
is an even, resp. odd function of $z$
according as $\nu$ is odd, resp. even, where $V$ is the univalent vertex belonging to the edge that contains $p$,
\item if $m$ is even, then the expression $\Psi^{(1)}_z(m,\nu)[\mu(\overline\Gamma',h',V)]_{z^2}^*$
is always an even function of $z$.
\end{itemize}

Next we verify that $\GS_{z^2}(\overline\Gamma,\bp,h)$ and $\GS_{z^2}(\overline
\Gamma',\bp,p,h')$ are regular at $z=1$ and $z=i$,
and hence they are symmetric Laurent polynomials in $y$. Indeed,
by Lemma \ref{l8}(1), the functions $[\mu]_{z^2}^{\pm}$ and $[\mu]_{z^2}^*$, $\mu\in\Z$,
and the functions $\Psi^{(2)}_z
(m,\nu_1,\nu_2)$, $\Psi^{(2)}_z(m,\nu_1,\nu_2)$ are regular at $z=1$.

The study of the point $z=i$ requires more work. We shall show that the total order of poles at $z=i$ does not
exceed the total order of zeroes
at $z=i$.

Let $(\overline\Gamma,\bp,h)
\in{\mathcal M}'_{1,(n_v,n_e)}(\R^2,\Delta)$ be a curve in the right-hand side of
(\ref{enew1}), and let $\overline\Gamma^{even}\subset\overline\Gamma$ be the subgraph introduced in
Definition \ref{dd1}(1).
The vertices of $\overline\Gamma^{even}$ are either univalent, or trivalent.
Then it follows from the regularity of the position of $\bp$ in $\overline\Gamma$ that
the closure of any component of $\overline\Gamma^{even}\setminus\bp$ contains at least as many unmarked vertices
as the marked ones;
hence, by Lemma \ref{l8}(2i) all the poles at $z=i$ cancel out with zeroes.

Suppose that $(\overline\Gamma',\bp,p,h')\in{\mathcal M}''_{0,(n_v-1,n_e+1)}(\R^2,\Delta)$,
and let $\overline\Gamma^{\;
\prime,even}\subset\overline\Gamma'$ be the subgraph introduced in Definition \ref{dd1}(2).
As in the preceding paragraph, the closure of any component of $\Gamma^{\prime,even}\setminus\bp$ contains at least as many
unmarked vertices as the marked ones. Assume that $E_p$ has a finite length.
If an endpoint of $E_p$ is a bivalent vertex of $\overline\Gamma^{\;\prime,even}$, then we can
take off this vertex gluing the adjacent edges of
$\overline\Gamma^{\;\prime,even}$ into one edge and preserving the claim of the previous sentence. In this case,
we derive the cancellation of poles at $z=i$ from Lemma \ref{l8}(2i) and the statements of Lemma \ref{l8}(2ii,2iii),
where the parameter $m$ is even.
If $d=1$ or $2$ endpoints of $E_p$ are univalent vertices of
$\overline\Gamma^{\;\prime,even}$, then again
the poles at $z=i$ cancel out due to Lemma \ref{l8}(2i) and the statements of Lemma \ref{l8}(2ii,2iii),
where either the parameter $\nu$ is even, or the parameter $m$ is odd, while precisely $d$ of the parameters $\nu_1,\nu_2$
are even. Assume that $E_p$ is an end. If $E_p$ has an odd weight, but its vertex belongs to
$\overline\Gamma^{\;\prime,even}$, then the argument of the preceding paragraph applies since then
$\Psi^{(1)}_z(m,\nu)$ vanishes at $z=-i$ by Lemma \ref{l8}(2ii). If $E_p$ has an even weight and it is a separate component of $\overline\Gamma^{\;\prime,even}$, then $m=\mu(\overline\Gamma',h',V)$ is even, where 
$V\in(\overline\Gamma')^0_\infty$ is a vertex of $E_p$, and hence the regularity of 
$\GS_{z^2}(\overline
\Gamma',\bp,p,h')$ at $z=i$ follows from the argument of the preceding paragraph and the fact that 
the product
$\Psi^{(1)}_z(m,\nu)[m]^*_{z^2}$ is regular at $z=i$ (see Lemma \ref{l8}(2i,2ii)).
If $E_p$ ends at a trivalent vertex of $\overline\Gamma^{\;\prime,even}$, then we separate the two 
other adjacent edges of $\overline\Gamma^{\;\prime,even}$ from $E_p$, glue them up into one edge and
obtain the required result as in the latter considered case.

The degree of $\GS_y(\Delta,1,(n_v,n_e))$  can be computed as in \cite[Proposition 2.11]{IM} with an additional term coming from $\prod_{V\in\overline\Gamma^0_\infty}[\mu(\overline\Gamma,h,V)]^*_y$.
\proofend

\subsection{Weights of elliptic broccoli curves and tropical descendant invariants}\label{char}

\begin{proposition}\label{cor2}
Given a balanced, nondegenerate multiset $\Delta\subset\Z^2\setminus\{0\}$ and integers $n_v>0$, $n_e\ge0$ such that
$2n_v+n_e=|\Delta|$, $n=n_v+n_e$, the value $\GS_1(\Delta,1,(n_v,n_e))$ is a positive integer.
Moreover, for any generic point $\bx\in\R^{2n}$, each summand in the right-hand side of
(\ref{enew1}) evaluated at $y=1$ is a positive integer.
\end{proposition}

{\bf Proof.}
The positivity is straightforward from Lemma \ref{l8}(1) and the absence of trivalent collinear unmarked vertices.
So, we have to explain the integrality. Consider the subgraph $G$ of $\Gamma$ (resp. $\Gamma'$) consisting of
finite edges (and their endpoints) that are incident to trivalent marked collinear vertices, whose two other incident edges
are ends of the same weight. Due to the general position of $\bx$, the second endpoint of any edge of $G$ is unmarked, and,
in the case of $G\subset\Gamma'$, if an unmarked vertex of $G$ is incident to two edges of $G$, then the third edge of
$\Gamma'$ attached to this vertex does not contain the marked point $p$. Thus, the integrality follows from
Remark \ref{rr1} and the fact that the Mikhalkin weight of a trivalent vertex incident to
at least one or two even edges
is divisible by $2$ or $4$, respectively.
\proofend

\begin{remark}\label{rr2}
If $\Delta$ is primitive, 
the invariant $\GS_1(\Delta,1,(n_v,n_e))$ can be regarded as a tropical elliptic descendant invariant
$\big\langle\tau_0(2)^{n_e}
\tau_1(2)^{n_v}\big\rangle^1_\Delta$.
\end{remark}

\begin{proposition}\label{cor1new}
Given a balanced, nondegenerate multiset
$\Delta\subset\Z^2\setminus2\{0\}$, 
integers $n_v>0$, $n_e\ge0$ such that
$2n_v+n_e=|\Delta|$, $n=n_v+n_e$, and a generic point
$\bx\in\R^{2n}$, the following holds:
\begin{enumerate}
\item[(1)] if $(\overline\Gamma,\bp,h)\in{\mathcal M}'_{1,(n_v,n_e)}
(\R^2,\Delta)$, $h(\bp)=\bx$, then $\GS_{-1}(\overline\Gamma,\bp,h)\ne0$ if and only
$(\overline\Gamma,\bp,h)$ is an elliptic broccoli curve; 
\item[(2)] if $(\overline\Gamma,\bp,h)\in{\mathcal M}''_{1,(n_v,n_e)}
(\R^2,\Delta)$, $h(\bp)=\bx$, and $\pi(\overline\Gamma,\bp,h)=
(\overline\Gamma',\bp,p,h')\in{\mathcal M}''_{0,(n_v-1,n_e+1)}
(\R^2,\Delta)$ then $\GS_{-1}(\overline\Gamma',\bp,p,h')\ne0$ if and only if
$(\overline\Gamma,\bp,h)$ is an elliptic broccoli curve. 
\end{enumerate}
\end{proposition}

{\bf Proof.}
Following the argument in the proof of Proposition \ref{p5} and using the computations of
Lemma \ref{l8}, one immediately obtains that
$\GS_{-1}(\overline\Gamma,\bp,h)$, resp. $\GS_{-1}(\overline\Gamma',\bp,p,h')$ vanishes
whenever $(\Gamma,\bp,h)$, resp. $(\Gamma',\bp,p,h')$ is not an elliptic broccoli curve.
The same reasoning yields that the weight $\GS_{-1}$ does not vanish for elliptic broccoli curves.
\proofend

\subsection{Lattice path algorithm}\label{sec-lp}
To efficiently compute the refined elliptic broccoli invariant, we provide here a
suitable version of the lattice path algorithm. It is a simplification of the algorithm from
\cite[Section 9]{MR}, in which we allow only trivalent vertices, but also consider elliptic curves.
For the reader's convenience, we provide here all details.

\smallskip
{\it (1) Initial data and general procedure.} Let us be given a balanced, nondegenerate multiset
$\Delta\subset\Z^2\setminus\{0\}$ and integers $n_v>0$, $n_e\ge0$ such that
$2n_v+n_e=|\Delta|$, $n=n_v+n_e$. Pick a vector $\overline a\in\R^2\setminus\{0\}$, which is not parallel or orthogonal
to any vector $u_1-u_2$, where $u_1,u_2\in P(\Delta)\cap\Z^2$, $u_1\ne u_2$, and consider a straight line
$L_{\overline a}$ through the origin, directed by $\overline a$.
Introduce a configuration of points $\bx=(x_1,...,x_n)\subset L_{\overline a}$ such that
$$x_i=M_i\overline a,\ i=1,...,n,\quad 0<M_1\ll M_2\ll...\ll M_n\ .$$

The linear functional $\varphi_{\overline a}(x)=\langle x,\overline a\rangle:\R^2\to\R$ defines a linear order
on the points $u\in P(\Delta)\cap \Z^2$, that is, $u\prec u'$ if $\varphi_{\overline a}(u)<\varphi_{\overline a}(u')$.
Denote by $u_{\min},u_{\max}\in P(\Delta)$ the extremal points. Any $\varphi_{\overline a}\;$-monotone sequence of points
$u_0=u_{\min}\prec u_1\prec...\prec u_r=u_{\max}\subset P(\Delta)\cap\Z^2$ is called a lattice path in $P(\Delta)$ (of length $r$).
Denote by ${\mathcal L}_n(\Delta,\overline a)$ the set of lattice paths in $P(\Delta)$ of length $n$.

The algorithm starts with a lattice path $G\in{\mathcal L}_n(\Delta,\overline a)$. Then we construct an enhanced lattice path
$\widehat G$, and inductively extend it (if possible) to
a certain subdivision of the polygon $P(\Delta)$. Such a subdivision uniquely determines a tropical
curve passing through $\bx$. We verify whether it is irreducible and has degree $\Delta$. Then it must be elliptic.
If the reconstructed tropical curve $(\overline\Gamma,\bp,h)$ is irreducible and has no collinear cycle, we assign to it a
refined weight $\GS_y(\overline\Gamma,\bp,h)$. If $(\overline\Gamma,\bp,h)$ contains a collinear cycle, we take
the rational curve
$\pi(\overline\Gamma,\bp,h)=(\overline\Gamma',\bp,p,h')$ and assign to it the refined weight
$\GS_y(\overline\Gamma',\bp,p,h')$. In both situations, the refined weights can be expressed in terms of
the pair (enhanced lattice path, subdivision).

Notice that if either a lattice path cannot be equipped with an enhancement, or a current subdivision does not cover $P(\Delta)$
and cannot be extended anymore, or the subdivision covers $P(\Delta)$ but does not define an irreducible tropical curve
of degree $\Delta$, we skip these outcomes.

Given a subdivision $S$ as a set of convex lattice polygons and segments, we let $|S|=\bigcup_{\sigma\in S}\sigma$.

\smallskip
{\it (2) Enhancement of a lattice path.} Let $G=(u_0,...,u_n)\in{\mathcal L}_n(\Delta,\overline a)$.
Denote by $S_0(G)$ the set of the segments $[u_{i-1},u_i]$, $i=1,...,n$.

Let $\partial P(\Delta)_+$ and $\partial P(\Delta)_-$ be the two components of of $\partial P(\Delta)\setminus\{u_{\min},u_{\max}\}$.
Each component of $P(\Delta)\setminus|S_0(G)|$ has a nonempty intersection either with
$\partial P(\Delta)_+$, or with $\partial P(\Delta)_-$. Denote by $P(\Delta)_+$ (resp., $P(\Delta)_-$) the 
the union of $\partial P(\Delta)_+$ (resp., $\partial P(\Delta)_-$) with the components of
$P(\Delta)\setminus|S_0(G)|$ intersecting $\partial P(\Delta)_+$ (resp., $\partial P(\Delta)_-$).

An enhancement $\widehat G$ of $G$ is a pair:
\begin{itemize}\item a lattice path of length $n+n_v$ extending $G$ with extra points $u'_i$, $i=1,...,n_v\in P(\Delta)\cap\Z^2$
such that $u_{i-1}\prec u'_i\prec u_i$ for all $i=1,...,n_v$;
\item a sequence of signs $\eps_i=\pm1$, $i=1,...,n_v$: if $u'_i\in P(\Delta)_+$ (resp., $u'_i\in P(\Delta)_-$) we set $\eps_i=1$
(resp., $\eps_i=-1$), if $u'_i\in[u_{i-1},u_i]$ we choose either
$\eps_i=1$, or $\eps_i=-1$.\end{itemize}
Denote by ${\mathcal L}_{n_v,n_e}(\Delta,\overline q)\subset
{\mathcal L}_n(\Delta,\overline a)$ the set of those lattice path which admit an enhancement $\widehat G$.
From now on we suppose that $G\in{\mathcal L}_{n_v,n_e}(\Delta,\overline a)$.

\smallskip
{\it (3) Initial subdivision.}
Define $S_1(\widehat G)$ to be the set of polygons $\conv\big\{u_{i-1},u'_i,u_i\big\}$, $i=1,...,n_v$, and the segments
$[u_{n_v+j-1},u_{n_v+j}]$, $j=1,...,n_e$. and let $\partial S_1(\widehat G)_+$ (resp.,
$\partial S_1(\widehat G)_-$) to be the lattice path consisting of the points $U_0,...,u_n$ and the points
$u'_i$, $1\le i\le n_v$, such that $\eps_i=1$ (resp., $\eps_i=-1$).

\smallskip
{\it (4) Step of the algorithm.} Let us be given $k\ge1$, a subdivision $S_k(\widehat G)$ and two
lattice paths $\partial S_k(\widehat G)_+$, $\partial S_k(\widehat G)_-$
such that $S_k(\widehat G)_\eps\subset\partial |S_k(\widehat G)|\cap\overline{P(\Delta)}_\eps$\;, $\eps=\pm1$, and
$\partial|S_k(\widehat G)|$ is the union of the broken lines induced by $\partial S_k(\widehat G)_{\pm}$.

If $\partial S_k(\widehat G)_+=\big\{v_0=u_{\min}\prec v_1...\prec v_r=u_{\max}\big\}$, $r\ge2$, we look for the minimal $i=1,...,r-1$
such that the points $v_{i-1},v_i,v_{i+1}$ are not collinear, and the triangle $T_i=\conv\big\{v_{i-1},v_i,v_{i+1}\big\}$ is not
contained in $|S_k(\widehat G)|$. If such $i$ does exist, we
\begin{itemize} \item either set $S_{k+1}(\widehat G)=S_k(\widehat G)\cup\{T\}$ and $\partial S_{k+1}(\widehat G)_+=\partial S_k(\widehat G)_+\setminus\{v_i\}$;
\item or, in case the parallelogram $\Pi_i=\conv\big\{v_{i-1},v_i,v_{i+1},v'_i\big\}$ lies inside $P(\Delta)$, set
$S_{k+1}(\widehat G)=S_k(\widehat G)\cup\{\Pi_i\}$ and $\partial S_{k+1}(\widehat G)=\big\{v_o\prec...\prec v_{i-1}\prec v'_i\prec v_{i+1}\prec...\prec v_r\big\}$.
\end{itemize} In both cases $\partial S_{k+1}(\widehat G)_-=\partial
S_k(\widehat G)_-$.

If we cannot perform the above step, we do the same exchanging all the signs.

\smallskip
{\it (5) Restoring a tropical curve.} Since the area of $|S_k(\widehat G)|$ strictly grows, the algorithm is finite.
Let $S_{fin}(\widehat G)$, $\partial S_{fin}(\widehat G)_+$, $\partial S_{fin}(\widehat G)_-$ be the outcome. We call the outcome
admissible if
\begin{itemize}\item $|S_{fin}(\widehat G)|=P(\Delta)$; in this case the broken lines induces by
$\partial S_{fin}(\widehat G)_{\pm}$ cover $\partial P(\Delta)$;
\item $\partial S_{fin}(\widehat G)_{\pm}$ match the degree $\Delta$; this means that, orienting the segments that join consecutive
points of $\partial S_{fin}(\widehat G)_{\pm}$ counter-clockwise and rotating all of them by $\pi/2$ clockwise, we
obtain the multiset of vectors $\Delta$>
\end{itemize}
Non-admissible outcomes are skipped in the count.

To an admissible outcome, we assign a tropical curve $(\overline\Gamma,\bp,h)$ of degree $\Delta$ with an $n$-tuple of marked points $\bp$
such that $h(\bp)=\bx$. In fact, we follow the above algorithm in its dual form and, moreover, in parallel, construct an
orientation of the components of $\overline\Gamma\setminus\bp$. Denote by $\overline b\in\R^2\setminus\{0\}$ the unit vector
orthogonal to $\overline a$ and oriented from $P(\Delta)_-$ towards $P(\Delta)_+$.

In the construction we shall use auxiliary objects, plane tropical precusrves, which are the following objects: Given
a marked plane tropical curve
$(\overline\Gamma,\bp,h)$ and any open bounded subset $\Gamma'\subset\Gamma$ containing $\bp$, we say that the triple
$(\Gamma',\bp,h\big|_{\Gamma'})$ is a marked plane tropical precurve (associated with $(\overline\Gamma,\bp,h)$).

Take an $n$-tuple of points $\bp=(p_1,...,p_n)$ and consider the tropical precurve
$(\Gamma_1,\bp,h_1)$ with $\Gamma_1$ being the union of graph germs $(\Gamma_1,p_i)$, $i=1,...,n$, which are
trivalent for $i=1,...,n_v$ and are bivalent for $i=n_v+1,...,n$ and with the map
$h_1:(\Gamma_1,\bp)\to\R^2$ determined by the conditions
\begin{itemize}\item $p_i\in\bp\mapsto h_1(p_i)=x_i\in\bx$, $i=1,...,n$;
\item if $n_v<i\le n$, the germ $(\Gamma_1,p_i)$ is mapped onto the germ of a straight line through $x_i$ orthogonal to
the segment $[u_{i-1},u_i]$;
\item if $1\le i\le n_v$, then the edges of $(\Gamma_1,p_i)$ are mapped to
the three rays rooted at $x_i$ and directed by the vectors orthogonal to the segments
$[u_{i-1},u_i]$, $[u_{i-1},u'_i]$, $[u'_i,u_i]$ and such that their scalar products with $\eps_i\overline b$ are respectively
negative, positive, and positive;
\item the differential $Dh_1$ along each edge the Euclidean length of the corresponding orthogonal segment, mentioned above.
\end{itemize} We then inductively proceed along the lattice path algorithm
extending the current tropical precurve.
Namely, having a tropical precurve $(\Gamma_k,\bp,h_k)$, we follow step (4) of the lattice path algorithm:
\begin{itemize}
\item
extend the edges $e,e'$ of $\Gamma_k$ corresponding to the
segments $[v_{i-1},v_i]$, $[v_i,v_{i+1}]$ induced by $\partial S_k(\widehat G)_+$ (or $\partial S_k(\widehat G)_-$), until
$h(e)$ and $h(e')$ hit each other;
\item if there is no third edge of $\Gamma_k$, which is mapped by $h$ to the same ray as
$E$ or $E'$, then either form a new
trivalent vertex of $\Gamma_{k+1}$ mapping the third edge
orthogonal to the segment $[v_{i-1},v_{i+1}]$, or slightly extend further the considered edges of the current
tropical precurve without
gluing them, while their $h$-images transversally intersect (cf. \cite[Figure 2.17
in Section 2.5.6]{IMS});
\item if there is an edge $e''$ of $\Gamma_k$, which is mapped to the same ray as, say, $e'$
(this can only happen when $e',e''$ are incident to the same collinear marked vertex), then
either we extend all three edges without gluing so that $h(e')=h(e'')$ intersects transversally
$h(e)$, or one of $e',e''$ joins with $e$ forming a trivalent vertex with the germ of a new edge
and the remaining edge
extends further, or all $e,e',e''$ join together forming a four-valent vertex with
the germ of a new edge. 
\end{itemize} Notice that the very last situation one obtains a precurve with a collinear cycle.

The final tropical precurve turns into a plane tropical curve as we
extend all open edges to unbounded rays and compactify them with univalent vertices.
We claim that, if it is irreducible then it is elliptic. Indeed, it is easy to see from the construction that each
component of $\overline\Gamma\setminus\bp$ is a tree containing exactly one univalent vertex of $\overline\Gamma$. Hence
$1-g(\overline\Gamma)=\chi(\overline\Gamma)=|\Delta|-2n_v-n_e=0$.

\begin{figure}
\setlength{\unitlength}{1cm}
\begin{picture}(14,5)(0,0)
\thinlines

\put(1,2.25){\line(0,1){1}}\put(2.5,1.5){\line(0,1){1}}\put(1,2.25){\line(2,-1){1.5}}
\put(1,3.25){\line(1,1){0.5}}\put(2.5,2.5){\line(1,1){0.5}}\put(3,3){\line(-2,1){1.5}}
\put(3,3){\line(-1,1){1}}\put(1.5,3.75){\line(2,1){0.5}}
\put(2,4){\line(0,-1){0.5}}\put(2,3.5){\line(-1,-1){0.5}}
\put(1,2.25){\line(-1,1){0.5}}\put(2.5,1.5){\line(1,0){0.5}}

\dashline{0.2}(0.5,2.75)(0,3.25)\dashline{0.2}(3,1.5)(3.5,1.5)

\put(4.5,2.25){\line(0,1){1}}\put(6,1.5){\line(0,1){1}}\put(4.5,2.25){\line(2,-1){1.5}}
\put(4.5,3.25){\line(1,1){0.5}}\put(6,2.5){\line(1,1){0.5}}\put(6.5,3){\line(-2,1){1.5}}
\put(6.5,3){\line(-1,1){1}}\put(5,3.75){\line(2,1){0.5}}
\put(5.5,4){\line(0,-1){0.5}}\put(5.5,3.5){\line(-1,-1){0.5}}
\put(4.5,2.25){\line(1,-1){1}}\put(5.5,1.25){\line(2,1){0.5}}

\put(8,2.25){\line(0,1){1}}\put(9.5,1.5){\line(0,1){1}}\put(8,2.25){\line(2,-1){1.5}}
\put(8,3.25){\line(1,1){0.5}}\put(9.5,2.5){\line(1,1){0.5}}\put(10,3){\line(-2,1){1.5}}
\put(10,3){\line(-1,1){1}}\put(8.5,3.75){\line(2,1){0.5}}
\put(8,2.25){\line(-1,1){0.5}}\put(9.5,1.5){\line(1,0){0.5}}

\dashline{0.2}(7.5,2.75)(7,3.25)\dashline{0.2}(10,1.5)(10.5,1.5)

\put(11.5,2.25){\line(0,1){1}}\put(13,1.5){\line(0,1){1}}\put(11.5,2.25){\line(2,-1){1.5}}
\put(11.5,3.25){\line(1,1){0.5}}\put(13,2.5){\line(1,1){0.5}}\put(13.5,3){\line(-2,1){1.5}}
\put(13.5,3){\line(-1,1){1}}\put(12,3.75){\line(2,1){0.5}}
\put(11.5,2.25){\line(1,-1){1}}\put(12.5,1.25){\line(2,1){0.5}}

\thicklines

\put(1,3.25){\line(2,-1){1.5}}\put(4.5,3.25){\line(2,-1){1.5}}
\put(8,3.25){\line(2,-1){1.5}}\put(11.5,3.25){\line(2,-1){1.5}}

\put(0.9,3.13){$\bullet$}\put(1.4,2.9){$\bullet$}\put(2.4,2.4){$\bullet$}
\put(4.4,3.13){$\bullet$}\put(4.9,2.9){$\bullet$}\put(5.9,2.4){$\bullet$}
\put(7.9,3.13){$\bullet$}\put(9.4,2.4){$\bullet$}
\put(11.4,3.13){$\bullet$}\put(12.9,2.4){$\bullet$}

\put(0.3,3.4){$u_{i-1}$}\put(3.8,3.4){$u_{i-1}$}\put(7.3,3.4){$u_{i-1}$}\put(10.8,3.4){$u_{i-1}$}
\put(1.3,2.5){$u'_i$}\put(4.8,2.5){$u'_i$}
\put(2.7,2.3){$u_i$}\put(6.2,2.3){$u_i$}\put(9.7,2.3){$u_i$}\put(13.2,2.3){$u_i$}
\put(2.1,1){$\partial P(\Delta)$}\put(9.1,1){$\partial P(\Delta)$}

\put(1.8,0){(a)}\put(5.3,0){(b)}\put(8.8,0){(c)}\put(12.3,0){(d)}

\end{picture}
\caption{Subdivisions related to collinear cycles}\label{fsub}
\end{figure}
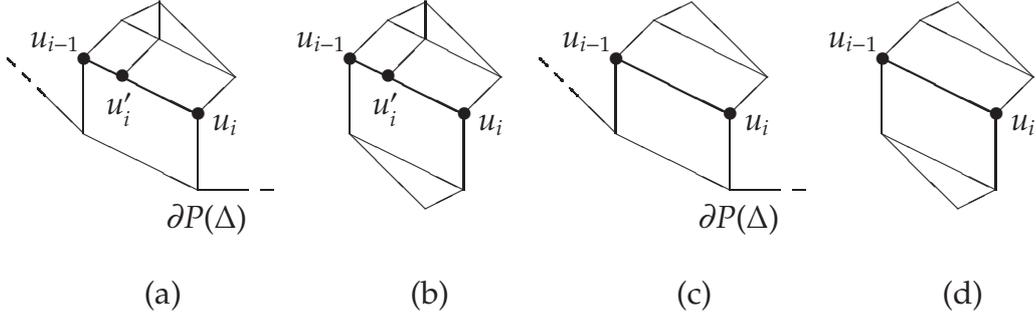

\smallskip
{\it (6) Refined broccoli weights.}
Given an admissible subdivision $S_{fin}$, corresponding to an irreducible tropical curve without
collinear cycle, we define its refined broccoli weight as the product of
\begin{itemize}\item $[\mu_i]^+_y$, where $\mu_i$ is the lattice area of the (possibly degenerate)
triangle $\conv\{u_{i-1},u'_i,u_i\}$, over all $i=1,...,n_v$,
\item and $[\mu]^-_y$, where $\mu$ runs over the lattice areas of all other triangles in $S_{fin}$,\end{itemize}
which is divided by $2^r$, $r$ being the number of segments among $\conv\{u_{i-1},u'_i,u_i\}$,
$i=1,...,n_v$.
If $S_{fin}$ corresponds to an irreducible tropical curve with a collinear cycle, it contains a special fragment like shown in
Figure \ref{fsub}(a,b) (where the number of the incline of parallelograms may vary).
We transform this fragment as shown in Figure \ref{fsub}(c,d), respectively, and then assign a refined weight to
the newly obtained subdivision $S'_{fin}$ to be the product of
\begin{itemize} \item either $\Psi^{(1)}(m,\nu)$, or $\Psi^{(2)}(m,\nu_1,\nu_2)$ as $S'_{fin}$ contains the fragment as in
Figure \ref{fsub}(c,d), respectively,
\item $[\mu_j]^+_y$, $j=1,...,n_v$, $j\ne i$,
\item $[\mu]^-_y$, where $\mu$ runs over all triangles of $S'_{fin}$ different from $\conv\{u_{j-1},u'_j,u_j\}$,
$j=1,...,n_v$, $j\ne i$, and lying outside the special fragment,
\item $[\mu]^*_y$, where $\mu$ runs over the lattice length of all the vectors in $\Delta$,
\end{itemize}
which is divided by $2^r$, $r$ being the number of segments among
$\conv\{u_{j-1},u'_j,u_j\}$,
$j=1,...,n_v$, $j\ne i$.

\begin{example}
The case of $\Delta=\{2\times(-1,0),3\times(0,1),2\times(1,0),3\times(0-1)\}$ is one of the simplest, when
the refined elliptic broccoli invariant is not constant. Here, $n_e+2n_v=|\Delta|=10$, and the lattice path algorithm
gives $$\GS_y(\Delta,1,(n_e,n_v))=2y+16-2n_v+2y^{-1},
\quad n_v=0,...,5\ .$$
We also notice that the lattice path algorithm with the vector $\overline a=(1,\eps)$, $0<\eps\ll1$,
counts only tropical curves in ${\mathcal M}'_{1,(n_v,n_e)}(\R^2,\Delta)$, while the procedure with the vector
$\overline a=(\eps,1)$ counts also tropical curves in ${\mathcal M}''_{1,(n_v,n_e)}(\R^2,\Delta)$.
\end{example}

{\ncsc Fachbereich Mathematik (AD)  \\[-21pt]

Universit\"at Hamburg\\[-21pt]

Bundesstrasse 55, 20146 Hamburg, Germany}\\[-21pt]

{\it E-mail address}: {\ntt franziska.schroeter@uni-hamburg.de}

\vskip10pt

{\ncsc School of Mathematical Sciences \\[-21pt]

Raymond and Beverly Sackler Faculty of Exact Sciences\\[-21pt]

Tel Aviv University \\[-21pt]

Ramat Aviv, 69978 Tel Aviv, Israel} \\[-21pt]

{\it E-mail address}: {\ntt shustin@post.tau.ac.il}

\end{document}